\documentclass[a4 paper, 11pt]{smfart}

\usepackage{smfthm}
\usepackage{smfenum}
\usepackage{amssymb}
\usepackage{amscd}

\usepackage[english,frenchb]{babel}  
\usepackage{textcomp}  
\usepackage[all]{xy}  

\usepackage{amssymb,amsmath,latexsym,epsfig,euscript}

\def\Conv {{ \mbox{\rm Conv}}}

\def\Vol {{\mbox{\rm Vol}}}

\def\ord {\mbox{\scriptsize ord}}
\def\Ord {\mbox{\rm ord}}

\def\PSL {{\rm \mbox{PSL}}}

\def\Res {\mbox{\rm R{\'e}s}}
\def\res {\mbox{\rm r{\'e}s}}
\def\elim {\mbox{\rm {\'e}lim}}

\def\cc {{ \mbox{\rm \scriptsize c}}}
\def\ccc {{ \mbox{\rm \tiny c}}}

\def\Ann {\mbox{\rm Ann}}

\def\Graph {\mbox{\rm Graphe}}

\def\irr {\mbox{\rm \scriptsize Irr}\,}
\def\FS {\mbox{\rm \scriptsize FS}\,}

\def\laurent {{\mbox{\rm \scriptsize Lau}}}
\def\daph {{\mbox{\rm \scriptsize DP}}}

\def\geo {{\mbox{\rm \scriptsize g{\'e}om}}}
\def\norm {{\mbox{\rm \scriptsize norm}}}
\def\ari {{\mbox{\rm \scriptsize arith}}}

\def\Spec   {\mbox{\rm Spec}\,}
\def\MV {\mbox{\rm MV}}
\def\MI {\mbox{\rm MI}}

\def\Card {\mbox{\rm Card}}
\def\card {{\mbox{\rm \scriptsize Card}}}
\def\Ch {{\mathcal{C}{\it h}}}

\def\Expo {\mbox{\rm Expo}}
\def\expo {\mbox{\rm \scriptsize Expo}}
\def\Ker {\mbox{\rm Ker}}

\def\Spec {\mbox{\rm Spec}}

\def\Stab  {\mbox{\rm Stab}}

\def\OK {{\cO_K}}
\def\Hnorm {\mathcal{H}_\norm}
\def\Harith {{\mathcal{H}_\ari}}
\def\hnorm {{\widehat{h}}}

\def\geom {{\mathcal A}}
\def\arith {{{\mathcal A}, \alpha}}
\def\cad {c.-{\`a}-d. }

\def\ov#1{{\overline{#1}}}
\def\un#1{{\underline{#1}}}
\def\wh#1{{\widehat{#1}}}

\def \C {\mathbb{C}}

\def \K {\mathbb{K}}
\def \N {\mathbb{N}}
\def \P {\mathbb{P}}
\def \Q {\mathbb{Q}}
\def \R {\mathbb{R}}
\def \T {\mathbb{T}}
\def \Z {\mathbb{Z}}

\def \cA {\mathcal{A}}
\def \cB {\mathcal{B}}
\def \cC {\mathcal{C}}

\def \cH {\mathcal{H}}
\def \cI {\mathcal{I}}

\def \cL {\mathcal{L}}
\def \cM {\mathcal{M}}

\def \cO {\mathcal{O}}

\def\gp {{\mathfrak{p}}}

\def\Qbar {{\overline{\Q}}}

\title[Hauteur normalis{\'e}e des vari{\'e}t{\'e}s toriques projectives]{Hauteur 
normalis{\'e}e des vari{\'e}t{\'e}s toriques projectives}

\author{Patrice Philippon} 

\address{Institut de Math{\'e}matiques - U.M.R. 7586 du CNRS,
Projet G{\'e}om{\'e}trie et Dynamique,
Case~7012, 
2 place Jussieu, 
75251 Paris Cedex 05,
France.}

\email{pph@math.jussieu.fr}

\urladdr{http://www.math.jussieu.fr/\~{}pph/}

\author{Mart{\'\i}n Sombra}

\address{ Universit{\'e} de Lyon 1, 
Laboratoire de Math{\'e}matiques 
Appliqu{\'e}es de Lyon, 
21 avenue~Claude~Bernard,
69622 Villeurbanne Cedex, 
France.}

\email{sombra@maply.univ-lyon1.fr}

\urladdr{http://maply.univ-lyon1.fr/\~{}sombra/} 



\thanks{P. Philippon a {\'e}t{\'e} partiellement financ{\'e} par une allocation de recherche de la {\em Fondation Alexander von Humboldt} pendant la r{\'e}alisation de ce travail.} 

\thanks{\newline M. Sombra a {\'e}t{\'e} partiellement financ{\'e} par une bourse post-doctorale Marie Curie du programme europ{\'e}en {\em Improving  Human Research Potential and the Socio-economic Knowledge Base}, contrat n\textordmasculine \ HPMFCT-2000-00709.} 

\date{Version du 18 Mars 2004.}

\begin{document}

\font\petitbf=cmbx9

\newdimen\fontdim
\def\og{\fontdim=\fontdimen6\font \font\fant=lasy10 at \fontdim
\ouvreg\fant\everypar={\ouvreg\fant}}
\def\ogg{\fontdim=\fontdimen6\font \font\fant=lasyb10 at \fontdim
\ouvreg\fant\everypar={\ouvreg\fant}}
\def\ouvreg#1{\leavevmode\hbox{#1\unskip\kern+0.20em(\kern-0.20em(\kern+0.20em}%
\ignorespaces\nobreak}
\def\fg{\fermeg\fant\everypar={}} \let\fgg=\fg
\def\fermeg#1{\relax \ifhmode \unskip\kern+0.20em\else \leavevmode \fi%
\hbox{#1)\kern-0.20em)\kern+0.20em\ignorespaces}}

\setlength{\baselineskip}{14pt}

%
%


\begin{abstract}
\setlength{\baselineskip}{12pt}
Nous pr{\'e}sentons une expression explicite pour la hauteur normalis{\'e}e
d'une vari{\'e}t{\'e} torique projective. Cette expression se d{\'e}compose comme
somme de contributions locales, chaque terme {\'e}tant l'int{\'e}grale d'une 
certaine fonction concave et affine par morceaux. 
Plus g{\'e}n{\'e}ralement, nous obtenons une expression explicite pour la multihauteur normalis{\'e}e d'un tore relative {\`a} plusieurs plongements monomiaux. 

L'ensemble de fonctions introduit se comporte comme un analogue arithm{\'e}tique 
du polytope classiquement associ{\'e} {\`a} l'action du tore. En plus des
formules pour les hauteur et multihauteurs, nous montrons que cet
objet se comporte de mani{\`e}re naturelle par rapport {\`a} plusieurs
constructions standard: 
d{\'e}composition en orbites, formation de joints, produits de Segre et
plongements de Veronese.

La d{\'e}monstration suit une d{\'e}marche indirecte: {\`a} la place de la
d{\'e}finition de la hauteur normalis{\'e}e, on s'appuie sur le calcul d'une 
fonction de Hilbert arithm{\'e}tique appropri{\'e}e. 
\end{abstract}


\begin{altabstract}
\setlength{\baselineskip}{12pt}
{\bf Normalized height of projective toric varieties.} 
We present an explicit expression for the normalized height of a projective toric variety. This expression decomposes as a sum of local contributions, 
each term being the integral of a certain function, 
concave and piecewise linear-affine. 
More generally, we obtain an explicit expression for the normalized
multiheight of a torus with respect to several monomial embeddings.

The set of functions introduced behaves as an arithmetic analog of the 
polytope classically associated with the torus action. 
Besides the formul\ae \ for the height and multiheight, we show that
this object behaves in a natural way with respect to several standard
constructions: 
decomposition into orbits, joins, Segre products and Veronese embeddings.

The proof follows an indirect way~: instead of the definition of the 
normalized height, we rely on the computation of an appropriate
arithmetic Hilbert function.

\end{altabstract}

\maketitle

\vspace*{-12mm}


\typeout{Introduction et resultats}

\section*{Introduction et r{\'e}sultats}

La {hauteur} d'une sous-vari{\'e}t{\'e} $X \subset \P^N(\Qbar)$ est une mesure 
de la complexit{\'e} binaire d'une repr{\'e}sentation de $X$, par exemple {\it via} 
sa forme de Chow. Cet invariant num{\'e}rique est donc relevant en g{\'e}om{\'e}trie alg{\'e}brique et en alg{\`e}bre commutative effectives, o{\`u} il joue un r{\^o}le important, notamment dans le cadre du Nullstellensatz 
effectif et de la r{\'e}solution de syst{\`e}mes d'{\'e}quations polynomiales, 
{\it voir} par exemple~\cite{KPS01}, \cite{GHP97}.  

La notion de hauteur de sous-vari{\'e}t{\'e}s g{\'e}n{\'e}ralise celle
appliqu{\'e}e aux points
par Siegel, North\-cott et Weil entre autres, pour {\'e}tudier
des questions d'approximation diophantienne. 
En dimension sup{\'e}rieure, cette notion est un analogue arithm{\'e}tique du
degr{\'e}, ou plus pr{\'e}cisement d'un bi-degr{\'e} dans $\P^1\times\P^N$. 
Par exemple, elle v{\'e}rifie des analogues du th{\'e}or{\`e}me de B{\'e}zout et de la 
formule de Hilbert-Samuel~\cite{BGS94}, \cite{Phi95}. 

\smallskip 

Il existe en fait diverses notions de hauteur, certaines 
pr{\'e}sentant des propri{\'e}t{\'e}s d'invariance par rapport {\`a} des 
op{\'e}rations g{\'e}om{\'e}triques significatives. 
En particulier, l'espace projectif peut {\^e}tre vu comme 
une compactification {\'e}quivariante du tore $(\P^N)^\circ:= 
\P^N \setminus \{ (x_0: \cdots : x_N) \, : \ x_0 \cdots x_N = 0\}$,
la structure de groupe de $(\P^N)^\circ$ permettant de d{\'e}finir une 
notion de  hauteur pour les sous-vari{\'e}t{\'e}s de $\P^N$ plus canonique que
les autres, appel{\'e}e {\it hauteur normalis{\'e}e}. Cette notion joue un
r{\^o}le central dans l'approximation diophantienne sur les tores, et tout
particuli{\`e}rement dans les probl{\`e}mes de Bogomolov et de Lehmer
g{\'e}n{\'e}ralis{\'e}s, {\it voir}~\cite{DP99}, \cite{AD02} et leurs r{\'e}f{\'e}rences,
{\it voir} {\'e}galement~\cite{Dav02} pour un aper{\c c}u historique.

Suivant~\cite{DP99}, la hauteur normalis{\'e}e peut se d{\'e}finir par un proc{\'e}d{\'e} \og {\`a} la Tate\fg. De fa{\c c}on pr{\'e}cise, pour $k \in \N$ 
on pose $ [k]: \P^N \to \P^N$, $(x_0: \cdots : x_N) \mapsto (x_0^k: \cdots : x_N^k)$ l'application puissance $k$-i{\`e}me; restreinte au tore $(\P^N)^\circ$ 
c'est l'application de multiplication par $k$. La {hauteur normalis{\'e}e}  
d'une sous-vari{\'e}t{\'e}  (r{\'e}duite et irr{\'e}ductible) $X \subset \P^N$ 
est par d{\'e}finition
$$
\wh{h}(X) := \deg(X) \cdot \lim_{k\to \infty} \, 
\frac{h({[k] \, X})}{k\, \deg([k]\, X)} \ \in \R_+
\enspace, 
$$ 
o{\`u} $\deg$ et $h$ d{\'e}signent  le degr{\'e} et la  hauteur projective, {\it
  voir}~\cite[\S \ 2]{DP99} ou le paragraphe~\ref{hauteur normalisee}
ci-dessous. Cette hauteur peut aussi se d{\'e}finir {\em via} la th{\'e}orie
d'Arakelov, comme la hauteur de 
la cl{\^o}ture de Zariski de $X$
dans $\P^N_{\cO_K}$ (l'espace projectif de dimension $N$
sur $\Spec(\cO_K)$) relative  au
fibr{\'e} en droites universel  $\cO(1)$ muni d'une m{\'e}trique 
hermitienne canonique~\cite{Zha95b}, \cite{Mai00}.

\smallskip

Contrairement au degr{\'e}, il y a peu de cas o{\`u} l'on sait calculer explicitement la hauteur d'une vari{\'e}t{\'e}, qui se r{\'e}v{\`e}le alors {\^e}tre toujours un nombre remarquable. Pour les hypersurfaces, la hauteur normalis{\'e}e s'{\'e}crit comme somme de contributions locales aux places 
archim{\'e}diennes, dont chacune est la  mesure de Mahler locale
d'une {\'e}quation de d{\'e}finition de $X$. Les valeurs des mesures de Mahler
sont tr{\`e}s {\'e}tudi{\'e}es, en particulier {\`a} cause de leurs connexions avec
les valeurs sp{\'e}ciales des fonctions $L$, {\it voir}~\cite{Boy98} et
ses r{\'e}f{\'e}rences. Une autre situation o{\`u} la hauteur normalis{\'e}e est
connue est celle des {\it vari{\'e}t{\'e}s de torsion}, \cad les sous-vari{\'e}t{\'e}s
telles que $X^\circ := X \cap (\P^N)^\circ$ soit le
translat{\'e} d'un sous-groupe par un point de torsion; dans ce cas la
hauteur 
normalis{\'e}e est nulle, et c'est le seul cas o{\`u} cela se produit. 

Pour la hauteur dite {\it projective}, le calcul d'exemples a {\'e}t{\'e}
surtout d{\'e}velopp{\'e} pour les espaces homog{\`e}nes sous un groupe de
Chevalley~\cite{KK99}, et de fa{\c c}on plus explicite pour les
$\mbox{SL}_N$-grasmaniennes et d'autres vari{\'e}t{\'e}s du m{\^e}me
type~\cite{Tam00}. 
En dehors de ces cas, il y a quelques r{\'e}sultats pour certaines hypersurfaces, {\it voir} par exemple~\cite{BY98}. 

\smallskip  

Les vari{\'e}t{\'e}s toriques sont l'une des principales sources d'exemples en 
g{\'e}om{\'e}trie alg{\'e}brique. Dans le pr{\'e}sent texte, nous {\'e}tudierons 
la hauteur normalis{\'e}e des {\it vari{\'e}t{\'e}s toriques projectives},
d{\'e}finies comme les sous-vari{\'e}t{\'e}s de $\P^N$ stables par rapport {\`a} une
action diagonale d'un tore dont une orbite est dense. {\`A} une telle
vari{\'e}t{\'e} $X_\arith$ (o{\`u} plut{\^o}t {\`a} l'action du tore et au choix d'un
point dans l'orbite principale) nous  
associons un ensemble de fonctions concaves et  affines par morceaux
$$
\Theta_\arith:= \Big( 
\vartheta_{\cA,\tau_{\alpha v}} \, : \ Q_\geom \to \R\Big)_{v \in M_K}
$$ 
index{\'e} par les places d'un corps de nombres appropri{\'e}, chaque fonction
$\vartheta_{\cA,\tau_{\alpha v}}$ {\'e}tant d{\'e}finie sur le polytope 
$Q_\geom \subset \R^n$ classiquement associ{\'e} {\`a} l'action du tore, 
au moyen du vecteur poids $\tau_{\alpha v}\in\R^{N+1}$, {\it voir} ci-dessous pour les d{\'e}tails.

Notre r{\'e}sultat principal (th{\'e}or{\`e}me~\ref{thm1} ci-dessous) affirme 
que la hauteur $\hnorm(X_\arith)$ s'exprime comme somme de contributions locales, dont chaque terme est l'int{\'e}grale de la fonction 
correspondante. 
De plus, cette expression s'{\'e}tend en une expression 
explicite pour la multihauteur normalis{\'e}e d'un  tore 
relative {\`a} plusieurs plongements monomiaux. 

L'ensemble de fonctions $\Theta_\arith$ se comporte comme un analogue
arithm{\'e}tique du polytope $Q_\geom$. En plus des formules 
pour la hauteur et multihauteur, nous montrons que cet objet se
comporte de mani{\`e}re naturelle par rapport {\`a} d'autres constructions
standard: 
d{\'e}composition en orbites, formation de joints, 
produits de Segre et plongements de Veronese. 

\bigskip

Soit $\T^n:= (\Qbar^\times)^n$ le tore alg{\'e}brique et  $\P^N$ l'espace projectif sur $\Qbar$, de dimensions $n$ et $N$ respectivement. Soit $\cA= ( a_0,  \dots, a_N) \in (\Z^n)^{N+1}$ une suite de $N+1$ vecteurs de $\Z^n$, on consid{\`e}re alors l'action diagonale de $\T^n$ sur $\P^N$
$$
*_\geom : \T^n \times \P^N \to \P^N
\quad \quad , \quad \quad  
(s, x) \mapsto (s^{a_0} \, x_0 : \cdots : s^{a_N} \, x_0)
\enspace .
$$ 
On s'int{\'e}resse {\`a} l'adh{\'e}rence de Zariski des orbites de cette action; pour un point $\alpha= (\alpha_0 : \cdots : \alpha_N) \in \P^N$ on pose
$$
X_\arith := \ov{\T^n *_\cA \alpha} \ \subset \P^N
$$
la {\it vari{\'e}t{\'e} torique projective} associ{\'e}e au couple $(\arith)$. 
Autrement-dit, $X_\arith$ est l'adh{\'e}rence de Zariski de l'image de l'application monomiale
$$
\varphi_{\cA,\alpha} := *_\geom|_\alpha:  \T^n \to \P^N 
\enspace , \quad \quad  
s \mapsto (\alpha_{0} \, s^{a_0} : \cdots :   
\alpha_{N} \, s^{a_N} )
\enspace .
$$ 
C'est donc une vari{\'e}t{\'e} torique projective au sens de~\cite{GKZ94}, \cad
une sous-vari{\'e}t{\'e} de $\P^N$ stable par rapport {\`a} l'action d'un tore
$\T^n$, avec une orbite dense $X_\arith^\circ:=\T^n *_\cA \alpha$. 

Pour simplifier l'exposition, on suppose 
$\alpha \in (\P^N)^\circ$ dans le reste de cette introduction, 
et on fixe d{\'e}sormais un syst{\`e}me de coordonn{\'e}es projectives 
de $\alpha$ dans $(K^\times)^{N+1}$ o{\`u} $K$ est un corps de nombres 
appropri{\'e}. 
Soit $L_\geom \subset \Z^n $ le sous-module engendr{\'e} par 
les diff{\'e}rences des vecteurs $a_0, \dots, a_N$. 
On supposera aussi $L_\cA=\Z^n$, ce qu'en particulier implique
$\dim(X_\cA)=n$. 
Nous renvoyons au 
paragraphe~\ref{monomiales} pour le cas g{\'e}n{\'e}ral. 

On notera $X_\geom$ la vari{\'e}t{\'e} torique associ{\'e}e {\`a}  $\geom \in
(\Z^n)^{N+1}$ et $( 1, \dots, 1)  \in (\Q^\times)^{N+1}$. Dans ce cas,
l'orbite principale $X_\geom^\circ$ est un sous-groupe du tore
$(\P^N)^\circ$; d'ailleurs tous les sous-groupes alg{\'e}briques connexes
de $(\P^N)^\circ$ sont de cette forme. 
Dans le cas g{\'e}n{\'e}ral
$$
X_\arith^\circ = \alpha \cdot X_\cA^\circ
$$
o{\`u} $\cdot$ d{\'e}signe la multiplication dans
$(\P^N)^\circ$. Autrement-dit, l'orbite principale de $X_\arith$ est
le translat{\'e} d'un sous-groupe. 

\medskip 

Suivant la philosophie g{\'e}n{\'e}rale de la th{\'e}orie des vari{\'e}t{\'e}s toriques,
la plupart des propri{\'e}t{\'e}s g{\'e}om{\'e}triques de ces vari{\'e}t{\'e}s peuvent se
traduire en des {\'e}nonc{\'e}s combinatoires sur  les vecteurs $a_0, \dots,
a_N \in \Z^n$ d{\'e}finissant 
l'action. 
En ce qui concerne la th{\'e}orie de l'intersection g{\'e}om{\'e}trique de ces
vari{\'e}t{\'e}s, le r{\'e}sultat le plus 
fondamental est que le degr{\'e} s'identifie au volume $n$-dimensionnel de
l'enveloppe 
convexe $Q_\geom:= \Conv(a_0, \dots, a_N)\subset \R^n$: 
$$
\deg(X_\geom) = n! \, \Vol_{n} \Big( Q_\cA \Big)
\enspace. 
$$ 
Le but principal du pr{\'e}sent texte est d'{\'e}tablir un analogue arithm{\'e}tique de ce r{\'e}sultat. Soit $M_K$ l'ensemble des places du corps $K$, pour chaque $v \in M_K$ on consid{\`e}re le vecteur $\tau_{\alpha v} := (\log|\alpha_0|_v,\dots,\log|\alpha_N|_v)\in\R^{N+1}$ et le polytope 
$$
Q_{\cA,\tau_{\alpha v}}:= \Conv \Big( (a_0,\log|\alpha_0|_v), 
\dots, (a_N,\log|\alpha_N|_v)  \Big) \ \subset \R^{n+1} 
\enspace, 
$$
dont la {\it toiture} au-dessus de $Q_\cA$ (\cad l'enveloppe sup{\'e}rieure) s'envoie bijectivement sur $Q_\cA$ par la projection standard $ \R^{n+1} \to \R^n$. 
On pose alors 
$$
\vartheta_{\cA,\tau_{\alpha  v}}  : Q_\cA \to \R
\enspace , \quad \quad 
x \mapsto \max \Big\{ y\in\R\, : \ (x,y) \in Q_{\cA,\tau_{\alpha v}} \Big\} 
$$ 
la param{\'e}trisation de cette toiture; c'est une fonction {\it concave} 
et {\it affine par morceaux}. 

\begin{thm}  \label{thm1} 
Soit $\geom \in (\Z^n)^{N+1}$ tel que $L_\geom= \Z^n$
et $\alpha \in (K^\times)^{N+1}$, alors
$$
\wh{h}(X_{\cA,\alpha}) = (n+1)! \, \sum_{v \in M_K} 
\frac{[K_v:\Q_v]}{[K:\Q]} \, \int_{Q_\cA} \,
\vartheta_{\cA,\tau_{\alpha v}} \ dx_1 \cdots dx_n 
\enspace .
$$
\end{thm} 

Notons que $\tau_{\alpha v} = 0$ pour presque tout $v$, donc cette
somme
ne contient qu'un nombre fini de termes non nuls.

Illustrons cette formule sur un exemple: 
consid{\'e}rons la courbe cubique  $C \subset \P^3$, adh{\'e}rence de l'image de l'application $ \T^1 \to \P^3$, 
$\displaystyle s \mapsto \Big(1 : 4\, s: \frac{1}{3}\, s^2: \frac{1}{2} \,
s^3\Big)$. 
Les figures suivantes montrent les polytopes associ{\'e}s et leurs
toitures, pour chaque place $v \in M_\Q$:


\vspace{49mm} 

\begin{figure}[htbp]

\begin{picture}(0,0) 

\put(30,50){
\epsfig{file=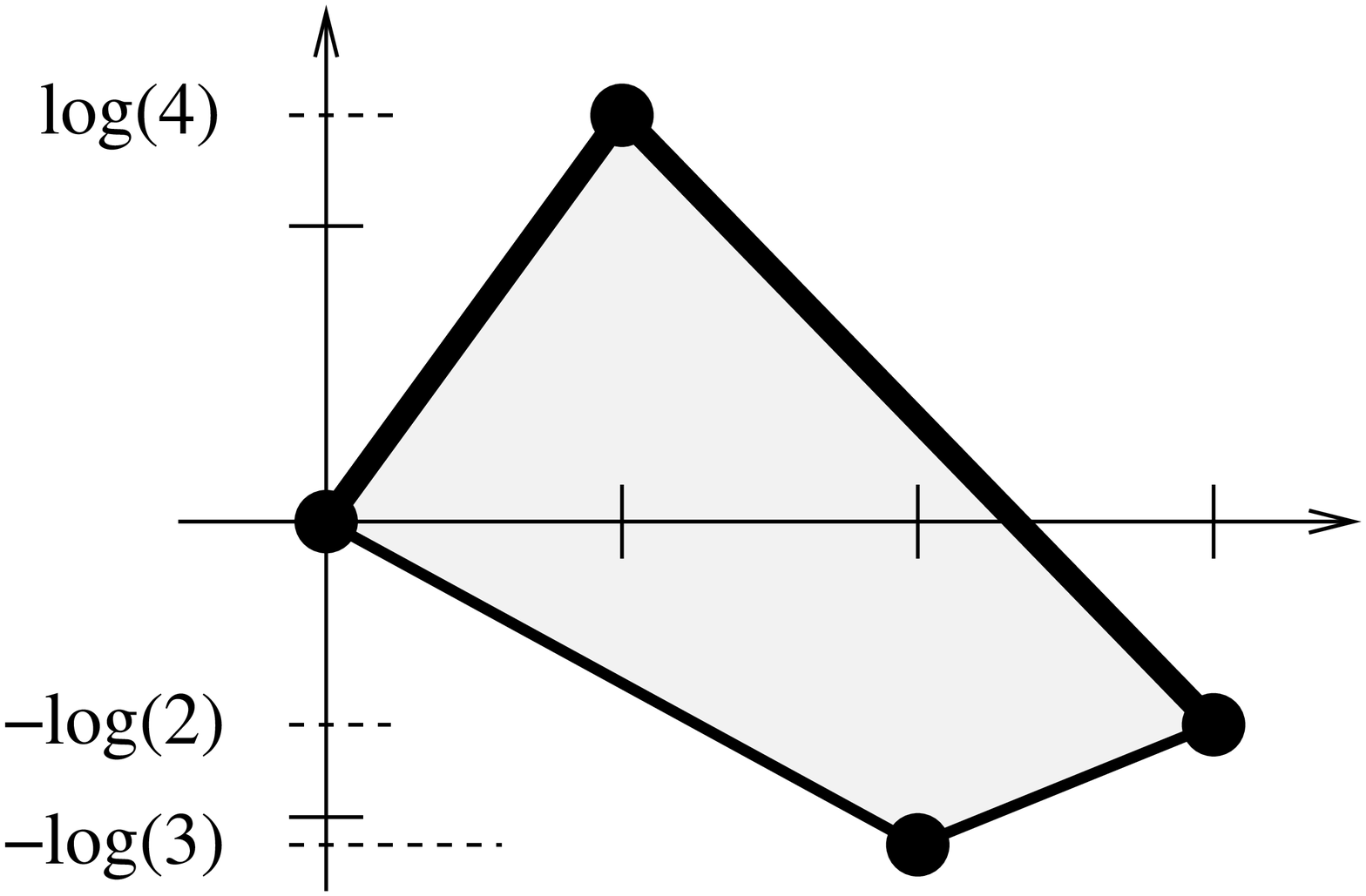, height= 32.5 mm} 
} 
\put(160,105){$v=\infty$} 

\put(260,43){
\epsfig{file=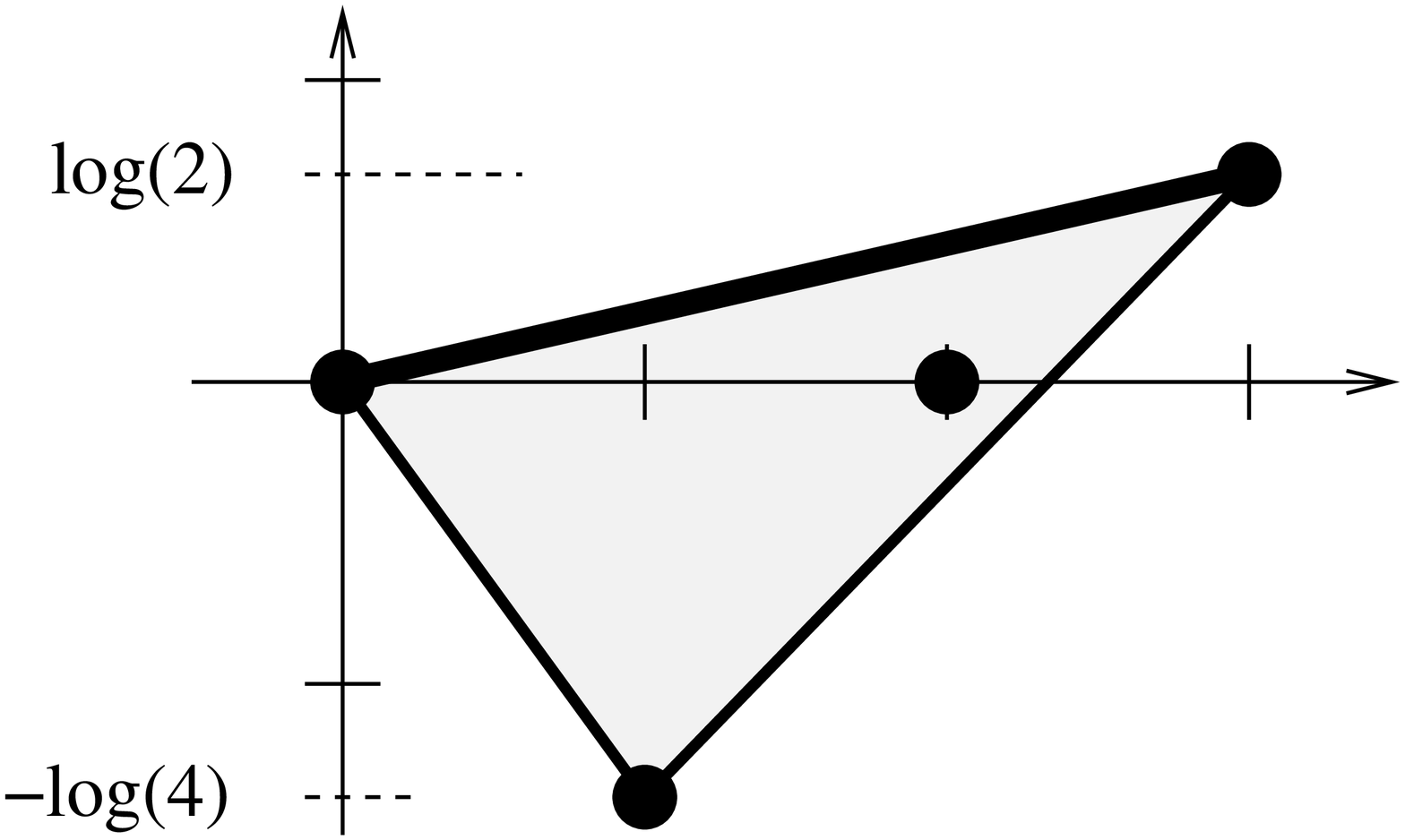, height= 30 mm} 
} 
\put(390,125){$v=2$} 

\put(34,-35){
\epsfig{file=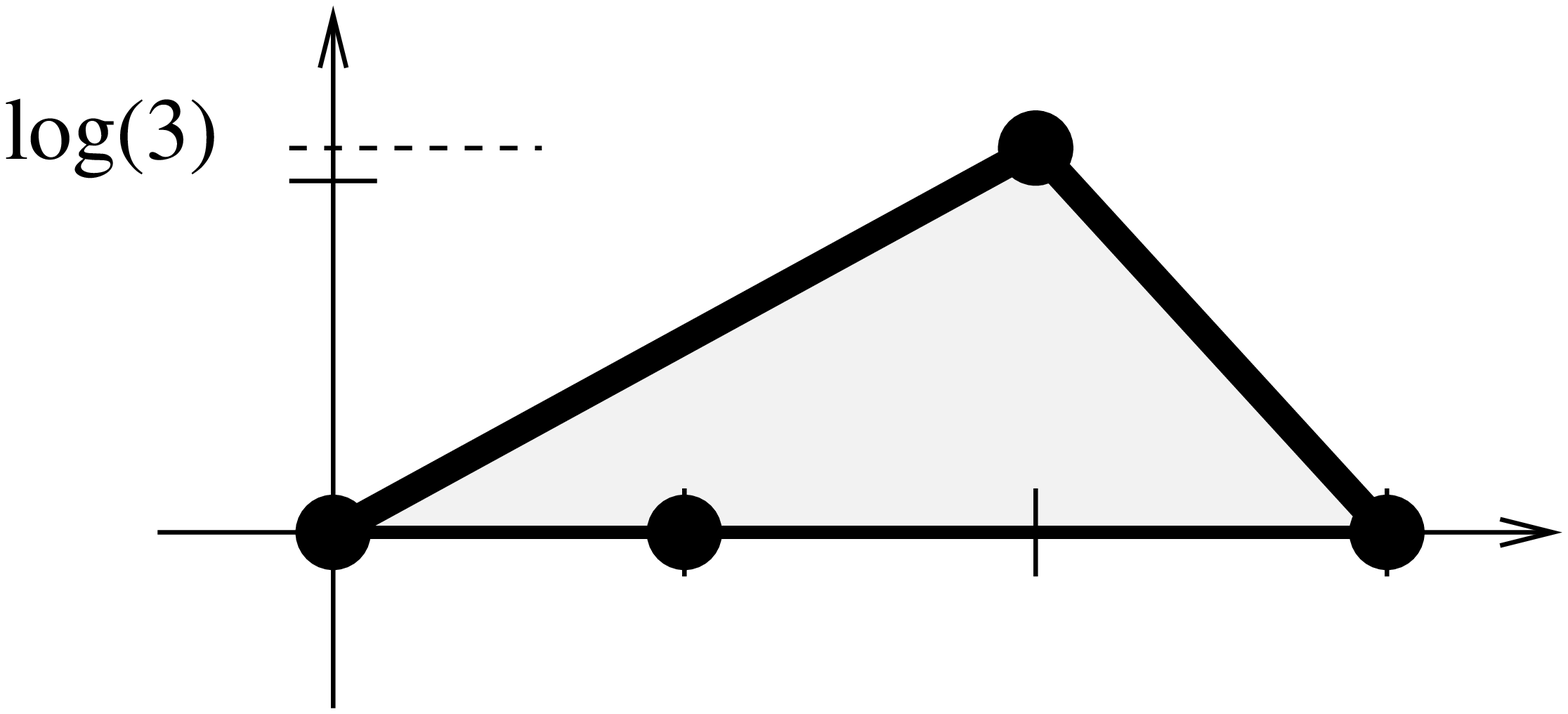, height= 22 mm} 
} 
\put(160,12){$v= 3$} 

\put(278,-35){
\epsfig{file=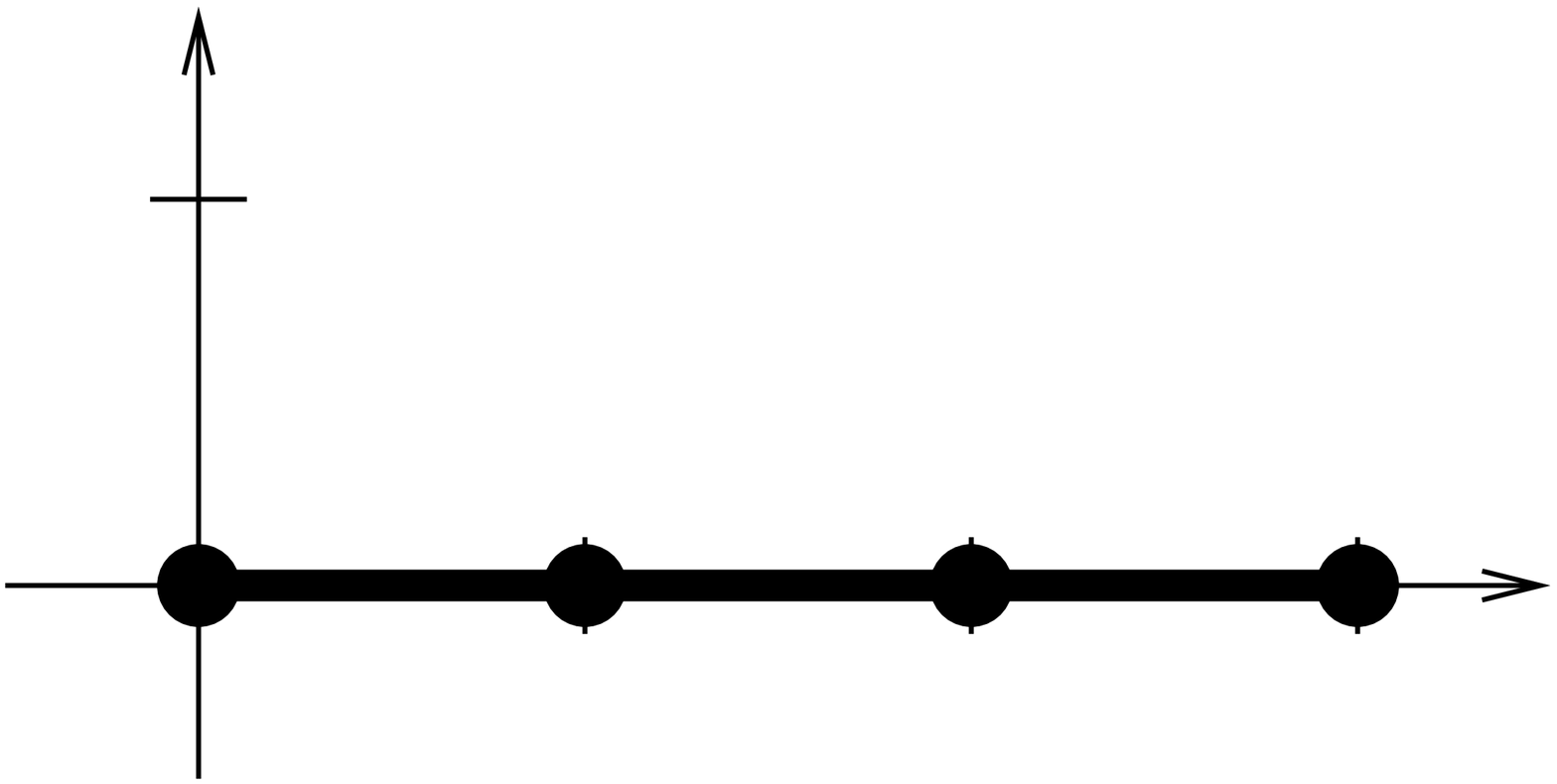, height= 22 mm} 
} 
\put(390,-2){$v\ne \infty, 2,3$}

\end{picture} 

\end{figure}

\vspace{9mm} 


\noindent Ainsi $\vartheta_v \equiv 0$ pour $v \ne \infty, 2,3$, d'o{\`u} 
$$
\hnorm(C) = 2!\,\bigg(\int_0^3\,\vartheta_\infty\, dx\ +\int_0^3\,\vartheta_2\, dx\ +\int_0^3\,\vartheta_3\, dx\bigg) = 7\, \log(2)+3\, \log(3)
\enspace. 
$$

Nous remarquons qu'en dehors des cas des points, hypersurfaces
binomiales et vari{\'e}t{\'e}s de torsion, aucun calcul de ce type n'{\'e}tait
pr{\'e}c{\'e}demment connu. {\'E}galement, ce r{\'e}sultat est bien plus qu'une simple
formule calculant $\hnorm(X_\arith)$. L'ensemble des fonctions
ainsi construit
$$
\Theta_\arith:= \Big( 
\vartheta_{\cA,\tau_{\alpha v}} \, : \ Q_\geom \to \R\Big)_{v \in M_K}
$$ 
se comporte comme un analogue arithm{\'e}tique du polytope $Q_\cA$. En plus de satisfaire {\`a} la formule pour la hauteur, cet objet est compatible avec la
d{\'e}composition en orbites, 
{\it voir}~\S~\ref{orbites}. 
Nous montrons aussi qu'il fournit une
traduction combinatoire naturelle des op{\'e}rations simples telles que
la formation de joints, produits de Segre et plongements de 
Veronese, {\it voir}~\S~\ref{joints}.

\smallskip 

Comme cons{\'e}quence imm{\'e}diate de ce r{\'e}sultat, on v{\'e}rifie 
$\hnorm(X_\arith) \in \log(\Qbar \cap \R_+^\times)$; 
on en d{\'e}duit que le nombre  $\wh{h}(X_\arith)$ est 
{\it transcendant} d{\`e}s qu'il est non nul, 
gr{\^a}ce au th{\'e}or{\`e}me de Baker ({\it voir} par
exemple~\cite[Thm.~1.6]{Wal00}).  
D'apr{\`e}s le th{\'e}or{\`e}me de Zhang~\cite{Zha95a} 
(\cad l'analogue du probl{\`e}me de Bogomolov
sur $\T^N$) 
ceci {\'e}quivaut {\`a} ce que 
$X_\arith$ ne soit pas une sous-vari{\'e}t{\'e} de torsion; 
ainsi on a 

\begin{cor} \label{transcendance} 
Soit $\geom \in (\Z^n)^{N+1}$ et $\alpha \in (K^\times)^{N+1}$
tel que que $X_\arith$ n'est pas une vari{\'e}t{\'e} de torsion, 
alors $\wh{h}(X_\arith) \notin \Qbar$. 
\end{cor}

Un couple  $(\arith)$ est dit  {\it sym{\'e}trique} si la vari{\'e}t{\'e} 
$X_\arith$ est sym{\'e}trique par rapport {\`a} 
l'inversion $[-1]: (\P^N)^\circ \to (\P^N)^\circ$, $(x_0:\cdots:x_N) \to (x_0^{-1}:\cdots:x_N^{-1})$, 
\cad si $X_{\cA,\alpha^{-1}} = X_{\arith}$. Dans ce cas, 
le r{\'e}sultat peut se reformuler en une expression dont les
termes locaux sont des volumes de polytopes
(Corollaire~\ref{hnormvol})
$$ 
\wh{h}(X_{\cA,\alpha}) = \frac{(n+1)!}{2}  \, \sum_{v \in M_K} 
\frac{[K_v:\Q_v]}{[K:\Q]} \, \Vol_{n+1} (Q_{\cA,\tau_{\alpha v}})
\enspace. 
$$

\medskip 

Plus g{\'e}n{\'e}ralement, il est possible 
d'attacher {\`a} une vari{\'e}t{\'e} une 
{\it multihauteur normalis{\'e}e} 
relative {\`a} plusieurs plongements projectifs.  
Soient 
$$
\cA_i \in (\Z^n)^{N_i+1}
\enspace , \quad \quad 
\alpha_i \in (K^\times)^{N_i+1}
\enspace, 
$$
tels que $L_{\cA_i}= \Z^n$ pour $ i=0, \dots, n$. Chaque couple
$(\cA_i , \alpha_i)$ d{\'e}finit un plongement monomial $\varphi_{i}:
\T^n \to \P^{N_i}$ et 
on note
$$
\wh{h} \Big((\cA_0, \alpha_0), \dots, (\cA_n, \alpha_n); 
\T^n \Big) \ \in \R_+ 
$$
la multihauteur normalis{\'e}e du tore $\T^n$ 
relative {\`a} ces plongements;  
la hauteur $\hnorm(X_\arith)$ correspond au cas $(\cA_i , \alpha_i) = (\arith)$ pour tout $i$. 
Nous donnons aussi une expression explicite pour cette
multihauteur. Pour cela, nous introduisons la notion 
d'{\it int{\'e}grale mixte} (ou {\it multi-int{\'e}grale}) 
$\MI(f_0, \dots, f_n)$ d'une famille de 
fonctions concaves $f_i: Q_i\to \R$ d{\'e}finies sur des  ensembles 
convexes compacts $Q_i \subset \R^n$. Ceci g{\'e}n{\'e}ralise l'int{\'e}grale
d'une fonction concave, on renvoie au 
paragraphe~\ref{demo thm2} pour sa d{\'e}finition et  
propri{\'e}t{\'e}s de base.

Pour $i=0, \dots, n$ et  $v \in M_K$ on d{\'e}signe par $\vartheta_{i, v}
: Q_{\cA_i} \to \R$ la fonction param{\'e}trant la toiture du polytope
$Q_{i,v} \subset \R^{n+1}$ 
au-dessus de
$Q_{\cA_i}$
associ{\'e} au poids 
$\tau_{\alpha_i \, v}:=(\log|\alpha_{i\, 0}|_v,\dots, $
$\log|\alpha_{i\, N_i}|_v)$.   

\begin{thm} \label{thm2} 
Soit $\cA_0 \in (\Z^n)^{N_0+1}, \dots, \cA_n \in (\Z^n)^{N_n+1}$ 
tels que $L_{\cA_i} = \Z^n$  pour $i=0, \dots, n$ 
et $\alpha_0 \in (K^\times)^{N_0+1}, \dots, \alpha_n \in (K^\times)^{N_n+1}$,
alors
$$
\wh{h} \Big((\cA_0, \alpha_0), \dots, (\cA_n, \alpha_n); \T^n \Big) 
= \sum_{v\in M_K} \frac{[K_v : \Q_v]}{[K:\Q]}
\, \MI \Big( \vartheta_{0, v}, \dots, \vartheta_{n, v} \Big)
\enspace.
$$
\end{thm} 
Ceci est un analogue arithm{\'e}tique de l'expression pour le multidegr{\'e} 
du tore comme un {\it volume mixte} (ou {\it multi-volume}) 
$\deg_{\cA_1, \dots, \cA_n} (\T^n) 
= \MV(Q_{\cA_1}, \dots, Q_{\cA_n})$~\cite[p.~116]{Ful93}.

\medskip

La d{\'e}monstration du th{\'e}or{\`e}me~\ref{thm1} suit une d{\'e}marche indirecte: au lieu d'utiliser la d{\'e}finition de la hauteur normalis{\'e}e, on 
s'appuie sur le calcul d'une fonction de Hilbert arith\-m{\'e}\-ti\-que appropri{\'e}e. L'un des principaux obstacles {\`a} surmonter dans une telle approche est de trouver une fonction de type Hilbert dont l'expression asymptotique soit li{\'e}e {\`a} la hauteur normalis{\'e}e; les diff{\'e}rentes variantes {\'e}tudi{\'e}es jusqu'{\`a} pr{\'e}sent sont li{\'e}es {\`a} la hauteur projective~\cite{GS91}, \cite{Zha95a}, \cite{AB95},
\cite{Ran01}. 

Soit $X \subset \P^N$ une $K$-sous-vari{\'e}t{\'e} de 
dimension $n$ et $I(X) \subset K[x_0, \dots, x_N]$ son id{\'e}al homog{\`e}ne de d{\'e}finition. Nous proposons ici une fonction de Hilbert arithm{\'e}tique 
${\cH}_{\rm norm}(X;\cdot)$ associant {\`a} un entier $D$ donn{\'e}, la
hauteur de Schmidt du $K$-espace lin{\'e}aire $I(X)_D  \subset  K[x_0,
\dots, \newline x_N]_D$. 
Se pose alors la question du comportement asymptotique de cette fonction, {\`a} laquelle nous apportons une r{\'e}ponse pour le cas des vari{\'e}t{\'e}s toriques: 

\begin{prop} \label{asymptotique} 
Soit $\cA \in (\Z^n)^{N+1}$ tel que $L_\cA = \Z^n$
et $\alpha\in (\Qbar^\times)^{N+1}$, alors 
$$
\cH_{\rm norm}(X_\arith;D) = 
\frac{\wh{h}(X_\arith)}{(n+1)!} \, D^{n+1} + O(D^n)
\enspace. 
$$ 
lorsque $D$ tend vers $ \infty$. 
\end{prop}

Par analogie avec les fonctions de Hilbert arithm{\'e}tiques connues, on pourrait conjecturer qu'un tel r{\'e}sultat reste valable pour une sous-vari{\'e}t{\'e} quelconque. Il serait int{\'e}ressant d'{\'e}tendre dans ce sens la proposition~\ref{asymptotique}, malheureusement notre approche reste confin{\'e}e au cas des vari{\'e}t{\'e}s toriques. 
 
\smallskip 

Pour {\'e}tablir cette proposition, notre d{\'e}marche 
consiste {\`a} {\'e}tudier l'asymptotique de la fonction $\cH_{\rm
  norm}(X_\arith;D)$, d'o{\`u} nous d{\'e}duisons une expression explicite 
du coefficient dominant; puis nous identifions ce coefficient avec la hauteur normalis{\'e}e de $X_\arith$. 
Tout d'abord on montre que cette fonction s'{\'e}crit comme somme de
contributions locales, dont chaque terme s'identifie {\`a} un certain {\it
  poids de Hilbert}: 
on a (Proposition~\ref{decom Hnorm}) 
$$
\Hnorm(X_\arith ; D) = \sum_{v \in M_K} 
\frac{[K_v:\Q_v]}{[K:\Q]} \, s_{\tau_{\alpha v}} (X_\geom; D)  
\enspace.
$$
Gr{\^a}ce {\`a} un r{\'e}sultat de D.~Mumford~\cite{Mum77} 
on sait que chaque poids de Hilbert satisfait {\`a} une formule
asymptotique dont le terme dominant est $\frac{D^{n+1}}{(n+1)!}$ fois
le {\it poids de Chow} $e_{\tau_{\alpha v}}(X_\cA)$ associ{\'e}. On en d{\'e}duit 
$
\displaystyle \Hnorm(X_\arith ; D) = \frac{c(X_\arith)}{(n+1)!} \, D^{n+1} + O(D^n)
$
avec
$$
c(X_\arith)  :=  \sum_{v \in M_K} \frac{[K_v:\Q_v]}{[K:\Q]} \, 
e_{\tau_{\alpha v}} (X_\geom)
\enspace. 
$$
Nous d{\'e}montrons en m{\^e}me temps une expression int{\'e}grale explicite pour les poids de Chow d'une vari{\'e}t{\'e} torique projective~(Proposition~\ref{poid-torique}).  

L'expression obtenue pour $c(X_\arith)$ montre que le quotient
$\frac{c(X_\arith)}{\deg(X_\arith)}$ est {\it lin{\'e}aire} par rapport
aux applications puissances $[k]: \P^N \to \P^N$. 
En outre, on montre que la constante $c(X_\arith)$ se compare {\`a} la
hauteur projective $h(X_\arith)$, en s'appuyant fortement sur un
r{\'e}sultat de H.~Randriambololona~\cite[Thm. A.2.2]{Ran01}. Ces deux conditions suffisent {\`a} montrer l'identit{\'e} $c(X_\arith) =  \hnorm(X_\arith)$, ce qui ach{\`e}ve d'{\'e}tablir {\`a} la 
fois le th{\'e}or{\`e}me~\ref{thm1} et la proposition~\ref{asymptotique} ci-dessus. 

\smallskip  

Comme sous-produit de notre d{\'e}monstration nous obtenons l'identit{\'e} attrayante
$$
\wh{h}(X_\arith) = \sum_{v \in M_K} 
\frac{[K_v:\Q_v]}{[K:\Q]} \, 
e_{\tau_{\alpha v}} (X_\geom)\enspace.   
$$
La notion de poids de Chow a {\'e}t{\'e} introduite dans le contexte 
de la th{\'e}orie g{\'e}om{\'e}trique des invariants~\cite{Mum77}. 
Les travaux r{\'e}cents de J-H.~Evertse et R.~Ferretti la font {\'e}galement intervenir de mani{\`e}re cruciale dans l'approximation diophantienne sur les vari{\'e}t{\'e}s projectives~\cite{EF02}. L'identit{\'e} ci-dessus fait {\`a} nouveau appara{\^\i}tre le poids de Chow en g{\'e}om{\'e}trie arithm{\'e}tique, mais dans un contexte diff{\'e}rent.

\smallskip

La hauteur normalis{\'e}e est un invariant des sous-vari{\'e}t{\'e}s du tore
particuli{\`e}rement int{\'e}\-res\-sant. 
Nous pr{\'e}senterons dans un prochain
texte les cons{\'e}quences 
des r{\'e}sultats obtenus pour l'approximation diophantienne.

\bigskip

Voyons maintenant l'organisation du pr{\'e}sent texte. 
Dans la section~\ref{generalites} on passe en revue quelques 
g{\'e}n{\'e}ralit{\'e}s sur les vari{\'e}t{\'e}s toriques projectives et sur 
les hauteurs et multihauteurs des sous-vari{\'e}t{\'e}s. 
Dans la section~\ref{hilbert} nous introduisons les fonctions de
Hilbert arithm{\'e}tiques et {\'e}tablissons leurs propri{\'e}t{\'e}s principales. 
Dans la section~\ref{normalisee} 
on 
explicite les poids de Hilbert et de Chow 
d'une vari{\'e}t{\'e} torique, puis 
on d{\'e}montre
l'expression asymptotique pour la fonction de Hilbert arithm{\'e}tique
d'une vari{\'e}t{\'e} torique, 
dont le terme principal s'exprime comme somme des poids de Chow. 
Dans
la section~\ref{multi-hauteurs} nous {\'e}tendons nos r{\'e}sultats 
au calcul des multihauteurs et des multipoids de Chow; 
nous {\'e}tudions aussi le comportement de la famille de fonctions
$\Theta_\arith$ introduite, par rapport {\`a} 
la d{\'e}composition en orbites sous l'action 
$*_\cA$, formation de joints, produits de Segre et 
plongements de Veronese.


\bigskip 

\noindent{\bf Remerciements.\hspace{1.5mm}---}\hspace{1mm} 
Nous remercions Vincent Maillot pour plusieurs discussions en rapport
avec travail, 
ainsi que Hughes Randriambololona 
pour nous avoir expliqu{\'e} ses r{\'e}sultats sur les fonctions de 
Hilbert arithm{\'e}tiques.


\bigskip 

\setcounter{tocdepth}{2}

\typeout{Matieres}

\tableofcontents



\renewcommand{\thesection}{\Roman{section}}

%
%


\typeout{Generalites}

\section{G{\'e}n{\'e}ralit{\'e}s}

\label{generalites}

\setcounter{equation}{0}
\renewcommand{\theequation}{\thesection.\arabic{equation}}

On d{\'e}signe par  $\Q$ le corps des nombres rationnels,
$\Z$ l'anneau des entiers rationnels,
$K$ un corps de nombres et
$\cO_K$ son anneau d'entiers.
On note $\R$ le corps des nombres r{\'e}els et 
$\C$ le corps des nombres complexes; on pose
$\R_+$ et $\R_+^\times$
 les ensembles des nombres r{\'e}els non-n{\'e}gatifs et positifs,
respectivement.
On note $\N$ et $\N^\times$ les entiers naturels avec et sans
0, respectivement. 
Pour $N, D \in \N$ 
on pose $\N^{N+1}_D:= \Big\{a\in\N^{N+1} \, : \ a_0+\dots+a_N=D \Big\}$.
On note encore $\K$ pour un corps arbitraire, alg{\'e}briquement clos.

Pour chaque premier rationnel $p$  on note $ | \cdot |_p $ la
valeur
absolue $p$-adique sur $\Q$ telle que
$|p|_p=p^{-1}$; 
on note aussi  $|\cdot|_\infty$
ou simplement $ | \cdot | $
la valeur absolue standard. 
Celles-ci forment un ensemble complet de valeurs absolues sur
$\Q$:
on identifie l'ensemble $M_\Q$
de ces valeurs absolues {\`a} l'ensemble $\{ \infty, \, p
\, : \, p \ \mbox{premier} \}  $.
Plus g{\'e}n{\'e}ralement, on d{\'e}signe par $M_K$ l'ensemble des valeurs absolues
de $K$ {\'e}tendant
les valeurs absolues de $M_\Q$, 
et
on note $M_K^\infty$ le sous-ensemble de $M_K$
des valeurs absolues
archim{\'e}diennes.

On note  $\T^n$ le tore alg{\'e}brique
et $\P^N$ l'espace projectif
sur $\Qbar$ ou sur $\K$, suivant le contexte. 
Une vari{\'e}t{\'e} est toujours suppos{\'e}e r{\'e}duite
et irr{\'e}ductible.


\typeout{Varietes toriques projectives}

\subsection{Vari{\'e}t{\'e}s toriques projectives}

\label{monomiales}

Dans ce paragraphe on donne quelques com\-pl{\'e}\-ments d'information
concernant les vari{\'e}t{\'e}s toriques, 
nos r{\'e}f{\'e}rences sont~\cite{GKZ94},~\cite{Ful93}. 

\medskip 

On se place sur un corps de base $\K$ alg{\'e}briquement clos. 
Soit $n,N\in \N$ des entiers, 
$\cA = ( a_0 , \dots, a_N) \in (\Z^n)^{N+1}$, 
$\alpha= (\alpha_0,\dots, \alpha_N)\in \P^N$; et notons par 
$X_\arith \subset\P^N:= \P^N(\K)$ la vari{\'e}t{\'e} torique projective
et  
$Q_\cA \subset \R^n$ le polytope associ{\'e}s, {\it voir} 
l'introduction.

Lorsque le point $\alpha$ est contenu dans un sous-espace 
standard $E \cong
\P^M$ de $\P^N$, la sous-vari{\'e}t{\'e} $X_\arith$
toute enti{\`e}re reste dans cet espace, puisque l'action $*_\cA$ est 
diagonale. 
Quitte {\`a} se restreindre {\`a} un sous-espace standard, 
on peut donc supposer s.p.d.g. $\alpha \in (\P^N)^\circ$ et on fixe 
d{\'e}sormais un
syst{\`e}me de coordonn{\'e}es projectives de $\alpha$ 
dans $(\K^\times)^{N+1}$.

\smallskip

Le couple $(\arith)$ 
peut s'interpr{\'e}ter 
comme une suite finie de mon{\^o}mes 
$\alpha_{0}  \,  s^{a_0}, \dots, \alpha_{N} \, s^{a_N}$ 
de l'anneau des polyn{\^o}mes de Laurent 
$\K[s_1^{\pm 1} , \dots, s_n^{\pm 1}]$, on dit alors 
que les $a_i$ sont ses {\em exposants} et 
les $\alpha_i$ ses {\em coefficients}. 

Posons 
$
\Expo(\cA) := \{ a_{i_0}, \dots, a_{i_M} \}  \ \subset \Z^n
$
l'ensemble des exposants {\em distincts} dans $\cA$, alors 
$X_\arith$ et $X_{\expo(\cA)}$ 
sont lin{\'e}airement isomorphes par l'application
$$
\P^N \to \P^M  
\enspace, \quad \quad 
(x_0 : \cdots : x_N ) \mapsto
({\alpha_{i_0}^{-1}} \, x_{i_0} : \cdots : 
{\alpha_{i_M}^{-1}} \, x_{i_M} ) \enspace,
$$
et donc pour leurs propri{\'e}t{\'e}s {\it g{\'e}om{\'e}triques} on 
peut toujours se ramener 
au cas classique o{\`u}
les $a_i$ sont tous distincts et $\alpha_i =1$ pour tout $i$.
Soit $L_\geom \subset \Z^n $ le sous-module engendr{\'e} par 
les diff{\'e}rences des vecteurs $a_0, \dots, a_N$; 
en particulier on trouve  
$
\dim(X_\arith) = \dim(L_\cA)$
car $L_\cA = L_{\expo(\cA)}$.

Le sous-module $L_\cA$ est un r{\'e}seau
de l'espace lin{\'e}aire engendr{\'e} $L_\cA \otimes_\Z \R \subset \R^n$; 
on consid{\`e}re la forme volume $\mu_\cA$ sur 
cet espace lin{\'e}aire, 
invariante par translations et normalis{\'e}e de sorte que 
$$
\mu_\cA(L_\cA \otimes_\Z \R/ L_\cA )=1 
\enspace,
$$ 
autrement-dit, de sorte que le volume d'un domaine fondamental soit 1. 
Posons $r:= \dim(L_\cA)$, 
le degr{\'e} de $X_\arith $ s'explicite alors comme le volume 
du polytope associ{\'e}~\cite[p.~111]{Ful93}
$$
\deg(X_\arith) = r! \, \mu_\cA(Q_\cA)
\enspace. 
$$
Soit maintenant  $\eta : \Z^r \hookrightarrow \Z^n$ une 
application lin{\'e}aire
injective telle que 
$\eta(\Z^r) = L_\cA$ et posons
$b_i:= \eta^{-1}(a_i) \in \Z^r$
puis 
$\cB:= (b_0, \dots, b_N) \in (\Z^r)^{N+1} $. 
Alors $X_{\cB, \alpha} = X_\arith$~\cite[Ch.~5,~Prop.~1.2]{GKZ94},
et ainsi on pourra supposer s.p.d.g. 
$L_\cA = \Z^n$.
 Dans cette situation, la forme volume  $\mu_\cA$ 
co{\"\i}ncide avec la forme volume euclidienne 
$\Vol_{n}$ de $\R^n$ et en particulier  
$$
\dim(X_{\cA,\alpha}) = n
\enspace , \quad \quad
\deg(X_{\cA,\alpha}) = n!\, \Vol_{n}(Q_\cA)\enspace.
$$
De plus, notons que sous l'hypoth{\`e}se $L_\cA=\Z^n$ 
l'application $\varphi_{\cA,\alpha}: \T^n \to X_{\cA, \alpha}^\circ$ 
est un isomorphisme g{\'e}om{\'e}trique, 
et 
en particulier $X_{\cA,\alpha}$ est une vari{\'e}t{\'e} rationnelle: 
notons $e_1, \dots, e_n$ la base standard de $\R^n$ et posons 
$$
e_i = \sum_{j=1}^N \lambda_{i\,j} \, (a_j - a_0) 
$$ 
avec $\lambda_{i\, j} \in \Z$, 
l'inverse $\varphi_\arith^{-1}: X_\arith^\circ \to \T^n$ 
est alors
$
\displaystyle
x \mapsto \bigg( \prod_{j=1}^N \bigg( \frac{\alpha_0 \, 
x_j}{\alpha_j \, x_0} \bigg)^{\lambda_{i \, j}} \, : \ 
1 \le i \le n \bigg)$.

Pour $c \in \Z^n$ et $\gamma \in \K^\times$ on consid{\`e}re les  translations
$c + \cA := (c+a_0, \dots, c + a_N) $ et 
$\gamma \cdot \alpha := (\gamma \, \alpha_0, \dots, \gamma \,
\alpha_N)$. 
On v{\'e}rifie imm{\'e}diatement que les sous-vari{\'e}t{\'e}s torique 
$X_{\arith}$ et $X_{c + \cA, \gamma \cdot \alpha}$ 
co{\"\i}ncident; on pourra donc supposer s.p.d.g. $a_0=0$ et $\alpha_0=1$.

\smallskip 

La sous-vari{\'e}t{\'e} $X_\arith$ est d{\'e}finie 
par des {\'e}quations {\it binomiales}. 
De fait, la correspondance 
$X \mapsto I(X)$ est une bijection
entre l'ensemble des sous-vari{\'e}t{\'e}s toriques de $\P^N(\K)$ et 
celui des
id{\'e}aux binomiaux premiers et homog{\`e}nes de
$ \K[x_0, \dots, x_N]$, {\it voir}~\cite{ES96}.


\typeout{Hauteur et multihauteur des varietes projectives}

\subsection{Hauteur et multihauteur des vari{\'e}t{\'e}s projectives}

\label{hauteur normalisee}

Dans ce paragraphe 
nous rappelons les d{\'e}finitions et propri{\'e}t{\'e}s de
 base
des hauteurs et multihauteurs projectives et normalis{\'e}es,
on renvoie {\`a}~\cite[\S~2]{DP99} pour la plupart des d{\'e}tails.

\medskip

On prend maintenant $\Qbar$ comme corps de base. 
Soit $X \subset \P^N=\P^N(\Qbar)$
une $K$-sous-vari{\'e}t{\'e} de dimension $n$; on rappelle d'abord la notion de
hauteur projective de $X$.

Soit
$f = \sum_a f_a \, U^a \in K[U_0, \dots,  U_n]$ un polyn{\^o}me
en $n+1$ groupes $U_i$ de $N_i+1$ variables chacun.
Soit
$S_{N_i+1}:= \Big\{ (z_0, \dots, z_{N_i}) \in \C^{N_i+1}
\, : \ |z_0|^2 + \cdots + |z_{N_i}|^2 = 1 \Big\} $
la sph{\`e}re unit{\'e} de $\C^{N_i+1}$,
{\'e}quip{\'e}e de
la mesure
$\mu_{N_i+1}$
de masse totale 1 invariante par rapport au groupe
unitaire $U(N_i+1)$.
Pour une place archim{\'e}dienne $v \in M_K^\infty$,
la {\it $S_{N_0+1} \times \cdots \times S_{N_n+1}$-mesure} 
de $f$ par rapport {\`a} $v$
est
$$
m_v( f)
:= \int_{S_{N_0+1} \times \cdots \times S_{N_n+1}}
\log |f|_v
\ \mu_{N_0+1} \times \cdots \times \mu_{N_n+1} \enspace.
$$
Pour une place ultram{\'e}trique $ v \in M_K \setminus M_K^\infty$ on 
note $
|f|_v:= \max \Big\{ |f_a|_v \, : \ a \in \N^{N_0+1}
\times \cdots \times \N^{N_n+1} \Big\}$, 
la {\it norme sup} de $f$ par rapport {\`a} $v$.
On pose alors
$$
h(f):=
\sum_{v \in M_K^\infty} \frac{[K_v: \Q_v]}{[K:\Q]} \,
m_v(f)
\ +
\hspace{-3mm}
\sum_{v \in M_K \setminus M_K^\infty} \frac{[K_v: \Q_v]}{[K:\Q]} \,
\log |f|_v \ \ + \sum_{i=0}^n \deg_{U_i}(f) \cdot
\bigg(\sum_{j=1}^{N_i} \frac{1}{2\, j} \bigg) \ \in \R_+
\enspace.
$$

Soit $\Ch_X \in K[U_0, \dots, U_n]$ la {\it forme de Chow} 
de $X$: 
pour $i=0, \dots, n$ on
consid{\`e}re une forme lin{\'e}aire g{\'e}n{\'e}rale
$
L_i:= U_{i\, 0} \, x_{ 0} + \cdots + U_{i \, N} \, x_{ N} 
\ \in K[U_i][x]
$ 
et on pose 
$$
\nabla_{X} := 
\Big\{ (u_0, \dots, u_n) \, 
: \ X \cap Z(L_0(u_0), \dots, L_n(u_n)) \ne \emptyset 
\Big\}
\
\subset (\P^{N})^{n+1} 
$$
pour l'ensemble des $(n+1)$-uplets de formes
lin{\'e}aires dont l'intersection avec $X$ est non vide. 
C'est une hypersurface
et la forme $\Ch_X$ est d{\'e}finie comme une {\'e}quation (unique {\`a} un 
facteur scalaire pr{\`e}s) de $\nabla_X$. 
C'est
un polyn{\^o}me
en $n+1$ groupes de $N+1$ variables chacun,
homog{\`e}ne de degr{\'e} $\deg(X) $ en chaque groupe $U_i$.

Suivant~\cite{Phi95}, la {\it hauteur projective} de $X$ est 
par d{\'e}finition 
$h (X)  := h (\Ch_X)$.
Alternativement, cette hauteur peut se d{\'e}finir {\it via} la
th{\'e}orie d'Arakelov
comme la hauteur $h_{\ov{\cO(1)}}(\Sigma)$,
o{\`u} $\Sigma$ est la cl{\^o}ture de Zariski de $X$
dans $\P^N_{\cO_K}$ 
relative au fibr{\'e} en droites universel
$\ov{\cO(1)}$ muni de la m{\'e}trique
de Fubini-Study~\cite[\S~3.1.3]{BGS94}.

\medskip

La structure de groupe de l'ouvert $(\P^N)^\circ \cong \T^N$ permet
d'introduire une hauteur normalis{\'e}e (ou canonique) pour les
sous-vari{\'e}t{\'e}s 
de l'espace projectif $\P^N$, vu comme une compactification 
{\'e}quivariante du tore $(\P^N)^\circ$.
Soit $ [k]:
\P^N \to \P^N,
(x_0: \cdots : x_N) \mapsto (x_0^k: \cdots : x_N^k) $
l'application puissance $k$-i{\`e}me, la {\it hauteur normalis{\'e}e} de  $X$
est alors d{\'e}finie par 
\begin{equation}\label{defhnormalisee}
\wh{h}(X) = \hnorm_{\P^N}(X) := \deg(X) \cdot \lim_{k\to \infty} \,
\frac{h({[k] \, X})}{k\, \deg([k]\, X)} \ .
\end{equation}
Une propri{\'e}t{\'e} fondamentale est la lin{\'e}arit{\'e} du quotient
$\displaystyle \frac{\hnorm}{\deg}$ relativement aux applications
puissances~\cite[Prop.~2.1(i)]{DP99}:
$$
\frac{\hnorm\Big([k]\, X\Big)}{\deg\Big([k]\, X\Big)} = k \cdot
\frac{\hnorm\Big(X\Big)}{\deg\Big(X\Big)}\enspace.
$$
Cette hauteur se compare {\`a} 
la hauteur projective~\cite[Prop.~2.1(v)]{DP99}:
$$
\Big|\hnorm(X) - h(X) \Big| \le \frac{7}{2} \, \log(N+1) \, \deg (X) 
\enspace .
$$

Dans le cas des points, la hauteur normalis{\'e}e co{\"\i}ncide
avec la hauteur de Weil~\cite[Prop.~2.1(vi)]{DP99}.
Pour les hypersurfaces, elle s'{\'e}crit comme somme de contributions locales aux places
archim{\'e}diennes, chaque terme {\'e}tant la  mesure de Mahler locale
d'une {\'e}quation de d{\'e}finition de $X$~\cite[Prop.~2.1(vii)]{DP99}.

\smallskip

Lorsque $X $ est contenue dans un sous-espace standard $E \cong
\P^M$ de $\P^N$, ses hauteurs normalis{\'e}es 
comme sous-vari{\'e}t{\'e} de $E$ et de $\P^N$ co{\"\i}ncident.
C'est une cons{\'e}quence simple de la d{\'e}finition de $\hnorm$ et
de la propri{\'e}t{\'e} analogue pour la hauteur projective.

Consid{\'e}rons donc le cas o{\`u} $X$ n'est contenue dans aucun
des sous-espaces standard, autrement dit
$X^\circ= X \cap (\P^N)^\circ \ne \emptyset$.
Le tore $\T^N$ agit sur $\P^N$
de mani{\`e}re naturelle {\it via} l'inclusion
$\T^N \cong (\P^N)^\circ \hookrightarrow \P^N$;
on pose alors
$ \Stab_{\T^N} (X):= \{ u \in \T^N \, : \ u \cdot X = X \} $
pour le {\it stabilisateur} de $X$ par rapport {\`a} cette action, 
qui est un sous-groupe de $\T^N$, et on note
$$
\sigma_X(k):= \Card\Big(\Stab_{\T^N}(X) \cap \mbox{\rm Ker} [k]\Big)
\enspace.
$$
Ainsi 
$\displaystyle \deg([k] \, X) = \frac{k^n}{\sigma_X(k)} \, 
\deg(X)$~\cite[Prop.~2.1(i)]{DP99}
et donc
$ \displaystyle
\wh{h}(X) =  \lim_{k\to \infty} \,
\frac{\sigma_X(k) }{k^{n+1}}\, h({[k] \, X})$\enspace.

\medskip

En g{\'e}n{\'e}ral, on peut munir le tore $\T^n$ d'une notion
de degr{\'e} et de hauteur des sous-vari{\'e}t{\'e}s relative 
{\`a} un plongement donn{\'e} $\varphi: \T^n
\hookrightarrow \P^N$, en posant 
$$
\deg_\varphi(Y):= \deg(\varphi(Y))
\enspace, \quad \quad 
h_\varphi(Y) := h(\varphi(Y))
$$
pour une sous-vari{\'e}t{\'e} $Y \subset \T^n$.
On peut normaliser cette hauteur par rapport {\`a} l'application 
de  
multiplication $[k]: \T^n \to \T^n$, $s \mapsto 
s^k= (s_1^k,
\dots, s_n^k)$  par la formule 
$$
\hnorm_\varphi(Y):= \deg_\varphi(Y) \cdot 
\lim_{k \to \infty} \frac{h_\varphi( [k] \, Y)}{k \, \deg_\varphi([k] \, Y)}
$$ 
lorsque cette limite existe;  
dans ce cas,
on a automatiquement 
$\displaystyle 
\frac{\hnorm_\varphi([k]\, Y)}{\deg_\varphi([k] \, Y)} 
= 
k\cdot \frac{\hnorm_\varphi(Y)}{\deg_\varphi( Y)}$. 
La hauteur normalis{\'e}e $\hnorm_{\P^N}$ 
correspond
{\`a} l'inclusion standard  $\T^N \hookrightarrow \P^N$, $(s_1, \dots, s_N) 
\mapsto (1: s_1: \cdots: s_N)$.  
Un autre exemple de cette situation appara{\^\i}t dans~\cite[\S~2]{DP99}, o{\`u} l'on 
consid{\`e}re 
pour des sous-vari{\'e}t{\'e}s 
de $\T^n$, 
la hauteur normalis{\'e}e  
induite par l'inclusion $\T^n \hookrightarrow (\P^1)^n 
{\hookrightarrow \atop \mbox{\scriptsize Segre}}  \P^{2^n-1}$; 
{\it voir} aussi le paragraphe~\ref{Demonstration du 
theoreme 1}.

Consid{\'e}rons maintenant 
un plongement projectif 
$\varphi : \T^n \hookrightarrow \P^N$ {\it {\'e}quivariant},  
\cad tel que l'action naturelle du tore sur l'image $\varphi(\T^n)$ 
s'{\'e}tende en une action $\T^n \times \P^N \to \P^N$. 
On montre que la hauteur normalis{\'e}e correspondante
est bien d{\'e}finie et se r{\'e}duit  {\`a} una hauteur $\hnorm_{\P^N}$.

\begin{prop} \label{h_varphi} 
Soit $\varphi: \T^n \hookrightarrow \P^N $ un plongement 
{\'e}quivariant
et $Y$ une sous-vari{\'e}t{\'e} de $\T^n$, 
alors il existe $\cB \in (\Z^n)^{N+1}$ tel que 
$ \displaystyle
\hnorm_{\varphi}(Y) = \hnorm_{\P^N}(\varphi_\cB(Y))$. 

\end{prop}

\begin{demo} 
D'apr{\`e}s~\cite[Ch.~5, Prop.~1.5]{GKZ94} il existe un vecteur 
$\cB \in (\Z^n)^{N+1}$ et une transformation lin{\'e}aire projective  
$U \in \PSL(N+1, \Qbar)$ 
tels que $\varphi(\T^n)= U \circ \varphi_\cB(\T^n)$. 
Ainsi
$ \displaystyle 
\Big| h(\varphi([k]\,Y)) - h(\varphi_\cB([k]\,Y)) \Big|
\le c \, \deg(\varphi_\cB([k]\,Y))
$
o{\`u} $c>0 $ ne d{\'e}pend ni de $k$ ni de $Y$~\cite[Lem.~2.7]{KPS01}, 
ce qui entra{\^\i}ne
$$
\lim_{k \to \infty}
\frac{h(\varphi([k]\, Y))}{k \, \deg( \varphi([k]\, Y))}
= \lim_{k \to \infty}
\frac{h(\varphi_\cB([k]\, Y))}{k \, \deg(\varphi_\cB([k]\, Y))}
= \frac{\hnorm_{\P^N}(\varphi_\cB(Y))}{\deg(\varphi_\cB(Y))}
\enspace, 
$$
car l'application $\varphi_\cB: \T^n \to (\P^N)^\circ $ est un
morphisme de groupes et donc commute avec la multiplication $[k]$. 
\end{demo} 

Pour $\cA \in (\Z^n)^{N+1}$ et $\alpha \in
(K^\times)^{N+1}$ 
on  peut prendre $\cB:= \cA$ dans l'{\'e}nonc{\'e} pr{\'e}c{\'e}dent, d'o{\`u}
$$
\hnorm_{\varphi_\arith}(Y) = \hnorm_{\P^N}(\varphi_\cA(Y))
\enspace
$$ 
 
Insistons sur le fait qu'il ne faut pas confondre 
$\hnorm_\varphi= \wh{h\circ \varphi}$ (la normalisation de
la hauteur induite par $\varphi$)
et
$\hnorm \circ \varphi$ 
(la hauteur 
induite par la hauteur normalis{\'e}e 
de $\P^N$ {\it via} le plongement $\varphi$).

\medskip

Dans la suite on g{\'e}n{\'e}ralise la hauteur normalis{\'e}e 
au cas
des sous-vari{\'e}t{\'e}s d'un produit d'espaces projectifs.
Pour cela on fait appel aux formes
r{\'e}sultantes d'id{\'e}aux multihomog{\`e}nes;
on renvoie {\`a}~\cite{Rem01a} et~\cite{Rem01b} pour une exposition
compl{\`e}te du sujet.

Pour $i=0, \dots, m$ on fixe
un groupe  $x_i= \{ x_{i \, 0}, \dots, x_{i\, N_i}\} $
de $N_i+1$
variables chacun 
et on consid{\`e}re
l'anneau $K[x_0, \dots, x_m]$, multigradu{\'e} par
$\deg(x_{i\, j}) := e_i$
o{\`u} $e_0, \dots, e_m \in \Z^{m+1}$
d{\'e}signent les vecteurs de la base standard de $\R^{m+1}$.

Soit $I \subset K[x_0, \dots, x_m]$ un id{\'e}al
multihomog{\`e}ne
de dimension de Krull
$m+n+1$
et 
$ d_0, \dots, d_n \in \N^{m+1} $;
pour chaque $ d_i \in \N^{m+1}$
on introduit un groupe de variables $U_i= \{U_{i\, 0},
\dots, U_{i\, N_i}\}$ et la forme g{\'e}n{\'e}rale de multidegr{\'e} $d_i$
$$
F_i := \sum_{\deg(a) = d_i} U_{i\, a} \, x^a \ \in K[U_i][x_0, \dots,
x_m]
\enspace.
$$
On pose alors
$$
\cM=M(I; \, d_0, \dots, d_m)
:= K[U_0, \dots, U_n][x_0, \dots, x_m]
/ \Big( I+(F_0, \dots, F_n) \Big)
$$
qui est un
$K[U]$-module
multigradu{\'e}; on note $\cM_q$ sa partie de multidegr{\'e} $q$
pour $q \in \Z^{m+1}$.

\smallskip

De fa{\c c}on g{\'e}n{\'e}rale, {\`a} tout $K[U]$-module $M$ de type fini et de 
torsion (\cad tel que 
$\Ann_{K[U]}(M) $ $ \ne 0$) est associ{\'e}e la forme
$$
\chi(M) :=\prod_{f\in \irr(K[U])}  f^{\ell(M_{(f)})} \ \in K[U] \setminus \{
0\}
$$
o{\`u} $f$ parcours l'ensemble des {\'e}l{\'e}ments irr{\'e}ductibles
de $K[U]$ et $\ell(M_{(f)})$ est la longueur du $K[U]_{(f)}$-module
$M_{(f)}$;
si $M$ n'est pas de torsion on pose $\chi(M):= 0$.

\smallskip 

D'apr{\`e}s~\cite[Lem.~3.2]{Rem01a}
on sait que dans la situation pr{\'e}c{\'e}dente, l'application
$q \mapsto \chi(\cM_q)$ est
constante pour $q \gg 0$.
Cette valeur commune
$$
\res_{d_0, \dots, d_n} (I):= \chi\Big(M(I; \, d_0, \dots, d_m)_q \Big)
\enspace, 
\quad \quad
q \gg 0
\enspace, 
$$
est par d{\'e}finition la {\it forme  r{\'e}sultante
de $I$ d'indice $d_0, \dots, d_n$}, selon~\cite[\S 3.2]{Rem01a}.

\smallskip

Soit $Z \subset \P^{N_0} \times \cdots \times \P^{N_m}$
une sous-vari{\'e}t{\'e} de dimension $n$,
son id{\'e}al  $I(Z) \subset K[x_0, \dots, x_m]$ est donc de dimension de
Krull $m+n+1$.
Pour  $c=(c_0, \dots, c_m) \in \N_{n+1}^{m+1}$
on pose
$$
d(c):= (\overbrace{e_0, \dots, e_0}^{c_0}, \dots, \overbrace{e_m,
  \dots, e_m}^{c_m}) \in (\N^{m+1})^{n+1} 
\enspace.
$$
La {\it multihauteur projective d'indice $c$ de $Z$}
est par d{\'e}finition
$h_c(Z) := h\Big(\res_{d(c)}(I(Z))\Big)$, {\it voir}~\cite[\S~2.3]{Rem01b}.
Soient $D_0, \dots, D_m \in \N^*$
des entiers et  $D:= (D_0, \dots, D_m)$, on consid{\`e}re
l'application
$$
\begin{array}{rll}
\Psi_D : & \P^{N_0} \times \cdots \times \P^{N_m} & \to  
\P^{ {D_0+N_0 \choose N_0} \cdots {D_m+N_m \choose N_m} -1} \\[2mm]
&(x_0, \dots, x_m) & \mapsto  \Big( x_0^{b_0} \cdots x_m^{b_m}
\, : \ b_0 \in \N^{N_0+1}_{D_0}, \dots, \ b_m \in \N^{N_m+1}_{D_m} \Big)
\enspace,
\end{array}
$$
composition des plongements de Veronese et de Segre.
On consid{\`e}re aussi l'application lin{\'e}aire diagonale
$$
\Delta_D : \P^{ {D_0+N_0 \choose N_0} \cdots {D_m+N_m \choose N_m} -1}
\to \P^{ {D_0+N_0 \choose N_0} \cdots {D_m+N_m \choose N_m} -1}, 
\quad
(x_b)_b \mapsto  \bigg({D_0 \choose b_0}^{1/2}
\cdots {D_m  \choose b_m}^{1/2} \, x_b \bigg)_b
\enspace,
$$
la composition $\Delta_D \circ \Psi_D$ est le {\it plongement
mixte remodel{\'e}}.
La hauteur projective de l'image
$\Delta_D \circ \Psi_D(Z)$ se
d{\'e}compose alors~\cite[p. 103]{Rem01b}
\begin{equation} \label{dech}
h\Big(\Delta_D \circ \Psi_D(Z)\Big) =
\sum_{c \in \N_{n+1}^{m+1}}
{n+1 \choose c} \, h_c(Z) \, D^{c}
\enspace.
\end{equation}

\smallskip

C'est la propri{\'e}t{\'e} suivante qui permet de d{\'e}finir les
multihauteurs normalis{\'e}es:

\begin{prop} \label{multihnorm}
Soit $Z \subset \P^{N_0} \times \cdots \times \P^{N_m}$
une sous-vari{\'e}t{\'e} de dimension $n$
et $c \in \N^{m+1}_{n+1}$, notons 
 $s:
\P^{N_0} \times \cdots \times \P^{N_m} \to
\P^{(N_0+1) \cdots (N_m+1) - 1}$
le plongement de Segre.
Alors la suite 
$ \displaystyle k \mapsto  \frac{h_c([k] \, Z)}{k\, \deg([k]\, s(Z))}$
converge lorsque $k$ tend vers l'infini. 
\end{prop}

Avec ces notations, on d{\'e}finit 
la {\it multihauteur normalis{\'e}e de $Z$ d'indice $c$}
par
\begin{equation} \label{hnorm_a}
\wh{h}_c(Z) := \deg(s(Z)) \cdot \lim_{k\to \infty} \,
\frac{h_c([k] \, Z)}{k\, \deg([k]\, s(Z))} \ \in \R_+
\enspace.
\end{equation}

\begin{lem} \label{stabilisateur}
Soit $\psi: \T^N \to \T^P$ un morphisme injectif de groupes 
et $Y$ une sous-vari{\'e}t{\'e} de $\T^N$, alors
$ \Stab_{\T^P} (\psi(Y)) = \psi(\Stab_{\T^N}(Y))$.
\end{lem}

\begin{demo}
Soit $u \in \Stab_{\T^N} (Y)$, alors 
$ \psi(u) \cdot \psi(Y) = \psi( u \cdot Y) = \psi(Y)$
d'o{\`u} on obtient 
$
\psi(\Stab_{\T^N}(Y)) \subset \Stab_{\T^P} (\psi(Y)) 
\enspace.
$

Dans l'autre direction, soit $v \in  \T^P  $
tel que $ v \cdot \psi(Y) = \psi(Y)$.
Alors
$v \in \psi(Y)\cdot \psi(Y)^{-1} \subset \psi(\T^{N})$, \cad qu'il existe
$u \in \T^{N}$ tel que $v = \psi(u)$.
Ainsi $ \psi(u \cdot Y ) = v \cdot \psi(Y) = \psi(Y)$
qui implique $u \in \Stab_{\T^N}(Y) $, {\`a} cause de l'injectivit{\'e} de $\psi$.
\end{demo}

\begin{demo}[D{\'e}monstration de la Proposition~\ref{multihnorm}]
On suppose s.p.d.g. $Z \cap \prod_{i=0}^m (\P^{N_i})^\circ  \ne
\emptyset$, sinon il suffit de
se restreindre {\`a} un produit
$\prod_{i=0}^m E_i$ de sous-espaces standard des $\P^{N_i}$.
Le plongement mixte remodel{\'e} 
$\Delta_D \circ \Psi_D$ 
restreint au tore
$\prod_{i=0}^m (\P^{N_i})^\circ$ 
est {\'e}quivariant, donc 
$$
\hnorm_{\Delta_D \circ \Psi_D} (Z) = 
\deg(\Psi_D(Z)) 
\cdot \lim_{k \to \infty}
\frac{h(\Delta_D \circ \Psi_D([k]\, Z))}{k \, \deg(\Psi_D([k]\, Z))}
= {\hnorm(\Psi_D(Z))}
$$
par la proposition~\ref{h_varphi}. 
On remarque aussi que $\Psi_D$ restreint {\`a}
$\prod_{i=0}^m (\P^{N_i})^\circ$ 
est un morphisme de groupes injectif 
et donc $\Stab(\Psi_D(Z)) = \Psi_D(\Stab(Z))$ 
comme cons{\'e}quence du lemme~\ref{stabilisateur}, 
puis 
$\sigma_{\Psi_D(Z)}(k) = \sigma_{Z}(k)$. 
Ainsi~\cite[Prop.~2.1(i)]{DP99}
$$
\deg(\Psi_D( Z)) = \frac{\sigma_{\Psi_D(Z)}(k)}{k^{n}} \,
\deg([k] \,\Psi_D( Z))= \frac{\sigma_{Z}(k)}{k^{n}} \,
\deg(\Psi_D( [k] \, Z)) 
\enspace,
$$
qui, joint {\`a} l'identit{\'e}~(\ref{dech}), 
entra{\^\i}ne
\begin{eqnarray*}
\hnorm(\Psi_D(Z))
& = &
\lim_{k \to \infty}
\Bigg( 
\frac{\deg( \Psi_D(Z))}  
{k \, \deg( \Psi_D([k]\, Z))}
\cdot h(\Delta_D \circ \Psi_D([k]\, Z))
\Bigg) 
 \\[2mm]
& = &
\lim_{k \to \infty}
\Bigg( 
\sum_{c \in \N_{n+1}^{m+1}}
{n+1 \choose c} \,
\frac{\sigma_{Z}(k)}{k^{n+1}} \, h_c([k]\,Z)  \, D^{c} \Bigg)
\enspace.
\end{eqnarray*}
Cette derni{\`e}re limite convergeant pour {\it tout} $D \in (\N^*)^{m+1}$,
on en d{\'e}duit facilement que chaque coefficient
$ \displaystyle 
\frac{\sigma_{Z}(k)}{k^{n+1}} \, h_c([k]\,Z)$ converge 
lorsque $k$ tend vers l'infini.
De plus, 
$\displaystyle
\frac{\sigma_{Z}(k)}{k^{n}} = \frac{\deg(s(Z))}{\deg([k] \, s(Z))}
$
aussi {\`a} cause du lemme~\ref{stabilisateur}, ce qui ach{\`e}ve la
d{\'e}monstration. 
\end{demo}

Comme cons{\'e}quence de cette d{\'e}monstration
on a {\'e}galement la d{\'e}composition
\begin{equation} \label{decomp hnorm} 
\hnorm(\Psi_D(Z)) =
\sum_{c \in \N_{n+1}^{m+1}}
{n+1 \choose c} \, \hnorm_c(Z) \, D^{c} \enspace.
\end{equation}
En particulier, cela montre que 
pour le cas $m=0$
on retrouve 
la hauteur normalis{\'e}e
$\hnorm(Z)$ comme la multihauteur $\hnorm_{n+1}(Z)$.

\smallskip 

Dans la suite, on s'int{\'e}ressera surtout au cas
d'une vari{\'e}t{\'e} $X$ de dimension $n$
munie d'applications rationnelles
$ \varphi_i : X \dashrightarrow \P^{N_i} $ ($i =0, \dots, n$).
On pose
$\Phi : X \to \prod_{i=0}^n\P^{N_i}$,
$x \mapsto (\varphi_0(x), \dots, \varphi_n(x))$
et on {\'e}crit
\begin{equation}\label{defmultihnormalisee}
\wh{h}({\varphi_0, \dots, \varphi_n}; X) :=
\hnorm_{(1, \dots, 1)}\Big(\Phi(X)\Big)
\end{equation}
pour la {\it multihauteur normalis{\'e}e de $X$ relative {\`a}} 
$\varphi_0, \dots, \varphi_n$.

Cette  multihauteur peut aussi se d{\'e}finir
par des m{\'e}thodes arakeloviennes,
en suivant la d{\'e}marche de \cite{Zha95b}, \cite{Mai00}, 
pour le cas d'une vari{\'e}t{\'e} projective $X \subset\P^N$:
notons $\Sigma$
la  cl{\^o}ture de Zariski de $X$
dans $\P^N_{\cO_K}$ et 
pour $i=0, \dots, n$ consid{\'e}rons le fibr{\'e} hermitien
$\ov{L}_i:= \varphi_i^*(\ov{\cO(1)})$, tir{\'e} en arri{\`e}re du
fibr{\'e} en droites
universel $\ov{\cO(1)}$, cette fois-ci
muni de la m{\'e}trique dite {\it canonique}.
Dans cette situation, la multihauteur ci-dessus
co{\"\i}ncide avec
la hauteur $h_{\ov{L}_0, \dots, \ov{L}_n}(X)$
telle qu'elle est d{\'e}finie dans~\cite[Prop.-D{\'e}fn. 5.5.1]{Mai00}.



%
%


\typeout{Fonctions de Hilbert arithmetiques} 

\section{Fonctions de Hilbert arithm{\'e}tiques} 

\label{hilbert}

\setcounter{equation}{0}
\renewcommand{\theequation}{\thesection.\arabic{equation}}

Soit $X\subset\P^N$ une sous-vari{\'e}t{\'e} d{\'e}finie sur un corps de nombres
$K$ et $I:= I(X)\subset K[x_0,\dots,x_N]$ son id{\'e}al de d{\'e}finition, 
on pose
$$
\cH_\geo(X; D) := \dim_K (K[x_0,\dots,x_N]/I)_D = \binom{D+N}{N} - \dim_K (I_D)
$$
pour la classique {\it fonction de Hilbert g{\'e}om{\'e}trique} de $X$. 
Ci-apr{\`e}s, on introduit un analogue arithm{\'e}tique de cette fonction. 
Posons $\ell:=\dim_K (I_D)$ et 
$$
\bigwedge^\ell K[x_0,\dots,x_N]_D
$$
la $\ell$-i{\`e}me puissance ext{\'e}rieure de $K[x_0,\dots,x_N]_D$. 
Pour $v
\in M_K$ et $f \in \bigwedge^\ell K[x_0,\dots,x_N]_D$ on note $|f|_v$
la norme sup des coefficients de $f$ pour la place $v$, par
rapport {\`a} la base standard 
$\Big\{ \bigwedge_{a \in L} x^{a} \, : \ L \subset \N_D^{N+1},
\ \Card(L)= \ell \Big\}$.

\begin{defn} \label{defn-Hnorm}
Soit $p_1, \dots, p_\ell $ une base de $I_D$ sur $K$, on pose
$$
\Hnorm(X;D):=  \sum_{v \in M_K} \frac{[K_v:\Q_v]}{[K:\Q]} 
\,  \log |p_1 \wedge \cdots \wedge p_\ell |_v
\enspace.   
$$
\end{defn} 

Cette d{\'e}finition ne d{\'e}pend pas du choix de la base $p_1, \dots,
p_\ell$, 
gr{\^a}ce {\`a} la formule du produit; elle est aussi invariante par extensions finies de $K$. 

La {\it fonction de Hilbert arithm{\'e}tique} introduite mesure, pour chaque valeur $D\in\N$ du param{\`e}tre, la complexit{\'e} 
binaire du $K$-espace lin{\'e}aire des formes de degr{\'e} $D$ de l'anneau
$K[X_0,\dots,X_N]$ modulo $I$.  
Dans le cas des vari{\'e}t{\'e}s toriques,
son comportement asymptotique est li{\'e} {\`a}
la {hauteur normalis{\'e}e} de la vari{\'e}t{\'e} en question
(Proposition~\ref{asymptotique}).
C'est gr{\^a}ce {\`a} ce r{\'e}sultat 
(dont on repousse la d{\'e}monstration au 
paragraphe~\ref{Demonstration du theoreme 1}) 
que la fonction $\Hnorm$ joue un r{\^o}le
crucial dans notre calcul de 
$\wh{h}(X_{\cA,\alpha})$.  
Cependant, il serait tr{\`e}s int{\'e}ressant de montrer que cette 
asymptotique reste valable pour une vari{\'e}t{\'e} projective quelconque.
On se hasarde {\`a} poser~: 

\begin{question} \label{hilbert-normalise} 
A-t-on
\begin{equation} \label{c(X)}
\cH_{\rm norm}(X;D) = \frac{c(X)}{(n+1)!} \, D^{n+1} + o(D^{n+1} ) 
\end{equation}
pour toute sous-vari{\'e}t{\'e} $X\subset\P^N$ 
de dimension $n$ 
avec, de plus, $c(X) = {\wh{h}(X)}$? 
\end{question}

\smallskip 

La fonction $\Hnorm$ admet la formulation duale suivante.
Explicitons d'abord
les contributions locales: {\'e}crivons $p_i= \sum_{a} p_{i\, a} \, x^a$ pour $i=1, \dots, \ell$ 
et posons  $p$ la matrice de taille $\ell \times \binom{D+N}{N}$ form{\'e}e 
par les coefficients de $p_1, \dots, p_\ell$. 
Pour chaque sous-ensemble $L \subset \N_D^{N+1}$ tel que $\Card(L) =
\ell$ 
on pose 
$$
p_L:= (p_{i\, a})_{1 \le i \le \ell \atop a \in L} \in \ K^{\ell \times \ell}
$$
le mineur de $p$ correspondant aux colonnes index{\'e}es par $L$; 
ainsi
$$
p_1 \wedge \cdots \wedge p_\ell = \sum_{ L; \ \card(L)= \ell} \det(p_L)  \, \bigwedge_{a \in L}  x^{a}
$$
et donc $|p_1 \wedge \cdots \wedge p_\ell |_v= \max 
\Big\{  | \det(p_L)|_v \, :  \ L \subset \N_D^{N+1}, \ 
\Card(L)= \ell \Big\}$. 

Soit  $K[x_0, \dots, x_N]_D^\vee$ l'espace 
dual (c.{\`a}.d. l'espace des fonctionnelles lin{\'e}aires 
$K[x_0, \dots, x_N]_D$ $ \to K$) muni de la base duale de celle des mon{\^o}mes. 
Soit $m:= \cH_\geo(X; D)$, pour $v
\in M_K$ et $g \in \bigwedge^m K[x_0,\dots,x_N]_D^\vee$ on note $|g|_v$
la norme sup des coefficients de $g$ pour la place $v$, par
rapport {\`a} la base duale 
$\Big\{ \bigwedge_{b \in L} (x^{b})^\vee \, : \ M \subset \N_D^{N+1},
\ \Card(M)= m  \Big\}$. 

\begin{prop} \label{Hnorm dual}
Soit  $q_1, \dots, q_m \in (K[x_0, \dots,
x_N]_D)^\vee$ 
une base sur $K$ de l'annulateur $\Ann(I_D)$, alors
$$
\Hnorm(X;D)=  \sum_{v \in M_K} \frac{[K_v:\Q_v]}{[K:\Q]} \,
\log |q_1 \wedge \cdots \wedge q_m|_v\enspace.  
$$
\end{prop} 

\begin{demo} 
Les espaces $\bigwedge^\ell K[x_0, \dots, x_N]_D$ et 
$\bigwedge^m K[x_0, \dots, x_N]_D^\vee$ s'identifient 
{\it via} leurs bases respectives~\cite[Ch.~3,~\S~11, Prop.~12]{Bou80}. 
Cette identification est d{\'e}finie par
$$
\Psi_\ell: \bigwedge_{a \in L} x^{a} \mapsto {\varepsilon(L)} \, 
\bigwedge_{b \in L^\ccc} (x^{b})^\vee 
$$
pour chaque sous-ensemble $L \subset \N_D^{N+1}$ 
tel que $\Card(L) = \ell$ et o{\`u}
$L^\cc:= \N_D^{N+1} \setminus L$,  $ {\varepsilon(L)}=\pm 1$. 
Ce dernier signe {\'e}tant irrelevant pour la suite de cette d{\'e}monstration, on renvoie {\`a} \cite[Ch.~3,~\S~11, Formule~(78)]{Bou80} pour sa d{\'e}finition exacte. 

Soit $p_1, \dots, p_\ell$ une base de $I_D$ sur $K$. 
D'apr{\`e}s le th{\'e}or{\`e}me de Brill-Gordan~\cite{Gor1873} 
({\it voir} aussi~\cite[Ch.~3,~\S~11, Prop.~15]{Bou80}) 
il existe $\lambda\in K^\times$ tel que
\begin{equation}\label{BrillGordan}
\Psi_\ell( p_1\wedge\cdots\wedge p_\ell) 
=\lambda \cdot q_1\wedge\dots\wedge q_m\enspace. 
\end{equation}
Ceci implique $\det(p_L) = \pm \lambda \, \det(q_{L^\ccc})$ 
pour tout $L$, d'o{\`u} 
$\log |p_1\wedge\cdots\wedge p_\ell|_v 
= \log|q_1\wedge\dots\wedge q_m|_v + \log|\lambda|_v $, 
et on conclut par sommation sur $v\in M_K$ avec la formule du produit. 
\end{demo} 

\begin{rem} 
Le partie gradu{\'e}e $I_D$  est, par d{\'e}finition,
$$
I_D= \bigg\{ f \, 
: \,  
f(\xi) = \sum_{a\in \N^{N+1}_D} f_a \, \xi^a = 0 \ \ \ \forall  \xi \in X \bigg\}\enspace.  
\vspace{-2mm} 
$$
Autrement-dit, $I_D= \Ann (v_{N,D} (X))$ o{\`u} 
$v_{N,D}: \P^N \to \P^{\binom{D+N}{N}} $ d{\'e}signe le plongement de 
Veronese de degr{\'e} $D$ en $N+1$ variables homog{\`e}nes.
Ainsi $\Ann(I_D)$ co{\"\i}ncide avec l'espace lin{\'e}aire engendr{\'e} par
l'ensemble
$v_{N,D} (X)$.

Notons encore  que  $\Hnorm$ 
se comporte bien par rapport aux plongements de Veronese~: 
on a $ \displaystyle 
{\cH}_{\rm norm}(v_{N,D}(X);D') = {\cH}_{\rm norm}(X; D\, D')$
pour $D' \in \N$. 

\end{rem} 

\smallskip 

On introduit maintenant une deuxi{\`e}me fonction de Hilbert arithm{\'e}tique, 
dont on aura besoin pour la d{\'e}monstration de la 
proposition~\ref{asymptotique}. 
Pour $f = \sum_{a} f_a \, x^a, \, 
g= \sum_{a} g_a \, x^a \in \C[x_0, \dots, x_N]_D$ on  consid{\`e}re le
produit scalaire 
$$
\langle f , g \rangle = \langle f , g \rangle_D:= \sum_{ a} \binom{D}{a}^{-1} \,  f_a \,\overline{g}_a\enspace. 
$$
{\`A} un facteur pr{\`e}s, ce produit co{\"\i}ncide avec la norme $\cL^2$ 
sur $\P^N(\C)$ par rapport {\`a} la m{\'e}trique de Fubini-Study: on a
\begin{equation}\label{produit L2}
\langle f , g \rangle= \binom{D+N}{N} \,  
\int_{\P^N(\C)} \frac{f(z) \, \overline{g(z)}}{\hspace{4mm}
||z||^{2\,D}} \, \omega_{\FS}^N
\enspace, 
\end{equation}
o{\`u} $\omega_{\FS}$ d{\'e}note la $(1,1)$-forme de Fubini-Study, 
normalis{\'e}e de sorte que l'espace ambiant $\P^N(\C)$ soit de masse totale $1$, {\it voir} par exemple 
\cite[\S~1.3.1]{Ran01}. 
La puissance ext{\'e}rieure $\bigwedge^\ell \C[x_0, \dots, x_N]_D$ est alors munie du produit scalaire 
\begin{equation*} 
\langle f_1 \wedge \cdots \wedge f_\ell , g_1 \wedge \cdots \wedge g_\ell
\rangle := \det \Big( \langle f_i, g_j \rangle \Big)_{1 \le i,j \le
  \ell} 
\enspace . 
\end{equation*} 
Soit $v \in M_K^\infty$ et $\sigma_v: K \hookrightarrow \C$ un plongement correspondant {\`a} $v$, alors pour $p_1, \dots, p_\ell \in K[x_0, \dots, x_N]_D$ on pose 
$$
||p_1 \wedge \cdots \wedge p_\ell||_v := ||\sigma_v(p_1)\wedge\cdots \wedge\sigma_v(p_\ell)||
$$ 
o{\`u} $|| \cdot ||$ est la norme associ{\'e}e au produit 
scalaire $\langle \cdot , \cdot \rangle$. 
Ainsi on d{\'e}finit la fonction 
\begin{eqnarray}\label{def Hproj}  
\qquad\Harith(X;D) &:= &\sum_{v\in
  M_K^\infty}\frac{[K_v:\Q_v]}{[K:\Q]}
\, \log||p_1 \wedge \cdots \wedge p_\ell||_v \\[-1mm]
&&+\sum_{v\in M_K\setminus M_K^\infty}\frac{[K_v:\Q_v]}{[K:\Q]} \, 
\log|p_1 \wedge \cdots \wedge p_\ell |_v \nonumber 
+ \frac{1}{2}\log (\gamma(N,D))\enspace, 
\end{eqnarray} 
o{\`u}  $p_1, \dots, p_\ell $ est une base de $I_D$ sur $K$ et 
$\gamma(N,D):= \prod_{a \in \N_D^{N+1}}{{D}\choose a}$. 
Comme pour $\Hnorm(X;D)$, cette d{\'e}finition  est ind{\'e}pendante du 
choix de la base et invariante par extensions finies du corps $K$.  

\smallskip 

Cette fonction est proche d'autres analogues arithm{\'e}tiques de la 
fonction de Hilbert, d{\'e}finis eux aussi comme une hauteur du
sous-espace lin{\'e}aire $I_D$. 
En particulier elle co{\"\i}ncide, {\`a} un d{\'e}calage explicite pr{\`e}s, avec  
la fonction de Hilbert arithm{\'e}tique propos{\'e}e par M.~Laurent pour 
les $\Q$-sous-vari{\'e}t{\'e}s dans \cite[pp.~224-226]{Lau92}, 
qu'on notera ici $\cH_\laurent$; on a 
$$
\cH_\ari(X ; D)= \cH_\laurent(X ; \, D)  + \frac{1}{2} \, {D+N
  \choose N}^{-1} \log(\gamma(N,D)) \cdot \cH_\geo(X;  D) \enspace.
$$
On remarquera que d'apr{\`e}s~\cite[Prop.~3]{Lau92} on a
$${{D+N}\choose N}^{-1}\log(\gamma(N,D)) 
= D \, \sum_{j=1}^N\frac{1}{j+1} -\frac{N}{2} \,  \log(D) 
+\frac{N+1}{2}\, \sum_{j=1}^N \frac{1}{j+1} - \frac{N}{2}\, \log(2\pi) +o_D(1)
\ .$$
En outre, $\Harith$ peut aussi se comparer avec la fonction de Hilbert
arithm{\'e}tique introduite par S.~David et P.~Philippon
dans~\cite[\S~4.3]{DP99}, qu'on notera ici $\cH_\daph$; on a 
$$
\cH_\daph(X ; D) \leq \Harith(X ; D) \leq \cH_\daph(X ; D) + c\, \frac{N+1}{2} {D+N \choose N}^{-1}\kern-3pt \log(\gamma(N,D)) \cdot \cH_\geo(X; D)   
$$
o{\`u} $c>0$ est une constante universelle telle que 
$ \displaystyle 
\pi(x):=\Card\Big( \{ p \ \mbox{premier} \, : \ p \le x \} \Big) 
\displaystyle 
\leq c  \cdot {x}/{\log(x)}$.
D'apr{\`e}s~\cite[Cor.~1]{RS62} on peut prendre $c=1,26$. 
La premi{\`e}re estimation est directe, tandis que pour la deuxi{\`e}me il
faut utiliser des formulations duales pour $\Harith$ 
(Proposition~\ref{expressionhilbertarith} ci-dessous) et $\cH_\daph$. 

\smallskip 

Comme pour $\Hnorm$,  la fonction $\Harith$ admet une
formulation duale. 
Le produit scalaire induit un isomorphisme lin{\'e}aire  
$$
\eta: \C[x_0, \dots, x_N]_D \to \C[x_0, \dots, x_N]_D^\vee 
\enspace , \quad \quad f \mapsto \langle \cdot, f \rangle\enspace, 
$$ 
et cette identification munit $\C[x_0, \dots, x_N]_D^\vee$ du 
produit scalaire {\it dual}, d{\'e}fini par 
$$
\langle  \eta(f), \eta(g) \rangle^\vee := \langle f, g \rangle 
$$
pour $f,g \in  \C[x_0, \dots, x_N]_D$. 
Pour $\theta = \sum_{ b } \theta_b \, (x^b)^\vee, \, 
\zeta= \sum_{b} \zeta_b \, (x^b)^\vee  
\in \C[x_0, \dots, x_N]_D^\vee$ on v{\'e}rifie ais{\'e}ment que ce produit s'{\'e}crit
$$
\langle  \theta, \zeta \rangle^\vee := \sum_{b} {D \choose
  b} \,  
\theta_b \,\overline{\zeta}_b
\enspace.$$
Ceci s'{\'e}tend de mani{\`e}re naturelle aux puissances ext{\'e}rieures
$\bigwedge^m (\C[x_0, \dots, x_N]_D)^\vee $ par
$$
\langle \theta_1 \wedge \cdots \wedge \theta_m , 
\zeta_1 \wedge \cdots \wedge \zeta_m
\rangle := \det \Big( \langle \theta_i, \zeta_j \rangle \Big)_{1 \le i,j \le
  m} 
\enspace . 
$$

\begin{prop}\label{expressionhilbertarith}
Soit $q_1, \dots, q_m \in K[x_0, \dots, x_N]_D^\vee$ une base de $\Ann(I_D)$ sur $K$, alors
$$
\Harith(X;D)= \sum_{v \in M_K^\infty} \frac{[K_v:\Q_v]}{[K:\Q]} \, 
\log ||q_1 \wedge \cdots \wedge q_m||_v^\vee \ +\kern-3pt \sum_{v \in
  M_K \setminus M_K^\infty} \frac{[K_v:\Q_v]}{[K:\Q]} \, 
\log |q_1 \wedge \cdots \wedge q_m |_v\enspace.  
$$
\end{prop}

Ce r{\'e}sultat {\'e}tend \cite[Thm.~3]{Lau92} au cas des sous-vari{\'e}t{\'e}s de
dimension positive; cependant la d{\'e}monstration reste (pour l'essentiel) inchang{\'e}e: 

\begin{demo} 
Soit $p_1 \dots, p_\ell$ une base de $I_D$, gr{\^a}ce au th{\'e}or{\`e}me de 
Brill-Gordan (Identit{\'e}~(\ref{BrillGordan}) ci-dessus) on a $\det(p_L)
= \pm \lambda \, \det(q_{L^\cc})$ pour tout sous-ensemble $L \subset
\N_D^{N+1}$ tel que $\Card(L) = \ell$ et o{\`u} $L^\cc:= {\N_D^{N+1} \setminus L}$. Donc pour $v \in M_K^\infty$  
$$
\begin{array}{rl}
\log ||p_1\wedge&\kern-8pt\cdots\wedge p_\ell||_v  \displaystyle = \frac{1}{2} \, \log \left(\sum_{L; \ \card(L) = \ell} \left( \prod_{a\in L} {D \choose a }^{-1} \right) \, | \det (p_L) |_v^2 \right) \\[4mm]
&\displaystyle = \frac{1}{2} \, \log \left(\sum_{L; \ \card(L) = \ell} \left( \prod_{a\in L} {D \choose a }^{-1} \right) \, | \det (q_{L^{\cc}}) |_v^2 \right) + \log |\lambda|_v \\[6mm]
&\displaystyle = \frac{1}{2} \, \log \left(\sum_{ M; \ \card(M) = m}\left( \prod_{b\in M} {D \choose b } \right) \, | \det (q_{M}) |_v^2 \right)  - \frac{1}{2} \, \log \Big( \gamma (N,D) \Big) + \log|\lambda|_v\\[6mm]
&\displaystyle = \log||q_1\wedge\dots\wedge q_M||_v^\vee-  \frac{1}{2}  \, \log \Big( \gamma (N,D) \Big) + \log|\lambda|_v  
\end{array}$$
tandis que  pour $v \in M_K \setminus M_K^\infty$ on a
$\log|p_1\wedge\cdots \wedge p_\ell|_v = \log|q_1\wedge\dots\wedge
q_m|_v + \log|\lambda|_v$. On conclut  par sommation sur $v\in M_K$ avec
la formule du produit.  
\end{demo} 

\medskip

On explicite maintenant quelques uns des r{\'e}sultats de la th{\`e}se de 
H.~Randriam\-bo\-lo\-lo\-na~\cite{Ran01}
qui nous permettrons d'expliciter le comportement asymptotique de
$\Harith$. 
Pour le reste de cette section, on utilisera librement le langage de la g{\'e}om{\'e}trie d'Arakelov, on renvoie {\`a} \cite{BGS94} et \cite[Partie I]{Ran01} pour les d{\'e}finitions et propri{\'e}t{\'e}s de base des objets.

Soit $X_K$ le sous-sch{\'e}ma ferm{\'e} int{\`e}gre de $\P^N_K$ d{\'e}fini par l'id{\'e}al 
$I_{K} := I \cap K[x_0, \dots, x_N]$ et soit $\Sigma$ la cl{\^o}ture de 
Zariski de $X_K$ dans $\P^N_{\cO_K}$. 
Autrement-dit, $\Sigma$ est le sch{\'e}ma projectif  d{\'e}fini par l'id{\'e}al
premier homog{\`e}ne 
$$
I_{\cO_K} := I \cap \cO_K[x_0, \dots, x_N]\subset \OK[x_0, \dots, x_N]
\enspace.
$$
Par construction, $\Sigma$ est un sous-sch{\'e}ma de $\P^N_{\cO_K}$ ferm{\'e}, int{\`e}gre, de dimension $n+1$ et plat sur $\Spec(\cO_K)$ car $I_{\cO_K} \cap \cO_K = \{0\}$. 

Soit ${\cO(D)}:= {\cO(1)}^{\otimes D}$ la $D$-i{\`e}me puissance
tensorielle du fibr{\'e} en droites universel ${\cO(1)}$ sur $\P^N_\OK$,
on pose 
$M:= \Gamma(\Sigma, \cO(D)|_\Sigma)$ le $\cO_K$-module des sections globales de $\cO(D)|_\Sigma$. 
Pour chaque plongement $\sigma: K \hookrightarrow \C$ on pose $\Gamma(\Sigma, \cO(D)|_{\Sigma})_\sigma := \Gamma(\Sigma, \cO(D)|_{\Sigma}) \otimes_\sigma \C$ et on consid{\`e}re l'application de restriction 
$$
\Gamma(\P^N_{\cO_K} , \cO(D))_\sigma \to \Gamma(\Sigma, \cO(D)|_{\Sigma})_\sigma\enspace.  
$$
L'espace $\Gamma( \P^N_{\cO_K} , \cO(D))_\sigma$ s'identifie
canoniquement {\`a} la partie gradu{\'e}e  $K_\sigma[x_0, \dots, x_N]_D$, on
le munit du produit $\cL^2$ par rapport {\`a} la m{\'e}trique de
Fubini-Study. Ce produit est un multiple du produit $\langle \cdot,
\cdot \rangle$ introduit pr{\'e}c{\'e}demment; on a~(Formule~(\ref{produit
  L2})\,)
$$
\langle f , g \rangle_{\cL^2 , \sigma} = {D+N \choose N}^{-1} \, \langle \sigma(f) , \sigma(g) \rangle 
$$ 
pour tout $f,g \in \Gamma(\P^N_{\cO_K} , \cO(D))_\sigma$. 
L'application de restriction est surjective pour $d \gg 0$; dans cette
situation on munit  $M_\sigma$ du produit scalaire {\it quotient}, c.{\`a}.d. on identifie cet espace avec l'orthogonal du noyau de l'application de restriction et on l'{\'e}quipe du produit induit. 

Cette construction munit $M$ d'une structure de $\OK$-module {\it hermitien} dont le degr{\'e} arith\-m{\'e}\-ti\-que est par d{\'e}finition
\begin{equation}\label{defdegarith}
\quad \wh{\deg}(\ov{M}) := \frac{1}{[K:\Q]} \, \Biggl( \log\mbox{\rm
  Card}\bigg(\bigwedge^s M /(f_1 \wedge  \cdots \wedge f_s)\bigg) 
- \sum_{\sigma: K \hookrightarrow \C} \log\Vert {f_1}\wedge\dots \wedge {f_s}\Vert_{\cL^2,\sigma}\Biggr) 
\end{equation}
o{\`u}  ${f_1},\dots,{f_s} \in M$ est une base de $\Gamma(X_K;\,
\cO(D)|_{X_K})$ sur $K$ et $\Vert\cdot\Vert_{\cL^2,\sigma}$ est la
norme associ{\'e}e au produit $\langle\cdot,\cdot\rangle_{\cL^2,\sigma}$.   
Dans le langage de~\cite{Ran01} ce degr{\'e} arithm{\'e}tique est la {\em $D$-i{\`e}me hauteur
  sch{\'e}matique 
de $\Sigma$}; 
on montre qu'il co{\"\i}ncide ({\`a} un d{\'e}calage pr{\`e}s) avec $\Harith(X; \, D)$ pour $D$ assez grand: 

\begin{lem} \label{Hnorm=degarith}
Avec les notations ci-dessus, il existe $D_0 \in \N$ tel que pour $D \ge D_0$ on ait 
$$
\cH_\ari(X;D) = \wh{\deg}\Big(\ov{\Gamma(\Sigma, \cO(D)|_{\Sigma})}\Big) 
- \frac{1}{2} \,\cH_\geo(X; D)\, \log{{D+N}\choose N}\enspace. 
$$
\end{lem} 
\begin{demo} 
Soit $\cI$ le faisceau d'id{\'e}aux de d{\'e}finition de $\Sigma$ et $\Gamma( \P^N_{\cO_K} , \cI \, \cO(D))$ le $\OK$-module des sections globales de $\cI \, \cO(D)$, muni des produits scalaires obtenus par restriction du produit $\cL^2$. 
Soit $M=\Gamma(\Sigma,\cO(D)|_{\Sigma})$, alors pour $D_0 \gg 0$ assez grand  et $D\geq D_0$
on a~\cite[Lem.~2.3.6]{Ran01}
\begin{equation} \label{suite-exacte} 
\wh{\deg}(\ov{M}) = \wh{\deg} \Big(\ov{\Gamma(\P^N_{\cO_K} , \cO(D))}\Big) 
-\wh{\deg} \Big(\ov{ \Gamma(\P^N_{\cO_K} , \cI \, \cO(D))}\Big) 
\end{equation} 
{\`a} cause de l'exactitude de la suite de $\OK$-modules hermitiens 
$$0 \to \ov{ \Gamma(\P^N_{\cO_K} , \cI \, \cO(D))} \to \ov{\Gamma(\P^N_{\cO_K} , \cO(D))} \to M \to 0\enspace,$$
qui r{\'e}sulte pour $D\geq D_0$ de l'exactitude de la suite des faisceaux m{\'e}tris{\'e}s correspondants.

En prenant la base des sections associ{\'e}e {\`a} la base des mon{\^o}mes, on v{\'e}rifie ais{\'e}ment  
$$
\wh{\deg} (\ov{\Gamma(\P^N_{\cO_K} , \cO(D))}) 
= \frac{1}{2} \, \log(\gamma(N,D)) + \frac{1}{2} \, 
{D+N \choose N} \, \log{D+N \choose N}\enspace.
$$
En outre, $\Gamma(\P^N_{\cO_K} , \cI \, \cO(D))$ s'identifie de mani{\`e}re naturelle {\`a} 
la partie gradu{\'e}e $\displaystyle (I_{\cO_K})_D \subset \OK[x_0, \dots, x_N]_D$. 
Soit $p_1, \dots, p_\ell \in (I_{\cO_K})_D$  une base de
$I_D$
sur $K$; de l'identit{\'e}~(\ref{suite-exacte}), de la d{\'e}finition~(\ref{defdegarith}) appliqu{\'e}e {\`a} $(I_{\cO_K})_D$ et du rapport entre les diff{\'e}rents produits scalaires on d{\'e}duit pour $D \ge D_0$ 
\begin{eqnarray*}
\wh{\deg}(\ov{M})&=&\frac{1}{2} \, \log(\gamma(N,D) ) 
+ \frac{1}{2} \, {D+N \choose N}  \, \log{D+N \choose N}\\[-1mm] 
&&+ \frac{1}{[K:\Q]} \, \Bigg( \sum_{\sigma: K \hookrightarrow \C} 
\log||p_1 \wedge \cdots \wedge p_\ell||_{\cL^2, \sigma} - \log \Card \Big( \bigwedge^\ell (I_\OK)_D / (p_1 \wedge  \cdots \wedge p_\ell)\Big) \Bigg) 
\end{eqnarray*}
\begin{eqnarray*}
  \hspace*{3mm} &= &\frac{1}{2} \, \log(\gamma(N,D)) + \frac{1}{2} \, \cH_\geo(X;  D)\, \log{D+N \choose N}  \\[-1mm] 
&&+\sum_{v \in M_K^\infty} \frac{[K_v:\Q_v]}{[K:\Q]} \, \log||p_1 \wedge
\cdots \wedge p_\ell||_{v} - \frac{1}{[K:\Q]} \, \log \Card \Big(\bigwedge^\ell (I_\OK)_D    /(p_1 \wedge  \cdots \wedge p_\ell)\Big)
\enspace. 
\end{eqnarray*}

On pose maintenant  $J:= \bigwedge^\ell (I_{\cO_K})_D$ et $j:= p_1
\wedge \cdots \wedge p_\ell \in J$; vu la d{\'e}finition de $\Harith$ on est  r{\'e}duit {\`a} d{\'e}montrer l'identit{\'e}
$$
\frac{1}{[K:\Q]} \, \Card \Big(J/j\OK) = - \hspace{-3mm} \sum_{v \in M_K \setminus M_K^\infty} 
\frac{[K_v:\Q_v]}{[K:\Q]} \, \log | j|_v\enspace. 
$$
On consid{\`e}re l'application $ J/j\OK  \to \prod_{\gp} ( J/j\OK)_\gp$
pour $\gp $ parcourant l'ensemble des id{\'e}aux premiers de $\OK$, c'est
un {isomorphisme} 
car $J$ est un $\OK$-module projectif de rang 1.
Soit $\gp$ un id{\'e}al premier de $\OK$,  le module localis{\'e} $J_\gp$ est
donc  {\it libre} de 
rang 1. 
Fixons un g{\'e}n{\'e}rateur $f_\gp $ et soit $j_\gp \in (\OK)_\gp$ tel que $j
= j_\gp \, f_\gp$, alors on a un isomorphisme
$$
(J/j\OK)_\gp \cong (\OK)_\gp/ j_\gp(\OK)_\gp\enspace. 
$$
Soit $\Ord_\gp$ la valuation et $|\cdot|_v:=
\mbox{Norme}_{K/\Q}(\gp)^{\ord_\gp(\cdot)} $
la valeur absolue ultram{\'e}trique associ{\'e}es {\`a} l'id{\'e}al $\gp$, on remarque que $|f_\gp|_{v} = 1$ car $(I_{\cO_K})_D$ est un sous-module satur{\'e} de $\OK[x_0, \dots, x_N]_D$, donc $|j_\gp|_{v} =|j|_{v}$. En passant au cardinaux on trouve 
$$
\Card \Big( ( J/j\OK)_\gp \Big)  = \Card \Big( (\OK)_\gp/ \gp^{\ord_\gp(j_\gp)} (\OK)_\gp \Big)  = \mbox{Norme}_{K/\Q}(\gp)^{\ord_\gp(j_\gp)} = |j|_{v}^{-[K_v : \Q_v]}\enspace,  
$$
ce qui entra{\^\i}ne l'identit{\'e} cherch{\'e}e. 
\end{demo} 

\medskip 

Maintenant, \cite[Thm.~A.2.2]{Ran01} 
({\it voir}
aussi~\cite{Ran03})
s'{\'e}nonce 
$ \displaystyle 
\wh{\deg}\Big(\ov{\Gamma(\Sigma, \cO(D)|_{\Sigma})}\Big) 
= \frac{h(X)}{(n+1)!} \, D^{n+1} + o(D^{n+1} )$
o{\`u} $h(X)$ est la hauteur projective de $X$. 
Avec le lemme pr{\'e}c{\'e}dent,  ceci implique 
\begin{equation} \label{asym Hproj}
\cH_\ari(X;\,D)= \frac{h(X)}{(n+1)!} \, D^{n+1} + o(D^{n+1} )
\end{equation} 
lorsque $D$ tend vers l'infini, car $\frac{1}{2}\, \cH_\geo(X; \, D)\, \log{{D+N}\choose N} = o(D^{n+1})$.

\begin{cor}\label{expressionhilbertnorm}
Soit  $X \subset \P^N$ une sous-vari{\'e}t{\'e} de dimension $n$, alors 
\begin{eqnarray*} 
(n+1)!\ \un{\lim}_{D\rightarrow\infty}\frac{1}{D^{n+1}}\, \Hnorm
(X;D)& 
\geq &h(X) - \frac{1}{2} \, (n+1)\log(N+1)\, \deg(X) \enspace,  \\[0mm] 
(n+1)!\ \ov{\lim}_{D\rightarrow\infty}\frac{1}{D^{n+1}}\, \Hnorm
(X;D)  & \leq & h(X) \enspace.
\end{eqnarray*} 
\end{cor}

\begin{demo} 
Soit $q_1, \dots, q_m \in K[x_0, \dots, x_N]_D^\vee$ une base de $\Ann(I_D)$ sur $K$. Pour $v \in M_K^\infty$ on a
\begin{eqnarray*} 
|q_1\wedge \cdots \wedge q_m|_v & = &  \max \Big\{ |\det(q_M)|_v \, :
 \ M \subset \N_D^{N+1}, \ \Card(M) = m \Big\}\\[0mm]  
&\le & \left(\sum_{ M; \ \card(M) = m} 
\left( \prod_{b\in M} {D \choose b } \right) \, | \det (q_{M}) |_v^2
 \right)^{1/2} =  ||q_1\wedge\dots\wedge q_m||_v^\vee
\enspace,    
\end{eqnarray*} 
et dans l'autre direction 
$$
||q_1\wedge\dots\wedge q_m||_v^\vee \le \left(\sum_{ M; \ \card(M) = m} 
 \ \prod_{b\in M} {D \choose b } \right)^{1/2}\, |q_1\wedge \cdots \wedge q_m|_v
\le (N+1)^{\frac{m\, D}{2}} \, |q_1\wedge \cdots \wedge q_m|_v
\enspace.$$
Par sommation sur $v\in M_K$ et gr{\^a}ce aux 
propositions~\ref{Hnorm dual} et~\ref{expressionhilbertarith} on
obtient l'encadrement 
$$
\Harith(X;D)- \frac{1}{2} \, \log(N+1)\ D\,  \cH_\geo(X;D) 
\ \leq \ \cH_{\rm norm}(X;D) \ \leq \ \cH_\ari(X;D)\enspace. 
$$
L'{\'e}nonc{\'e} d{\'e}coule de ces estimations et de l'asymptotique~(\ref{asym Hproj})
pour $\Harith$. 
\end{demo} 

\smallskip

Dans le cas o{\`u} ${\cH}_{\rm norm}(X;D)$ admet un comportement
asymptotique du type (\ref{c(X)}), il suit de ce corollaire
l'encadrement 
\begin{equation}\label{encadrementhilbert}
c(X) \le h(X) \le c(X) + \frac{n+1}{2} \, \log(N+1)\, \deg(X)\ . 
\end{equation} 
En particulier, lorsque la question \ref{hilbert-normalise} a une r{\'e}ponse positive, cela am{\'e}liore \cite[Prop. 2.1 (v)]{DP99} d'un facteur $2/7$.



%
%


\typeout{Hauteur et fonction de Hilbert arithmetique des varietes toriques} 

\section{Hauteur et fonction de Hilbert arithm{\'e}tique des vari{\'e}t{\'e}s toriques} 

\label{normalisee} 

\setcounter{equation}{0}
\renewcommand{\theequation}{\thesection.\arabic{equation}}


\typeout{Poids de Hilbert et de Chow} 

\subsection{Poids  de Hilbert et de Chow} 

\label{poids}

Dans ce paragraphe on se place, comme au~\S~\ref{monomiales}, 
sur un corps de base $\K$. 
Soit $X \subset \P^N=\P^N(\K)$ une sous-vari{\'e}t{\'e} de dimension $n$
et $\tau= (\tau_0, \dots, \tau_N) \in \R^{N+1}$ un vecteur {\it poids}. 
Consid{\'e}rons la forme de Chow $\Ch_X \in \K[U_0, \dots, U_n]$ 
de $X$ et une variable additionnelle $t$; 
si l'on {\'e}crit 
$$
\Ch_X\Big( t^{\tau_j} \, U_{i\,j} \, : \ 0\le i \le n, \ 0 \le j \le N \Big) = t^{e_0} \, F_0 + \cdots +t^{e_M} \, F_M
$$
avec $F_0, \dots, F_M \in \K[U_0, \dots, U_n] \setminus \{ 0\}$ 
et $e_0 > \cdots > e_M$, le {\em $\tau$-poids de Chow de $X$} 
est d{\'e}fini par $e_{\tau}(X):= e_0$. 
Cette notion est introduite dans~\cite[p.~61]{Mum77} dans le contexte
de la th{\'e}orie g{\'e}om{\'e}trique des invariants, en relation avec
l'{\'e}tude de la 
stabilit{\'e} des vari{\'e}t{\'e}s projectives.

\smallskip 

Posons 
$$
\ell_\tau(x):= \tau_0\, x_0 + \cdots + \tau_N \, x_N  \in \R[x_0, \dots, x_N] 
$$
la forme lin{\'e}aire lin{\'e}aire associ{\'e}e {\`a} $\tau$. 
Pour $D\in \N$ et $I= I(X) \subset \K[x_0,\dots, x_N]$ 
l'id{\'e}al de d{\'e}finition de $X$, 
le {\em $\tau$-poids de Hilbert de $X$ en degr{\'e} $D$} (ou
{\it $\tau$-poids de $\Big(\K[x_0, \dots, x_N]/I \Big)_D$} 
selon la terminologie de~\cite{Mum77}) est 
$$
s_\tau(X;D):= \max_J \quad \sum_{\lambda \in J} \ell_\tau(\lambda) 
$$
o{\`u} le maximum est pris sur tous les 
$J \subset \N^{N+1}_D$
tels que l'ensemble des mon{\^o}mes associ{\'e}s 
$\{ x^\lambda \, :\ \lambda \in J\}$ induise une base de 
la partie gradu{\'e}e $\Big(\K[x_0, \dots, x_N]/I \Big)_D$. 
Cette quantit{\'e} est aussi appel{\'e}e {\em $D$-i{\`e}me poids de 
Hilbert de $I$ par rapport {\`a} $\tau$} dans~\cite[\S~4.1]{EF02}. 
Un r{\'e}sultat de Mumford~\cite[Prop.~2.11]{Mum77} montre que 
l'asymptotique de cette fonction pour $D$ tendant vers l'infini 
est 
\begin{equation}\label{asymptopoids}
s_\tau(X;D)= \frac{e_{\tau}(X)}{(n+1)!} \, D^{n+1} + O(D^n)
\enspace.
\end{equation}

\smallskip

En principe, on sait calculer les poids de Chow d'une vari{\'e}t{\'e} 
donn{\'e}e puisqu'on peut calculer sa forme de Chow, pourtant 
cela reste difficile en pratique. 
Dans la suite on donne une expression simple des poids
de Hilbert et de Chow dans le cas des 
vari{\'e}t{\'e}s toriques. 

Pour un vecteur  $\cA= (a_0, \dots, a_N) \in (\Z^n)^{N+1}$ 
et un poids $\tau \in \R^{N+1}$ on consid{\`e}re le polytope
$$
Q_{\cA,\tau}:= \Conv \Big((a_0,\tau_0), \dots, (a_N,\tau_N) \Big) \ \subset \R^{n+1}\enspace, 
$$
dont l'enveloppe sup{\'e}rieure s'envoie bijectivement sur $Q_\cA$ par la projection standard $\R^{n+1} \to \R^n$. Soit alors 
$$
\vartheta_{\cA,\tau}  : Q_\cA \to \R \quad , \qquad x \mapsto 
\max \Big\{ y \in \R \, : \ (x,y) \in Q_{\cA,\tau} \Big\}\enspace, 
$$ 
la param{\'e}trisation de cet enveloppe sup{\'e}rieure
au dessus de $Q_\cA$.  
C'est une fonction  {\it concave} et  
{\it affine par morceaux}, donc en particulier Riemann int{\'e}grable
sur tout polytope. 

Consid{\'e}rons l'application lin{\'e}aire  
$$
M_\cA:= \Z^{N+1} \to \Z^n
\enspace, \quad \quad 
\lambda\mapsto 
\lambda_0
\, a_0 + \cdots + \lambda_N \, a_N
\enspace.
$$ 
On pose 
$ \cA(D): = M_\cA
\Big( \N_D^{N+1} \Big) \subset \Z^n $
l'image de $\N^{N+1}_D$ par $M_\cA$; 
ainsi $\cA(1)= \Expo(\cA)$.

\begin{prop} \label{poid-torique}
Soit $\cA \in ( \Z^n)^{N+1}$ tel que $L_\cA= \Z^n$, $\tau \in \R^{N+1}$ et $D \in \N$, alors 
\begin{eqnarray*} 
s_{\tau}(X_\cA;D) &= &\sum_{c \in \cA(D) } 
\max \Big\{ \ell_\tau (\lambda)  \, : \  \lambda \in \N^{N+1}_D, \ M_\cA(\lambda) =c \Big\}\enspace, \\[0mm] 
e_{\tau}(X_\cA) & = & (n+1)! \, \int_{Q_\cA} \, \vartheta_{\cA,\tau} \ dx_1 \cdots  dx_n\enspace.
\end{eqnarray*}  
\end{prop}

\begin{demo}
L'application $\varphi_\cA:\T^n=(\K^\times)^n \to X_\cA^\circ $ 
est un isomorphisme {\`a} cause de l'hypoth{\`e}se $L_\cA=\Z^n$. 
Au niveau des alg{\`e}bres cela se traduit par un isomorphisme
$$
(\varphi_\cA)^*: \K[x_0, \dots, x_N] /I(X_\cA) \to \K[s_1^{\pm 1}, \dots, s_n^{\pm 1} ] \quad , \quad \quad x_i \mapsto s^{a_i} 
$$
mettant en correspondance les mon{\^o}mes de $\K[x_0, \dots, x_N]$ avec
ceux de $\K[s_1^{\pm 1}, \dots, s_n^{\pm 1} ] $. En identifiant ces
ensembles avec $\N^{N+1}$ et $\Z^n$ respectivement, on v{\'e}rifie que
cette correspondance 
se traduit en l'application lin{\'e}aire 
$M_\cA: \N^{N+1} \to \Z^n$.
Donc, un ensemble $J \subset \N^{N+1}_D$  correspond {\`a} une base monomiale de  
$ \Big( \K[x_0, \dots, x_N] /I(X_\cA) \Big)_D$ si et seulement si  
$ M_\cA: J \to  \cA(D) $ est une {\it bijection}, ce qui d{\'e}montre la premi{\`e}re formule de la proposition.

\smallskip 

Consid{\'e}rons maintenant la fonction auxiliaire
$$
\rho(D) := \sum_{\gamma \in Q \cap 
\frac{1}{D} \Z^n} \vartheta_{\cA, \tau} (\gamma) \enspace,
$$
c.{\`a}.d.  $D^n$ fois la somme de Riemann de la fonction  $\vartheta_{\cA, \tau}$ sur le polytope $Q:= Q_\cA$. On va comparer le poids de Hilbert $s_{\tau}(X_\cA;D)$ {\`a} cette fonction, ce qui donnera l'asymptotique, puis l'expression int{\'e}grale cherch{\'e}e pour le poids de Chow.

On suppose s.p.d.g. $\tau \in \R_+^{N+1}$. 
Alors, pour $c \in \cA(D) $ 
on a $ \displaystyle 
\gamma:= \frac{c}{D} \in Q \cap \frac{1}{D} \, \Z^n$ 
et pour tout $\lambda \in \N^{N+1}_D $ tel que 
$M_\cA(\lambda)= c$ on a $\displaystyle \ell_\tau(\lambda)  \le D\, \vartheta_{\cA, \tau}(\gamma)$ car
$$
\vartheta_{\cA, \tau}(\gamma) = \max \Big\{ \ell_\tau(t) \, : 
\ t \in \R_+^{N+1}, \ t_0+\dots+t_N = 1, \ M_\cA(t) = \gamma  \Big\}\enspace. 
$$ 
De plus, $\vartheta_{\cA, \tau}(\gamma) \ge 0 $ pour tout $\gamma \in Q \cap \frac{1}{D} \, \Z^n$ {\`a} cause de l'hypoth{\`e}se $\tau \in \R_+^{N+1}$, donc $s_{\tau}(X_\cA;D) \le \ D\, \rho(D)$. On en d{\'e}duit, en appliquant l'asymptotique~(\ref{asymptopoids}), 
$$
\frac{e_{\tau}(X_\cA)}{(n+1)!} \quad = \quad \lim_{D\to \infty} 
\frac{s_{\tau}(X_\cA;D)}{D^{n+1}} \quad \le \quad \lim_{D\to \infty} 
\frac{\rho( D)}{D^{n}} \quad = \quad \int_Q \vartheta_{\cA,\tau} 
\, dx_1 \cdots  dx_n\enspace. 
$$

\smallskip 

Pour l'autre direction on suppose, {\`a} nouveau s.p.d.g., $\tau_0 \le \tau_i$ pour tout $i$ et $a_0 = \tau_0 =0$. D'abord on va  montrer qu'il existe  $D_0\in \N$ et $v_0 \in \Z^n$ tels que pour tout $D$ on ait
$$ 
D\, Q   \cap \Z^n \subset \cA(D+D_0) - v_0\enspace. 
$$

Soit $\Lambda_\cA$ un domaine fondamental quelconque du 
r{\'e}seau  $\Ker(M_\cA) \cap \Z^{N+1}$ de l'espace lin{\'e}aire
$\Ker(M_\cA)$; on prend aussi $w_0 \in \Z^{N+1}$ tel que $\Lambda_\cA + w_0 \subset \R_+^{N+1}$. Soit $c  \in D\, Q   \cap \Z^n$, qu'on {\'e}crira 
$$
c  = \mu_0 \, a_0 + \cdots + \mu_N \, a_N  
= D \, \Big( t_0 \, a_0 + \cdots + t_N \, a_N \Big)  
$$
avec $\mu:= (\mu_0 , \dots, \mu_N)  \in \Z^{N+1}$ 
et $ t:= (t_0, \dots, t_N) \in \R_+^{N+1}$ tel que $t_0+\dots+t_N =
1$; 
ainsi $\mu - D\, t \in \Ker(M_\cA)$. 
Soit  $\nu\in \Ker(M_\cA) \cap \Z^{N+1}$ 
tel que $\mu - D \, t - \nu \in \Lambda_\cA$, 
on pose alors $\lambda:=  \mu - \nu +w_0 \in \Z^{N+1}$. On a $\lambda   \in \Lambda_\cA + w_0 + D\, t  \subset \R_+^{N+1}$, d'o{\`u} on d{\'e}duit facilement $\lambda \in \N^{N+1}_{D+D_1}$ avec $D_1\le D_0:=\max\{  \xi_0+\dots+\xi_N \, : \ \xi \in \Lambda_\cA + w_0\}$. Finalement, on pose $v_0: = M_\cA(w_0) \in \Z^n$ et on a bien  
$$
c= M_\cA(\mu) = M_\cA(\lambda) -v_0 \in \cA({D+D_0}) - v_0 \enspace,  
$$
car $\cA({D+D_1}) \subset \cA({D+D_0})$ {\`a} cause de l'hypoth{\`e}se $a_0 =0$.

\smallskip

Maintenant, soit   $\gamma \in Q \cap \frac{1}{D} \, \Z^n$ et $t \in \R_+^{N+1}$ tels que $t_0+\dots+t_N = 1$, $\gamma = M_\cA(t)$ et de plus $\vartheta_{\cA, \tau}(\gamma) = \ell_\tau(t)$. On pose alors $c:= D \, \gamma  \in D\, Q  \cap \Z^n$ et, suivant l'argument pr{\'e}c{\'e}dent, on prend $\lambda \in \N^{N+1}_{D+D_0}$ tel que 
$c = M_\cA(\lambda) - v_0$ et $\lambda - D\, t \in \Lambda_\cA + w_0$. Alors
$$
D \, \vartheta_{\cA, \tau}(\gamma) \quad = \quad D\, \sum_{i=1}^N \tau_i \, t_i \quad 
=  \quad \sum_{i=1}^N \tau_i  \lambda_i  - \sum_{i=1}^N \tau_i \, (\lambda_i - D\, t_i) \quad = \quad \ell_\tau (\lambda) + O(1)\enspace, 
$$ 
d'o{\`u} $D\, \rho(D)\le s_\tau ({X_\cA};D+D_0) + O(D^n)$. En appliquant {\`a} nouveau l'asymptotique~(\ref{asymptopoids}) on obtient
$$
\frac{e_{\tau}(X_\cA)}{(n+1)!} 
\quad = \quad 
\lim_{D\to \infty} 
\frac{s_{\tau}(X_\cA;D+D_0)}{D^{n+1}} 
\quad \ge \quad   
\lim_{D\to \infty} 
\frac{\rho(D)}{D^{n}} 
\quad = \quad  
\int_Q \, \vartheta_{\cA,\tau} \ dx_1
\cdots  dx_n\enspace. 
$$
\end{demo} 

\begin{exmpl} 
Consid{\'e}rons la courbe rationnelle normale 
$C \subset \P^5$, adh{\'e}rence de l'image de
l'application $s \mapsto (1:s:s^2:s^3:s^4:s^5)$, avec le poids 
$\tau:= (-3,0,1,-1,0, -2) \in \R^6$. 
La figure suivante illustre le polytope associ{\'e} 
et sa toiture: 
\vspace{-3mm}
\begin{figure}[htbp]
$$\epsfig{file=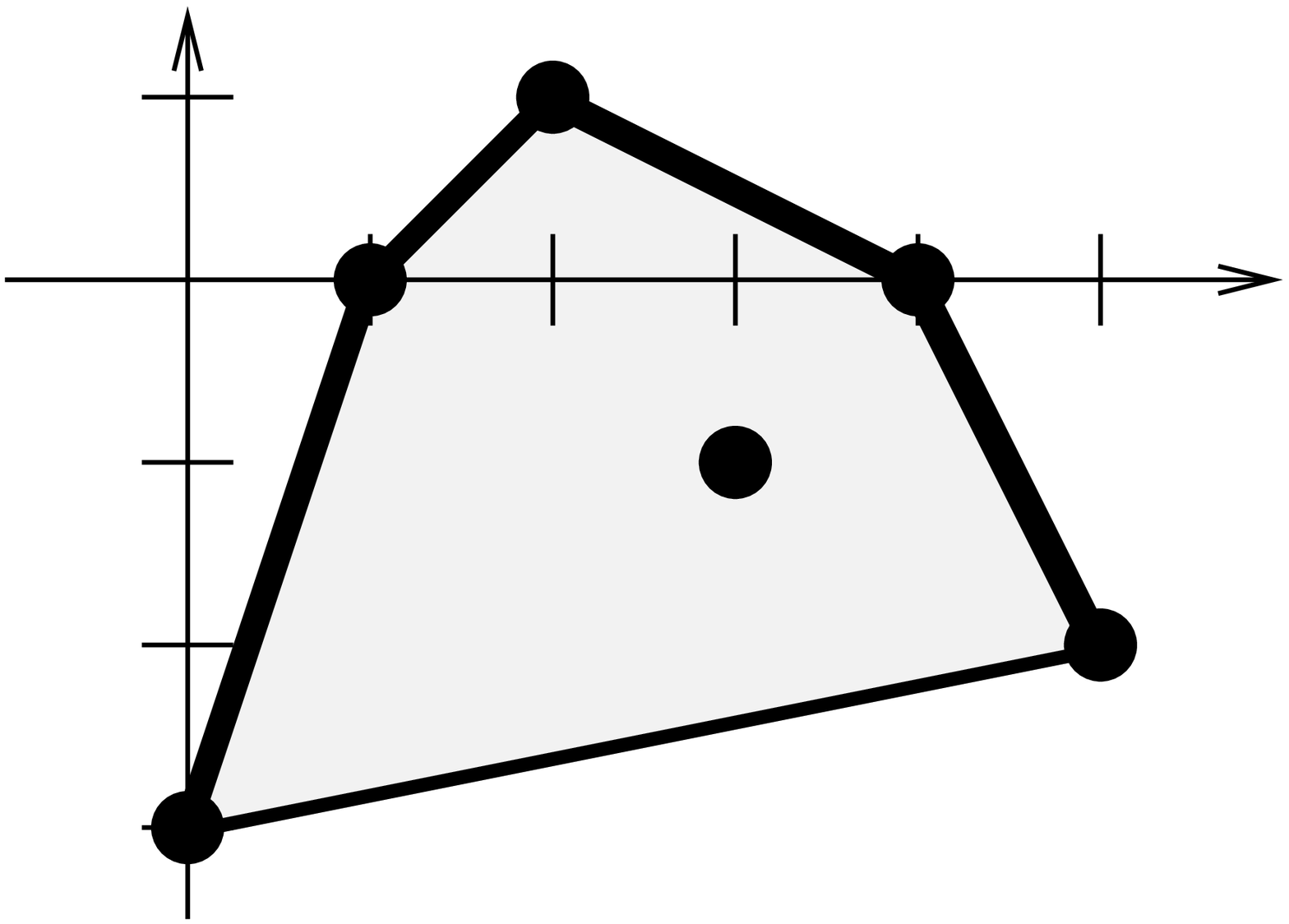, height= 32 mm}$$
\end{figure}
\vspace{-4mm} 

\noindent 
Le poids de Chow correspondant 
est 2 fois l'int{\'e}grale de la param{\'e}trisation de 
sa toiture, 
soit $\displaystyle e_{\tau}(C) = -2$. 
\end{exmpl}

Pour l'expression int{\'e}grale des poids de 
Chow de $X_\cA$, on remarque qu'alternativement
on peut l'obtenir aussi comme cons{\'e}quence 
des r{\'e}sultats de I.M.~Gelfand, 
M.M.~Kapranov et A.V.~Zelevinski sur le polytope de Newton du 
$\cA$-r{\'e}sultant~\cite[Ch.~7 et 8]{GKZ94}. 
Dans les lignes suivantes on indique bri{\`e}vement comment faire cela;  
on supposera une certaine familiarit{\'e} avec les objets et notations 
de cette r{\'e}f{\'e}rence.

Tout d'abord on remarque qu'il suffit de d{\'e}montrer l'{\'e}nonc{\'e} pour un choix g{\'e}n{\'e}rique (au sens de Zariski) du vecteur $\tau$ dans $\Q^{N+1}$, puisque les termes consid{\'e}r{\'e}s sont continus par rapport {\`a} $\tau$. Par lin{\'e}arit{\'e}, on peut supposer de plus $\tau \in \Z^{N+1}$.

Soit $T$ la sous-division poly{\'e}drale du polytope $Q_\cA$ induite 
par $\tau$ comme dans~\cite[Ch.~7,~Exemple~1.1,~p.~215]{GKZ94}; c'est une {\it triangulation} car  $\tau$ est suppos{\'e} g{\'e}n{\'e}rique. 
Le th{\'e}or{\`e}\-me~3.3 de \cite[Ch.~8,~p.~261]{GKZ94} implique $e_{\tau}(X_\cA) = \langle\tau, \varphi_{T}\rangle $, o{\`u} 
$ \varphi_{T} \in \R^{N+1}$ est la {\em fonction caract{\'e}ristique de
  $T$} 
({\it voir}~{\it loc. cit.}, \S~7.1.D, p.~220) et $ \langle \cdot, \cdot \rangle $ 
d{\'e}signe le produit scalaire standard de $\R^{N+1}$. 
La condition $\tau \in \Z^{N+1} \cap \mbox{Int}(C(T))$ est satisfaite 
par d{\'e}finition de $C(T)$, 
{\it voir}~{\it loc. cit.}, Ch.~7,~D{\'e}fn.~1.4, p.~219. 
En recollant avec la d{\'e}finition de $\varphi_{T}$ (Ch.~7,~Identit{\'e}~(1.4), 
p.~220 de {\it loc. cit.}) ceci implique
\vspace{-3mm} 
$$
e_{\tau}(X_\cA) = \sum_{i=0}^N \bigg(\tau_i \, \sum_{\sigma} n! \,
\Vol_{n}(\sigma)\bigg)
\enspace, 
$$
o{\`u} la deuxi{\`e}me somme porte sur les faces $\sigma \in T$ de 
dimension $n$ dont $a_i$ est un sommet. 
Chacun de ces $\sigma$ {\'e}tant un simplexe, on v{\'e}rifie ais{\'e}ment
$$
\frac{1}{n+1} \, \Vol_{n}(\sigma) \bigg(\sum_{i \, ; \, a_i  \in \mbox{\scriptsize S}(\sigma)}  \tau_i \bigg) 
= \int_\sigma \, \vartheta_{\cA,\tau} \ dx_1 \cdots  dx_n\enspace,  
$$
o{\`u} $\mbox{S}(\sigma)$ d{\'e}note l'ensemble des sommets du simplexe $\sigma$. Des alin{\'e}as pr{\'e}c{\'e}dents il suit enfin
\begin{eqnarray*}
\frac{e_{\tau}(X_\cA)}{(n+1)!} 
&=& \sum_{\sigma \in T} \frac{1}{n+1} \, \bigg(\sum_{i \, : \  a_i  \in \mbox{\scriptsize S}(\sigma)}  \tau_i \bigg) \, \Vol_{n}(\sigma) \\[0mm]
&=& \sum_{\sigma \in T} \int_\sigma\,\vartheta_{\cA,\tau} \ dx_1\cdots  dx_n 
=\int_{Q_\cA} \, \vartheta_{\cA,\tau} \ dx_1 \cdots  dx_n
\enspace.
\end{eqnarray*} 

\medskip

\begin{cor} \label{corpoid-torique}
Avec les notations de la proposition~\ref{poid-torique} on a
$$e_\tau(X_{\cA}) + e_{-\tau}(X_{\cA}) = (n+1)!\, 
{\rm Vol}_{{n+1}}(Q_{{\cA},\tau})\enspace.$$
\end{cor}
\begin{demo} Il suffit de remarquer que $\vartheta_{{\cA},\tau}(x) + \vartheta_{{\cA},-\tau}(x)$ est {\'e}gal {\`a} l'{\'e}paisseur de $Q_{{\cA},\tau}$ au-dessus du point $x\in Q_{\cA}$.
\end{demo}

Ce r{\'e}sultat nous permet de calculer la hauteur d'une vari{\'e}t{\'e} 
torique d{\'e}finie sur un corps de base $F:= \K(t)$ 
(corps des fractions rationnelles en une variable $t$ sur $\K$). 
Pour $\alpha \in F^\times$ et $v \in \K$ 
on pose $\Ord_v (\alpha)$ l'ordre d'annulation de la fonction $\alpha$ 
au point $v$. 
On pose {\'e}galement $\Ord_\infty(\alpha):= \deg(\alpha)$
l'ordre d'annulation de $\alpha$ {\`a} l'infini. 
Avec ces conventions on a  
\begin{equation} \label{fproduit} 
\sum_{v \in \{ \infty\} \cup \K} \Ord_v (\alpha) = 0
\enspace. 
\end{equation} 
On pose maintenant, pour $\gamma= (\gamma_0: \cdots : \gamma_M) \in
\P^M(F)$,
$$
h(\gamma):= \sum_{v \in \{ \infty\} \cup \K} 
\max\{ 
\Ord_v (\gamma_0), \dots, \Ord_v (\gamma_M) \} \ \in \N 
\enspace,
$$
qui est bien d{\'e}finie gr{\^a}ce {\`a} la formule~(\ref{fproduit}). 
Finalement, la hauteur d'une $F$-sous-vari{\'e}t{\'e} $X \subset \P^N(\ov{F})$
est par d{\'e}finition, celle de sa forme de Chow:
$h(X):= h(\Ch_X)$, par analogie avec le cas des corps de nombres.

\begin{prop} 
Soit $\cA \in (\Z^n)^{N+1}$ tel que $L_\cA = \Z^n$, 
$\tau\in\Z^{N+1}$. 
Soit  
$X_{\cA,\tau} \subset \P^N(\ov{F})$ l'adh{\'e}rence de Zariski 
de l'application monomiale
$$
\T^n(\ov{F}) \to \P^N(\ov{F}) \enspace, 
\quad \quad 
s \mapsto (t^{\tau_0}\, s^{a_0}:\cdots:t^{\tau_N}\,s^{a_N})
\enspace,
$$
alors $h(X_{\cA,\tau}) 
= (n+1)!\,{\rm Vol}_{{n+1}}(Q_{\cA,\tau})$. 
\end{prop} 

\begin{demo} 
Soit $\Ch_{X_\cA} \in \K[U_0, \dots, U_n] $ la forme de Chow de
$X_\cA$, 
la vari{\'e}t{\'e} torique correspondant {\`a}  $\tau = {\mathbf 0}_{N+1}$. 
Alors
$$
\Ch_{X_{\cA, \tau}} = 
\Ch_{X_\cA}\Big( t^{-\tau_j} \, U_{i\,j} \, : \ 0\le i \le n, 
\ 0 \le j \le N \Big) \in F[U_0, \dots, U_n]
$$
est un polyn{\^o}me de Laurent en $t$,
dont sa hauteur est {\'e}gale {\`a} la somme de ses 
valuations $t$-adique et $\displaystyle\frac{1}{t}$-adique, 
soit:
$$
h(X_{\cA, \tau}) = 
\Ord_0 (\Ch_{X_{\cA, \tau}})
+ 
\Ord_\infty (\Ch_{X_{\cA, \tau}}) 
= 
e_{\tau}(X_\cA) + e_{-\tau}(X_\cA) = (n+1)! \,{\rm Vol}_{{n+1}}(Q_{\cA,\tau}) 
$$
d'apr{\`e}s le corollaire~\ref{corpoid-torique}. 
\end{demo}


\typeout{Demonstration du theoreme 1} 

\subsection{D{\'e}monstration du th{\'e}or{\`e}me~\ref{thm1}}

\label{Demonstration du theoreme 1}

Dans ce paragraphe on explicite 
la fonction de Hilbert arithm{\'e}tique 
et la hauteur normalis{\'e}e d'une vari{\'e}t{\'e} torique. 
On montre que ces invariants num{\'e}riques 
se d{\'e}composent en somme finie de contributions locales, 
chaque terme {\'e}tant un poids de Hilbert
et un poids de Chow, respectivement. 

\medskip 

On reprend $\Qbar$ comme corps de base. 
Soit $\cA \in(\Z^n)^{N+1}$
tel que $L_\cA= \Z^n$ et $\alpha \in(K^\times)^{N+1}$. 
Soit 
$X_{\cA,\alpha} \subset \P^N $ la sous-vari{\'e}t{\'e} torique 
associ{\'e}e et notons 
$\tau_{\alpha v}:=(\log|\alpha_0|_v, \dots, \log|\alpha_N|_v) \in
\R^{N+1}$
le poids de $\alpha$ relatif {\`a} une place $v \in M_K$.

\begin{prop} \label{decom Hnorm} 
$\displaystyle 
\Hnorm(X_{\cA,\alpha}; D)
=\sum_{v\in M_K}\frac{[K_v:\Q_v]}{[K:\Q]} \, s_{\tau_{\alpha
    v}}(X_\cA; D) $ 
pour tout $D\in\N$. 
\end{prop}

\begin{demo} 
Soit $\Big\{ (x^\lambda)^\vee \, : \ \lambda \in \N^{N+1}_D \Big\} $ la base 
de $K[x_0, \dots, x_N]_D^\vee$
duale de celle
des mon{\^o}mes, et pour $c\in \cA(D)$ on pose 
$$
q_c := \sum_{\lambda\in  \N^{N+1}_D\,; \ M_\cA(\lambda) =c } 
\alpha^\lambda \, (x^\lambda)^\vee 
\enspace.
$$
Ainsi, pour une forme  $f = \sum_\lambda f_\lambda\, x^\lambda \in K[x_0, \dots, x_N]_D$ 
 et un point  $s  \in \T^n$ on a 
$$
f\circ\varphi_{{\cA,\alpha}}(s) = \sum_{\lambda \in \N^{N+1}_D} f_\lambda \, \alpha^\lambda \,
s^{M_\cA(\lambda)} = 
\sum_{c \in \cA(D) } 
q_c(f) \, s^c \enspace. 
$$
Soit $I:= I(X_\arith) \subset K[x_0, \dots, x_N]$ l'id{\'e}al de 
d{\'e}finition de $X_{\cA,\alpha}$;  
on voit ainsi que 
$f \in I_D$ si et seulement si $q_c(f) =0$ pour tout $c$. 
Donc l'ensemble
$
\Big\{ q_c \, : \  c \in \cA(D) \Big\} $
forme une base de l'annulateur $\Ann(I_D)$ et 
alors la  proposition~\ref{Hnorm dual}
implique
$$
\cH_{\rm norm}(X_{\cA,\alpha};D)= 
\sum_{v \in M_K} \frac{[K_v:\Q_v]}{[K:\Q]} \, 
\log \bigg| \bigwedge_{c \in \cA(D)} q_c \bigg|_v \enspace.
$$
Pour $v\in M_K$ on a 
$$
\log |q_c|_v = 
\log \max \Big\{ |\alpha^\lambda|_v \, : \,  \lambda \in \N^{N+1}_D, \, 
M_\cA(\lambda)=c
\Big\} = \max \Big\{ \ell_{\tau_{\alpha v}}(\lambda) \, : \, \lambda \in
\N^{N+1}_D, 
\, M_\cA(\lambda) =c\Big\} 
$$
o{\`u} $\ell_{\tau_{\alpha v}}$ d{\'e}signe la forme lin{\'e}aire
associ{\'e}e au poids ${\tau_{\alpha  v}}$.
Le fait que les supports des composantes non nulles des vecteurs $q_c$ 
soient disjoints 
et la proposition~\ref{poid-torique} 
entra{\^\i}nent 
$$
\log \bigg| \bigwedge_{c} q_c \bigg|_v 
= \sum_c \log |q_c|_v 
= s_{\tau_{\alpha v}} (X_\cA)
\enspace, 
$$
ce qui ach{\`e}ve la
d{\'e}monstration. 
\end{demo} 

\medskip

Maintenant on {\'e}nonce et d{\'e}montre le th{\'e}or{\`e}me~\ref{thm1} 
dans sa version g{\'e}n{\'e}rale, \cad pour un vecteur  $\cA$ quelconque. 
On rappelle que $\mu_\cA$ d{\'e}signe la forme volume 
sur l'espace lin{\'e}aire $L_\cA \otimes_\Z \R$ associ{\'e}e 
{\`a} $\cA$, {\it voir}~\S~\ref{monomiales}.

\begin{thm} \label{thm1general}
Soit $\geom \in (\Z^n)^{N+1}$ et 
$\alpha \in (K^\times)^{N+1}$. 
Posons
$r:= \dim(Q_\cA)$, alors
$$
\wh{h}(X_{\cA,\alpha}) = (r+1)! \, \sum_{v \in M_K} 
\frac{[K_v:\Q_v]}{[K:\Q]} \, 
\int_{Q_\cA} \,
\vartheta_{\cA,\tau_{\alpha v}} \ \mu_\cA 
\enspace .
$$

\end{thm} 

\begin{demo}[D{\'e}monstration du 
th{\'e}or{\`e}me~\ref{thm1general} et de la proposition~\ref{asymptotique}]

Suppossons pour le moment \newline $L_\cA = \Z^n$.
Notons que $\tau_{\alpha  v} = 0$ pour presque tout $v$,
donc la d{\'e}composition de $\Hnorm$ ne contient
qu'un nombre fini de poids de Hilbert non nuls. 
On d{\'e}duit de la proposition pr{\'e}c{\'e}dente et de l'asymptotique~\cite[Prop.~2.11]{Mum77} 
$$
\Hnorm(X_{\cA,\alpha}; D) = \frac{c(X_\arith)}{(n+1)!} \, D^{n+1} 
+ O(D^n)
$$ 
avec $ \displaystyle c(X_{\cA,\alpha}):= \sum_{v} \frac{[K_v:\Q_v]}{[K:\Q]} \,
e_{\tau_{\alpha v}}(X_\cA)$.

Soit  $k \in\N^*$, 
alors
$ [k] \, X_{\cA,\alpha} = X_{\cA, [k]\alpha} $
puisque $[k] \, X_\cA^\circ = X_\cA^\circ$; en particulier
$\deg([k]\, X_{\cA,\alpha})= \deg(X_\arith)$.
On a aussi  
$\tau_{[k]\alpha \, v} = 
k \, \tau_{\alpha  v} $, ce qui implique  
$c([k] \, X_{\cA,\alpha} ) = k \cdot c(X_{\cA,\alpha} )$. 
De plus
$$ 
0 \le h([k] \, X_{\cA,\alpha}) -  c([k] \, X_{\cA,\alpha}) 
\le \frac{n+1}{2} \, \log(N+1)\, \deg(X_\arith)
\enspace, 
$$
gr{\^a}ce {\`a} l'encadrement~(\ref{encadrementhilbert}), 
donc 
$ \displaystyle
c(X_{\cA,\alpha}) =
\deg(X_\arith) \, \lim_{k\rightarrow\infty} 
\frac{h([k]\, X_{\cA,\alpha})}{k \, \deg([k] \, X_\arith)} 
= \wh{h}(X_{\cA,\alpha})$. 
On conclut par application de la proposition~\ref{poid-torique}.

\smallskip 

Consid{\'e}rons maintenant le cas g{\'e}n{\'e}ral. 
Soit  $\eta : \R^r \hookrightarrow \R^n$ une application lin{\'e}aire
injective 
d{\'e}finie sur $\Z$ telle que 
$\eta(\Z^r) = L_\cA$. 
On pose 
$b_i:= \eta^{-1}(a_i) \in \Z^r$
puis 
$\cB:= (b_0, \dots, b_N) \in (\Z^r)^{N+1} $,  
alors $X_{\cB, \alpha} = X_\arith$~\cite[Ch.5,~Prop.~1.2]{GKZ94}. 
On v{\'e}rifie $\eta(Q_\cB) = Q_\cA$, 
$\vartheta_{\cB, \tau_{\alpha v} } = \vartheta_{\cA, \tau_{\alpha v}} \circ \eta$
 et 
$\eta^*(\mu_\cA) = dx_1 \cdots dx_r$, 
d'o{\`u}
$$
\int_{Q_\cB} \,
\vartheta_{\cB,\tau_{\alpha v}} \ dx_1 \cdots dx_r
=
 \int_{Q_\cA} \,
\vartheta_{\cA,\tau_{\alpha v}} \ \mu_\cA
\enspace.  
$$
On a $L_\cB= \eta^{-1}(L_\cA) = \Z^r$,
ce qui nous ram{\`e}ne au cas pr{\'e}c{\'e}demment consid{\'e}r{\'e}. 
\end{demo} 

Alternativement, 
on peut reformuler le th{\'e}or{\`e}me~\ref{thm1} en 
$$
\wh{h}(X_{\cA,\alpha}) = (n+1)! \, \int_{Q_\cA} \theta_{\cA,\alpha} \ 
dx_1  \cdots dx_n
$$
avec $\displaystyle \theta_{\cA,\alpha}:= \sum_{v} 
\frac{[K_v:\Q_v]}{[K:\Q]} \, \vartheta_{\cA,\tau_{\alpha v}}$.
Cette fonction est concave et affine par morceaux, 
car elle est somme
d'un nombre fini de fonction de ce type; de plus  
$\theta_\arith $ est {\`a} valeurs {\it positives}: 
pour $0 \le i \le N$ on a 
$
\vartheta_{\cA,\tau_{\alpha v}}(a_i) \ge \log|\alpha_i|_v 
$ 
et donc 
$$
\theta_{{\cA,\alpha}}(a_i) \ge \sum_{v \in M_K} 
\frac{[K_v:\Q_v]}{[K:\Q]} \, \log|\alpha_i|_v = 0
$$
par la formule du produit. 
De l{\`a} on d{\'e}duit 
$\theta_\arith(x) \ge 0 $ pour tout $x \in Q_\cA$, 
{\`a} cause de la concavit{\'e}.

\begin{exmpl} 

Soit $\cA:=(a_0,\dots,a_{N-1},a_N)\in (\Z^n)^{N+1}$ 
tel que $L_\cA= \Z^n$ et  
$Q_\cA= \Conv ( a_0, \dots, a_{N-1})$, 
et soit aussi $\alpha=(1,\dots,1, \alpha_N)\in (K^\times)^{N+1}$. 
Pour chaque $v \in M_K$ on a 
$$
\int_{Q_\cA} \vartheta_{\cA,\tau_{\alpha v}} \, dx_1 \cdots dx_n 
=  \frac{1}{n+1} \, \Vol_n(Q_\cA) \, \log \max \{ 1, |\alpha_N|_v\}
=  \frac{1}{n+1} \, \Vol_n(Q_\cA) \, \log |\alpha|_v
$$ 
d'o{\`u}   
$$ 
\wh{h}(X_{\cA,\alpha}) = n! \, \Vol_n( Q_\cA) \,  \hnorm(\alpha) 
= \deg(X_\arith) \,  \hnorm(\alpha)
\enspace.
$$ 
Par exemple,  
consid{\'e}rons la courbe plane
$C \subset \P^2$, adh{\'e}rence de l'image de  $ \T^1 \to \P^2$, 
$\displaystyle s \mapsto (1 : 3 \, s^2: s^3)$. 
Les figures suivantes montrent les polytopes associ{\'e}s et leurs
toitures: 


\vspace{25mm} 

\begin{figure}[htbp]

\begin{picture}(0,36) 

\put(30,50){\epsfig{file=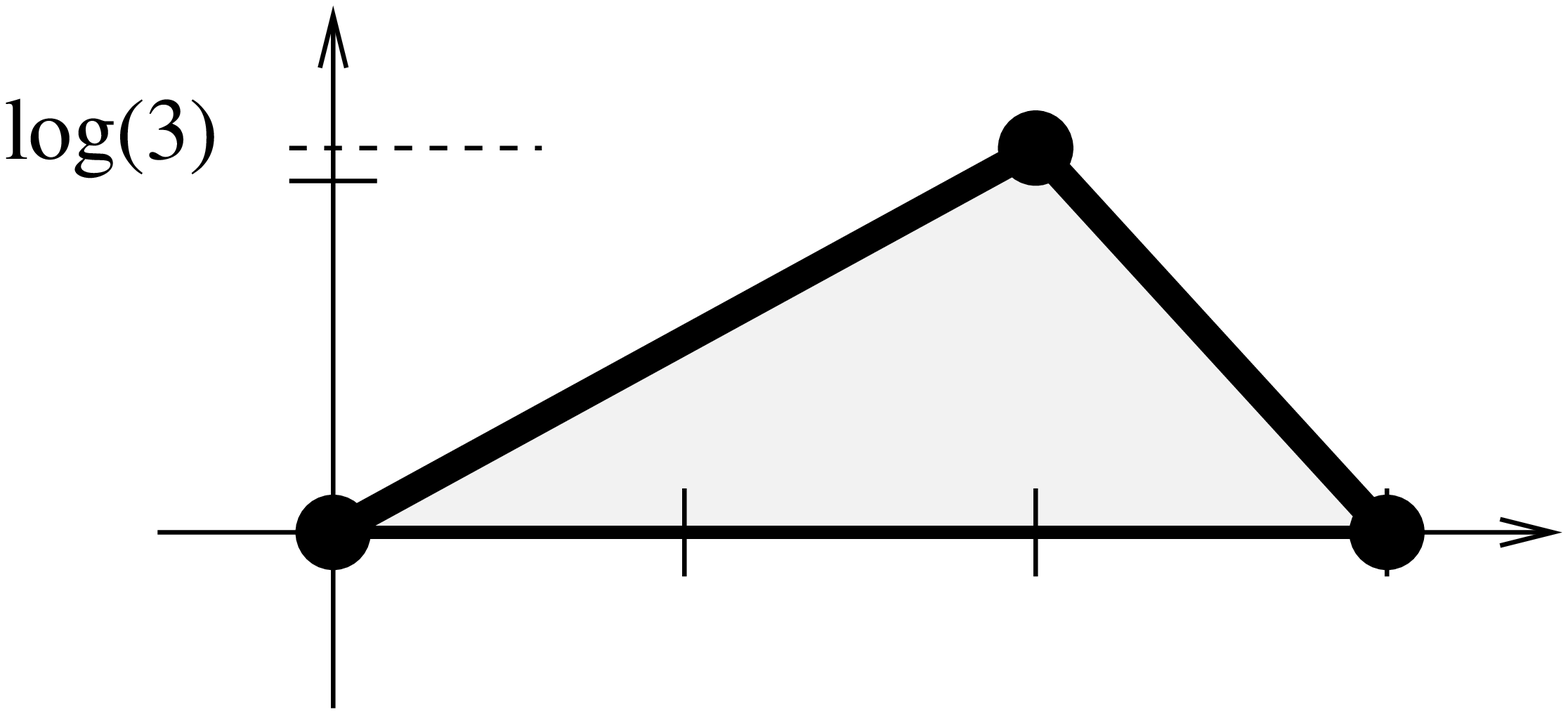, height= 21 mm}}
\put(160,78){$v=\infty$} 


\put(265,29){\epsfig{file=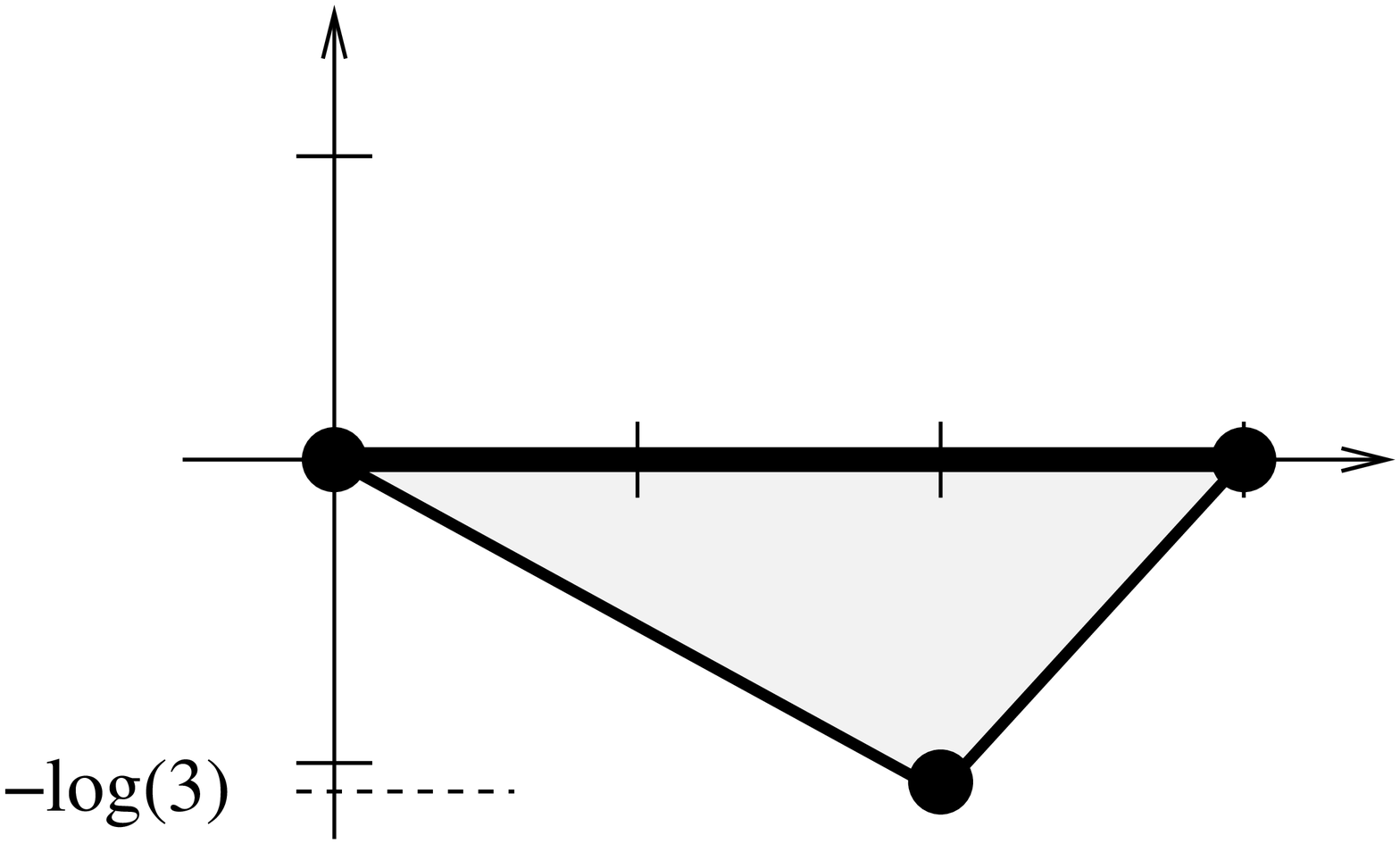, height= 28.5 mm}} 
\put(395,78){$v=3$}

\end{picture} 

\end{figure}

\vspace{-18mm} 


\noindent 
et 
$\vartheta_v \equiv 0$
pour $v \ne \infty, 3$; 
donc $\hnorm(C) = 3 \, \log(3)$.

\end{exmpl} 
  
Cet exemple contient le cas des {\it hypersurfaces} binomiales, 
qui  est toutefois plus simple {\`a} traiter directement. 
Soit $f= x^a - \lambda \, x^b \in K[x_0, \dots, x_N] $
une forme irr{\'e}ductible, 
alors pour chaque place archim{\'e}dienne 
$v \in M_K^\infty$ on peut calculer la mesure de
Mahler correspondante {\`a} l'aide de la formule de Jensen, on a 
$
\mu_v(f):= \mu(x^a - |\lambda|_v \, x^b )= \log \max \{ 1,
|\lambda|_v\}$. 
Ainsi 
$$
\hnorm(Z(f)) = 
\sum_{v\in M_K^\infty}\frac{[K_v:\Q_v]}{[K:\Q]} \, \mu_v(f) \ + 
\sum_{v\in M_K\setminus M_K^\infty}\frac{[K_v:\Q_v]}{[K:\Q]} 
\, \log|f|_v 
=h((1:\lambda)) \enspace.
$$

\medskip 

Maintenant soit $[-1]: (\P^N)^\circ \to (\P^N)^\circ$, 
$(x_0:\cdots:x_N) \mapsto  (x_0^{-1}:\cdots:x_N^{-1})$
l'application d'inversion. 
On obtient une expression 
pour la somme de la  hauteur de
$X_{\cA,\alpha}$ et de sa vari{\'e}t{\'e} sym{\'e}trique, 
dont les
termes locaux sont des volumes de polytopes:

\begin{cor} \label{hnormvol}
Soit $\cA \in (\Z^n)^{N+1}$ tel que $L_\cA = \Z^n$ et $\alpha \in
(K^\times)^{N+1}$, alors 
$$
\wh{h}(X_{\cA,\alpha})+ \wh{h}([-1] \, X_{\cA, \alpha}) = 
(n+1)! \, \sum_{v
  \in M_K} \frac{[K_v:\Q_v]}{[K:\Q]} \, {\rm
  Vol}_{{n+1}}(Q_{\cA,\tau_{\alpha v}})
\enspace.
$$
\end{cor}

\begin{demo} 
On a 
$
[-1] \, X_{\cA,\alpha}= X_{\cA,[-1] \, \alpha}$ et 
$\tau_{[-1] \alpha \  v} = - \tau_{\alpha v}$; le r{\'e}sultat est donc
une une cons{\'e}quence imm{\'e}diate du th{\'e}or{\`e}me~\ref{thm1} et du 
corollaire~\ref{corpoid-torique}. 
\end{demo}

Consid{\'e}rons un espace ambiant   $\T^N$ muni de la hauteur 
normalis{\'e}e $\hnorm_\iota$ 
relative aux plongement $\iota: \T^N \hookrightarrow \P^{2^N-1}$, 
composition de l'inclusion $\T^N \hookrightarrow (\P^1)^N$ 
avec le plongement de Segre $(\P^1)^N \hookrightarrow \P^{2^N-1}$, 
{\it voir}~\S~\ref{hauteur normalisee}. 
Le plongement $\iota $ est {\'e}quivariant, donc 
pour une sous-vari{\'e}t{\'e} $Y \subset \T^N$
on a
$\hnorm_\iota(Y) = \hnorm_{\P^N}(\ov{\iota(Y)})$
d'apr{\`e}s la proposition~\ref{h_varphi}. 
De plus, notons que cette hauteur est invariante par rapport {\`a} l'inversion:
on a $[-1] \, Y = \sigma(Y)$ o{\`u}  $\sigma: ( (x_1:y_1), \dots, (x_N:y_N)) \mapsto 
( (y_1:x_1), \dots, (y_N:x_N) )$, donc 
$\hnorm_\iota([-1] \, Y) = \hnorm_\iota(Y)$.

Soit maintenant $\cB = (b_1, \dots, b_N) \in (\Z^n)^{N}$ tel que 
$L_\cB=\Z^n$ et
$\beta= (\beta_1, \dots, \beta_N) \in (K^\times)^{N}$. 
Consid{\'e}rons 
l'application monomiale
$$
\T^n \to \T^N \enspace, 
\quad \quad 
s \to (\beta_1 \, s^{b_1} , \dots, \beta_N \, s^{b_N}) 
$$
et posons $Y_{\cB, \beta} \subset  \T^N$ son image. 
On v{\'e}rifie alors que $\ov{s(Y_{\cB, \beta})} \subset \P^{2^N -1}$ est 
la vari{\'e}t{\'e} torique
correspondant aux vecteurs
$$
\bigg( \sum_{i \in I} b_i \, : \ I \subset \{ 1, \dots, N\} \bigg) \in
(\Z^n)^{2^N} 
\enspace, \quad \quad  
\bigg( \prod_{i \in I} \beta_i \, : \ I \subset \{ 1, \dots, N\} \bigg) \in
(K^\times)^{2^N} 
\enspace.
$$ 
Pour $v \in M_K$ et $i=1, \dots, N$ on consid{\`e}re le segment 
$P_{i\, v} := \Conv \Big( \mathbf{0}, (b_i , \log|\beta_i|_v) \Big)
\subset \R^{n+1}$; 
ainsi le polytope associ{\'e} {\`a} la vari{\'e}t{\'e} torique 
$s(Y_{\cB, \beta})$ pour la place $v$
est donn{\'e} par la somme de Minkowski
$P_{0, v} + \cdots + P_{N, v} $.  
Le corollaire~\ref{hnormvol} entra{\^\i}ne
\begin{equation} 
\wh{h}_\iota (Y_{\cB, \beta})
=
\frac{(n+1)!}{2}  \, \sum_{v
  \in M_K} \frac{[K_v:\Q_v]}{[K:\Q]} \, {\rm
  Vol}_{{n+1}}
(P_{0, v} + \cdots + P_{N , v} )
\enspace.
\end{equation}
Comme illustration de cette formule, consid{\'e}rons la courbe 
$C \subset \T^3 \hookrightarrow (\P^1)^3 
{\hookrightarrow \atop \mbox{\scriptsize Segre} }\P^7$ image de l'application
$\displaystyle s \mapsto  \Big({2} \, s, \frac{1}{3}\, s, \frac{3}{2} \, s\Big) $. 
Les figures suivantes montrent les segments $P_{i,v}$
(en trait discontinu)  
et les polytopes associ{\'e}s, pour chaque
place $v \in M_\Q$: 


\vspace{45mm} 

\begin{figure}[htbp]

\begin{picture}(0,0) 

\put(18,55){\epsfig{file=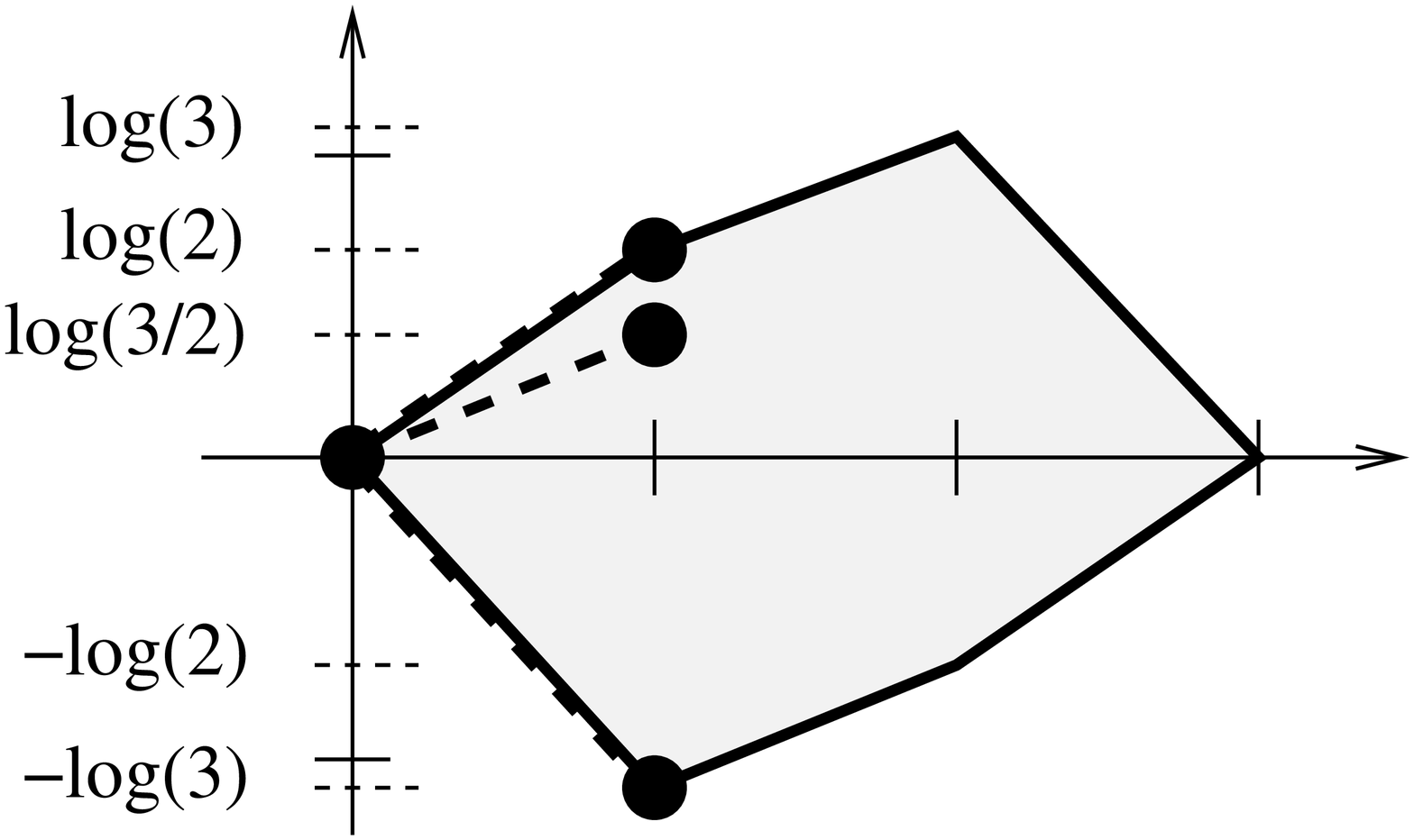, height= 29 mm}} 
\put(155,75){$v=\infty$} 

\put(261,55){\epsfig{file=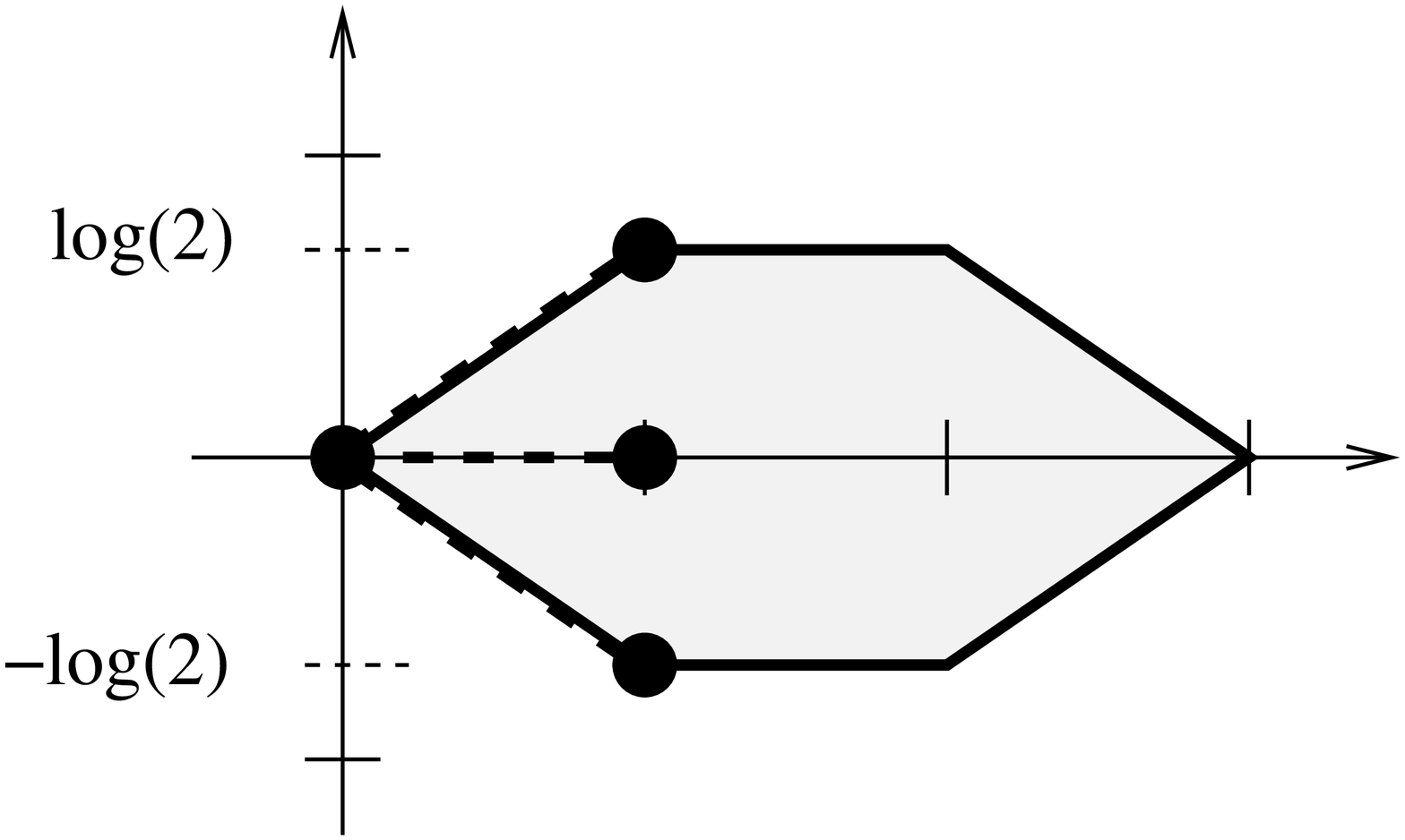, height= 29 mm}} 
\put(390,75){$v=2$} 

\put(20,-33){\epsfig{file=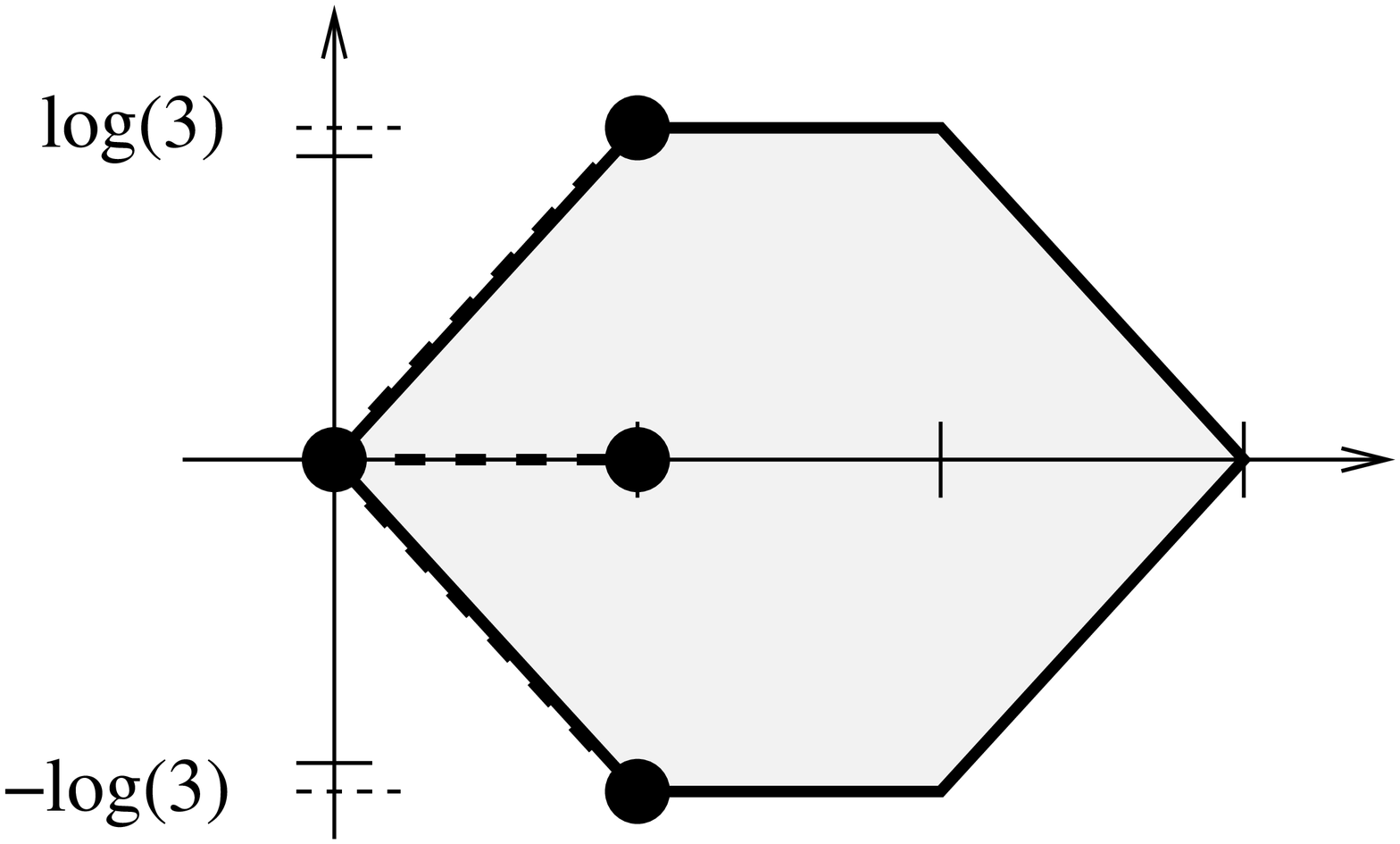, height= 29 mm}} 
\put(155,-12){$v= 3$} 

\put(279,-33){\epsfig{file=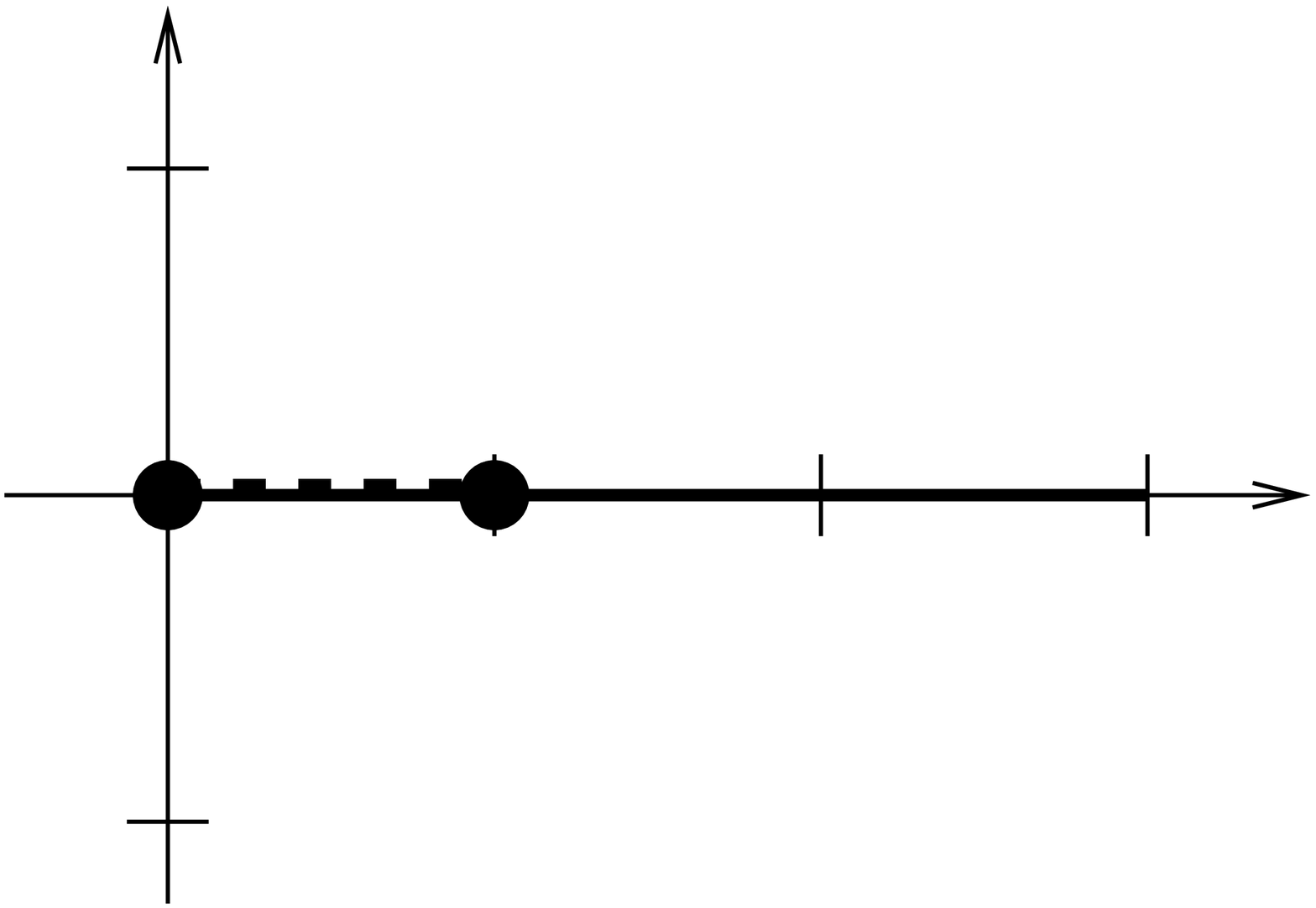, height= 29 mm}} 
\put(390,-12){$v\ne \infty, 2,3$}

\end{picture} 

\end{figure}

\vspace{9mm} 


\noindent donc $\hnorm_\iota (C) = 6 \, \log(2) + 6 \, \log(3)$.



%
%


\typeout{Sur quelques constructions standard} 

\section{Sur quelques constructions standard} 

\label{multi-hauteurs} 

\setcounter{equation}{0}
\renewcommand{\theequation}{\thesection.\arabic{equation}}

Dans la suite nous {\'e}tendons les r{\'e}sultats de la section pr{\'e}c{\'e}dente au cas
d'un tore  muni de 
plusieurs plongements monomiaux. 
D'abord, on introduit les multipoids de Chow et on {\'e}tablit ses
propri{\'e}t{\'e}s essentielles; 
notre approche s'appuie sur les 
formes
r{\'e}sultantes d'id{\'e}aux multihomog{\`e}nes,
d{\'e}j{\`a} consid{\'e}r{\'e}es
au paragraphe~\ref{hauteur normalisee}.

Puis, on introduit la notion d'int{\'e}grale mixte d'une famille
de fonctions concaves. 
On montre ensuite que le multipoids du tore $\T^n$ relatif {\`a} plusieurs
plongements
monomiaux s'{\'e}crit comme une int{\'e}grale mixte. 
Finalement, on explicite la multihauteur normalis{\'e}e de $\T^n$ relative {\`a}
des plongements monomiaux. 
Cette multihauteur se d{\'e}compose comme somme de contributions locales, 
chaque terme {\'e}tant l'int{\'e}grale mixte de
l'ensemble de fonctions concaves et affines
par
morceaux associ{\'e}es {\`a} ces plongements monomiaux, pour la place $v \in
M_K$ correspondante.

On {\'e}tudie {\'e}galement le comportement de la famille de fonctions
$\Theta_\arith$ par rapport {\`a} d'autres constructions: 
d{\'e}composition de $X_\arith$  en orbites sous l'action $*_\cA$, 
formation de joints, produits de Segre et plongements de Veronese.


\typeout{Decomposition en orbites} 

\subsection{D{\'e}composition en orbites} 

\label{orbites}

Soit $\cA \in (\Z^n)^{N+1}$ et $\alpha \in (K^\times)^{N+1}$, 
les orbites de l'action $*_\cA$ sur $X_\arith$ sont en correspondance
avec  
l'ensemble $F(Q_\cA)$  des faces du polytope $Q_\cA$. 
Pour chaque face $ P$ on consid{\`e}re un point 
$\alpha_P := (\alpha_{P, \, 0} : \cdots : \alpha_{P, \, N} ) \in \P^N$ 
d{\'e}fini par $  \alpha_{P, \, j}  :=\alpha_j$ 
si \,$ a_j \in P $  et  $  \alpha_{P, \, j}
:= 0$ sinon; 
la bijection est  donn{\'e}e par~\cite[Ch.~5, Prop.~1.9]{GKZ94},  
\cite[\S~3.1]{Ful93} 
$$
P \mapsto X_{\arith, P}^\circ := \T^n *_\cA \alpha_P \ \subset \P^N 
\enspace. 
$$
On a la d{\'e}composition 
$ \displaystyle 
X_\arith  = \bigsqcup_{P \in F(Q_\cA) } X_{\arith, P}^\circ 
$, posons  $N(P):= {\rm Card}\{ i \, : \ a_i \in P\} -1$ et
$$
\cA(P) = (a_i \, : \ a_i \in P) \in(\Z^n)^{N(P)+1}
\enspace,\quad \quad 
\alpha(P) := ( \alpha_i \, : \  a_i \in P) \in (K^\times)^{N(P)+1} \enspace.
$$
On v{\'e}rifie que $X_{\arith, P}^\circ \subset \P^N$ 
est l'orbite principale d'une vari{\'e}t{\'e} torique 
contenue
dans un sous-espace standard $E \cong \P^{N(P) }$. 
Restreinte {\`a} ce sous-espace, elle
s'identifie {\`a} la sous-vari{\'e}t{\'e} 
$ X_{\cA(P),\alpha(P)}^\circ \subset \P^{N(P) }$, 
de dimension {\'e}gale {\`a} la dimension (r{\'e}elle) de la face  $P$.

\smallskip 

Finalement, chaque orbite est en correspondance avec une certaine face $P$
du polytope $Q_\cA$. 
Similairement la famille de 
fonctions 
$\Theta_{\cA(P), {\alpha(P)}} $
associ{\'e}e {\`a} l'orbite, est la {\it restriction} {\`a} la face $P$ 
de la famille $\Theta_\arith$ 
associ{\'e}e {\`a} l'action $*_\cA$ et au point $\alpha$:

\begin{prop} \label{theta-orbite}
Pour chaque $v \in M_K$ on a 
\enspace
$ \displaystyle 
\vartheta_{\cA(P), \tau_{\alpha(P) \, v}} 
= \vartheta_{\cA, \tau_{\alpha \, v}}  \Big|_P$ \enspace.
\end{prop}


\typeout{Multipoids de Chow}

\subsection{Multipoids de Chow} 

\label{multipoids}

On se place {\`a} nouveau sur un corps de base $\K$
alg{\'e}\-bri\-que\-ment clos.
Soient $X$ une vari{\'e}t{\'e} (non plong{\'e}e)  
de dimension $n$ et  
$
\varphi_i : X \dashrightarrow \P^{N_i} $ des
applications rationnelles, 
{\it plongements} sur un ouvert dense 
de $X$ ($i=0, \dots, n$). 
On pose
$\Phi : X \to \prod_{i=0}^n\P^{N_i}$,
$x \mapsto (\varphi_0(x), \dots, \varphi_n(x))$, 
et on consid{\`e}re la forme r{\'e}sultante
$$
\Res_{\varphi_0, \dots, \varphi_n} := 
\res_{e_0, \dots, e_n}\Big( I(\Phi(X)) \Big) \in 
\K[U_0, \dots, U_n]
$$ 
de l'id{\'e}al premier multihomog{\`e}ne 
$I(\Phi(X)) \subset \K[x_0, \dots, x_n]$
d'indice 
$e_0, \dots, e_n$, 
dans la notation du paragraphe~\ref{hauteur normalisee}.

\begin{defn} \label{defn multipoids} 
Soient  $\tau_0 \in \R^{N_0+1}, \dots, 
\tau_n \in \R^{N_n+1} $ des poids et  $t$ une variable 
additionnelle.  
Le {\em multipoids de Chow de 
 $X$ relatif {\`a} $(\varphi_0, \tau_0) , \dots, (\varphi_n, \tau_n)$} 
est 
$$
e_{\tau_0 , \dots,  \tau_n}(\varphi_0, \dots, \varphi_n;  X):= 
\deg_t \Big( \Res_{\varphi_0, \dots, \varphi_n} 
(t^{\tau_{i\, j}} \, U_{i\, j}\, 
: \ 0 \le i \le n, \ 0 \le j \le N_i) \Big) 
\ \in \R
\enspace.   
$$
\end{defn}

Remarquons que pour $\tau_0 := (1, \dots, 1) $ 
et $\tau_i := 0 $ pour $ 1 \le i \le n$, 
$$
e_{\tau_0 , \dots,  \tau_n}(\varphi_0, \dots, \varphi_n;  X):= 
\deg_{U_0} ( \Res_{\varphi_0, \dots, \varphi_n} )
$$
ne d{\'e}pend pas du plongement $\varphi_0$~\cite[Prop.~3.4]{Rem01a};
c'est le {\it multidegr{\'e}} de $X$ relatif {\`a} $\varphi_1, 
\dots, \varphi_n$. 

\smallskip 

Soient $\varphi: X \to \P^N$, $\psi: X \to \P^P$ des applications 
projectives, on pose
$$
\varphi \oplus \psi: X \to \P^{(N+1) \, (P+1) -1 }
$$
la composition de $X \to \P^N \times \P^P, 
x \mapsto (\varphi(x), \psi(x))$ avec le plongement
de Segre. 
On obtient ainsi une structure de semi-groupe commutatif sur
l'ensemble des
applications projectives de $X$. 
Pour $\tau \in \R^{N+1}$, $\omega \in \R^{P+1}$ on pose 
$\tau \oplus \omega := \Big( \tau_i + \omega_j \, : \ 
0\le i \le N, \ 0 \le j \le P \Big)
\in \R^{(N+1) \, (P+1) } $, ainsi
$(\varphi, \tau) \oplus (\psi, \omega) := ( \varphi \oplus \psi, 
\tau 
\oplus \omega )$.

\begin{lem} \label{proprie-multipoids} 
\

\smallskip 

{\noindent (a)} \ 
L'application $\displaystyle ( (\varphi_0, \tau_0), \dots,  (\varphi_n,
\tau_n)) 
\mapsto e_{\tau_0 , \dots, \tau_n}(\varphi_0, \dots, \varphi_n;  X) $ 
est sym{\'e}trique et lin{\'e}aire en chaque variable
par rapport {\`a} $\oplus$ 
\enspace; 

\smallskip 

{\noindent (b)} \  
$ \displaystyle e_{\tau , \dots, \tau}(\varphi, \dots, \varphi;   X) 
= e_{\tau}(\ov{\varphi(X)}) $
\enspace;  

\noindent (c) \
$
\displaystyle  e_{\tau_0 , \dots, \tau_n}(\varphi_0, \dots, \varphi_n;
 X) = 
\frac{1}{(n+1)!} \, \sum_{j=0}^n (-1)^{n+1 - j} 
\hspace{-3mm} \sum_{0 \le i_0 < \cdots < i_j \le n} 
e_{\tau_{i_0} \oplus \cdots \oplus \tau_{i_j}}
\Big(\varphi_{i_0}  \oplus 
\cdots \oplus 
\varphi_{i_j} (X) \Big) \ . 
$

\end{lem} 

Le multipoids de Chow g{\'e}n{\'e}ralise
donc    
le poids de Chow consid{\'e}r{\'e} au~\S~\ref{poids}. 
La partie (c) permet de ramener le calcul des 
multipoids de Chow  {\`a}
celui de poids de Chow standard.

\begin{demo}
La sym{\'e}trie est {\'e}vidente, puisque la forme 
$\Res_{\varphi_0, \dots, \varphi_n}$ est 
sym{\'e}trique ({\`a} un facteur scalaire pr{\`e}s) 
par rapport aux permutations 
des variables. 
Donc pour la partie (a), il suffit d'{\'e}tablir la lin{\'e}arit{\'e}
par rapport 
{\`a} la premi{\`e}re variable.

Soient 
 $\varphi: X \dashrightarrow \P^N$ et  
$\psi: X \dashrightarrow \P^P$ des plongements sur un ouvert dense
de $X$, et $\tau \in \R^{N+1}$, $\omega \in \R^{P+1}$ des poids. 
Soient $x=\{ x_0, \dots, x_N\} $ et $y = \{ y_0, \dots, y_P\}$ des
coordonn{\'e}es homog{\`e}nes pour $\P^N$ et $\P^P$ respectivement, 
et notons  $X_{\varphi, \psi}$ l'image de l'application 
$$
X \to \P^N \times \P^P \times \P^{N_1} \times \cdots \times \P^{N_n}
\enspace, \quad \quad s \mapsto \Big(\varphi(s), \psi(s), \varphi_1(s), 
\dots, \varphi_n(s)\Big) 
$$
et 
$I_{\varphi,\psi} \subset  \K[x,y,x_1, \dots, x_n]$
son id{\'e}al. 
Consid{\'e}rons maintenant
des groupes de variables $V= \{ V_0, \dots, V_N\}$, 
$W=\{ W_0, \dots, W_P\}$ et posons 
$V \cdot W:= 
\Big\{ V_i \, W_j \, : \ 0 \le i \le N, \ 0 \le j \le P \Big\}$.  
On a la factorisation~\cite[Prop.~3.5]{Rem01a} 
\begin{eqnarray} \label{factorisation} 
\res_{e+e', e_1, \dots, e_n}(I_{\varphi,\psi})
&& \hspace{-7mm} 
( V \cdot W,  U_1, \dots, U_n ) \nonumber \\[1mm]  
&& 
\hspace{-8mm}
= \lambda \cdot \res_{e, e_1, \dots, e_n}(I_{\varphi,\psi})
( V, U_1, \dots, U_n) 
\cdot  
\res_{e', e_1, \dots, e_n}(I_{\varphi,\psi})
( W, U_1, \dots, U_n) \ne 0  
\end{eqnarray}
o{\`u} $e,e', e_1, \dots, e_n$ d{\'e}signent les vecteurs de la 
base standard de $\R^{n+2}$
et $\lambda$ un {\'e}l{\'e}ment de $\K^\times$.  

\smallskip 

On v{\'e}rifie 
$\res_{e+e', e_1, \dots, e_n}(I_{\varphi,\psi})
= \Res_{\varphi \oplus \psi, \varphi_1, \dots, \varphi_n}$; 
cela vient de la th{\'e}orie des formes r{\'e}sul\-tan\-tes et 
peut se d{\'e}montrer avec des arguments tr{\`e}s proches 
de ceux de~\cite[pp.~102--103]{Rem01b}. 
Similairement 
$$
\res_{e, e_1, \dots, e_n}(I_{\varphi,\psi})
= \Res_{\varphi, \varphi_1, \dots, \varphi_n}
\enspace, \quad \quad 
\res_{e', e_1, \dots, e_n}(I_{\varphi,\psi})
= \Res_{\psi, \varphi_1, \dots, \varphi_n}
\enspace ;
$$ 
pour obtenir cela 
il faut en plus tenir compte du fait que 
les projections naturelles
$ \P^N  \times \P^P \times \prod_{i=1}^n \P^{N_i} 
\to  \P^N   \times \prod_{i=1}^n \P^{N_i}$
et $ \P^N  \times \P^P \times \prod_{i=1}^n \P^{N_i}
\to \P^P \times \prod_{i=1}^n \P^{N_i} $
sont des plongements sur un ouvert dense de $X_{\varphi, \psi}$. 

\smallskip 

Soit maintenant $U= \{ U_{i\,j} \, : \ 0 \le i \le N, \ 0 \le j \le P \} $ 
un groupe de $(N+1)(P+1)$ 
variables. 
On d{\'e}duit de ce qui pr{\'e}c{\`e}de et de la factorisation~(\ref{factorisation})
$$
\begin{array}{l} 
e_{\tau \oplus \omega, 
\tau_1,  
\dots, \tau_n}
(\varphi \oplus \psi, \varphi_1, \dots, \varphi_n;   X) 
=  
\deg_t 
\Big( \Res_{\varphi \oplus \psi, \varphi_1, \dots, \varphi_n} 
( t^{\tau \oplus \omega} \, U, 
t^{\tau_1}\, U_1, \dots,
t^{\tau_n} U_n ) \Big) \\[2mm] 
\hspace*{1mm} = 
\deg_t \Big ( 
\Res_{\varphi \oplus \psi, \varphi_1, \dots, \varphi_n} 
 ( (t^{\tau} \, V) \cdot (t^\omega \, W), 
t^{\tau_1}\, U_1, \dots,
t^{\tau_n} U_n ) \Big) 
\\[2mm] 
\hspace*{1mm} =  
\deg_t 
\Big( \Res_{\varphi, \varphi_1, \dots, \varphi_n} 
(t^{\tau}\,V , t^{\tau_1}\, U_1, \dots, t^{\tau_n} \,
U_n)  \Big)  
+\deg_t 
\Big(  
\Res_{\psi, \varphi_1, \dots, \varphi_n}
(t^{\omega}\, W , t^{\tau_1} \, U_1, 
\dots, t^{\tau_n} \, U_n)  \Big) \\[2mm] 
\hspace*{1mm} = 
e_{\tau, \tau_1,  
\dots, \tau_n}(\varphi, \varphi_1, \dots, \varphi_n;   X) 
+ 
e_{\omega, \tau_1,  
\dots, \tau_n}(\psi, \varphi_1, \dots, \varphi_n;   X) \enspace. 
\end{array}
$$

\smallskip

Passons {\`a} la partie (b), soit $I \subset \K[x_0, \dots, x_n]$ l'id{\'e}al de l'image du morphisme
$X \to (\P^N)^{n+1}$, $s \mapsto (\varphi(s), \dots, \varphi(s))$, 
on sait que 
$\res_{e_0, \dots, e_n}(I) = \elim_{e_0, \dots, e_n}(I)^f$  
pour une certaine puissance $f \ge 1$, 
o{\`u} 
$ \elim_{e_0, \dots, e_n}(I) $
d{\'e}signe la 
{\it forme {\'e}liminante} 
de $I$ d'indice $e_0, \dots, e_n$, {\it voir}~\cite[p.~74]{Rem01a}. 
Or,
par d{\'e}finition, cette forme {\'e}liminante co{\"\i}ncide 
avec la 
forme de Chow de $\ov{\varphi(X)} \subset \P^N$. 

En outre $\deg_{U_i}(\res_{e_0, \dots, e_n}(I)) = 
\deg(\varphi(X))$~\cite[Prop.~3.4]{Rem01a}, d'o{\`u} 
$\res_{e_0, \dots, e_n}(I)= \Ch_{\ov{\varphi(X)}}$. 
On en conclut 
$\displaystyle e_{\tau , \dots, \tau}(\varphi, \dots, \varphi;   X) 
= e_{\tau}(\ov{\varphi(X)}) $
car les d{\'e}finitions respectives co{\"\i}n\-ci\-dent. 
 
\smallskip

L'alin{\'e}a (c)
s'ensuit formellement 
de ces propri{\'e}t{\'e}s et de l'identit{\'e} alg{\'e}brique
\begin{equation} \label{polarisation} 
(n+1)! \, L(z_0, \dots, z_n) = 
\sum_{j=0}^n (-1)^{n+1 - j} 
\sum_{0 \le i_0 < \cdots < i_j \le n} 
\ell(z_{i_0}+ \cdots +z_{i_j})
\end{equation} 
valable pour 
une fonction multilin{\'e}aire sym{\'e}trique $L$ sur un semi-groupe 
ab{\'e}lien $(S, +)$ quelconque
et $\ell(z):=
L(z,\dots, z)$, et qu'on laisse au lecteur interess{\'e} le soin de d{\'e}montrer.  
\end{demo}


\typeout{Multihauteurs des plongements monomiaux}

\subsection{Multihauteurs des plongements monomiaux}

\label{demo thm2}

Soient $f: Q \to \R$ et $g: R \to \R$ des fonctions concaves
d{\'e}finies sur des ensembles convexes $Q, R \subset \R^n$ 
respectivement. 
On pose 
$$
f \boxplus g : Q+R \to \R
\enspace , \quad \quad
x \mapsto \max \{ f(y) + g(z) \, : \ y \in Q, \ z \in R, \ y+z =x\}
\enspace,  
$$
qui est une fonction concave d{\'e}finie sur la somme de Minkowski $Q+R$;  
on obtient ainsi une structure de semi-groupe commutatif sur
l'ensemble des  fonctions concaves.
Par la suite on supposera que toutes les fonctions concaves
consid{\'e}r{\'e}es sont 
d{\'e}finies sur des ensembles convexes {\it compacts}.

\begin{defn} \label{integrale mixte} 
Pour une famille de fonctions concaves 
$
f_0 : Q_0 \to \R, \dots,  f_n: Q_n\to \R$ l' {\em int{\'e}grale mixte} 
(ou {\em multi-int{\'e}grale}) est d{\'e}finie {\em via} la formule
\begin{equation*} 
\MI(f_0, \dots, f_n) := 
\sum_{j=0}^n (-1)^{n+1 - j} 
\sum_{0 \le i_0 < \cdots < i_j \le n} 
\int_{Q_{i_0} + \cdots + Q_{i_j}} f_{i_0} \boxplus 
\cdots \boxplus 
f_{i_j} \, dx_1 \cdots dx_n \enspace.
\end{equation*} 
\end{defn} 

Cette notion est proche de celle de {\it volume mixte} 
(ou {\it multi-volume}) $\MV(Q_1, \dots, Q_n)$ d'une famille
d'ensembles convexes $Q_1, \dots, Q_n \subset \R^n$, 
dont une d{\'e}finition possible est  
\begin{equation*} 
\MV(Q_1, \dots, Q_n) := 
\sum_{j=1}^n (-1)^{n - j} 
\sum_{1 \le i_1 < \cdots < i_j \le n} 
\Vol_n(Q_{i_1} + \cdots + Q_{i_j}) \enspace. 
\end{equation*} 
On sait que ceci g{\'e}n{\'e}ralise le volume d'un ensemble convexe, 
car on a 
$ \MV(Q, \dots, Q) = n! \, \Vol_n(Q)$.
Le volume mixte est sym{\'e}trique et lin{\'e}aire en chaque
variable $Q_i$ par rapport {\`a} la somme de Minkowski. 
On renvoie {\`a}~\cite[\S~7.4]{CLO98} pour les propri{\'e}t{\'e}s 
de base de cette notion. 
Souvent on notera $\MV_n(Q_1, \dots, Q_n)$
pour souligner qu'il s'agit du volume mixte $n$-dimensionnel.

\smallskip 

Les propri{\'e}t{\'e}s analogues de l'int{\'e}grale mixte sont 
r{\'e}sum{\'e}es dans la proposition suivante. 
Pour une fonction concave $f: Q \to \R$ et
$\mu \le
\min(f) $ on consid{\`e}re le polytope 
$$
Q_{f,\mu} := \Conv\Big( \Graph(f), Q \times \{\mu\}\Big) = 
\Conv\Big( (x, f(x)), (x, \mu) \, : \ x \in Q \Big) \ 
\subset \R^{n+1} \enspace;  
$$
on a $ \displaystyle 
\int_Q f \, dx_1 \cdots dx_n = \Vol_{n+1} (Q_{f,\mu}) + \mu \, \Vol_n
(Q)$ \enspace  .

\begin{prop} \label{MI}
\

\smallskip 

\noindent (a) \  
$ \MI(f_0, \dots, f_n)$ est sym{\'e}trique et
lin{\'e}aire en chaque variable $f_i$
par rapport {\`a} $\boxplus$
\enspace; 

\smallskip 

\noindent (b) \ 
$ \displaystyle \MI(f, \dots, f) = (n+1)! \, \int_Q f \, dx_1 \cdots
dx_n $
\enspace;  

\smallskip 

\noindent (c) \ 
$\MI(f_0, \dots, f_n) \ge 0 $ pour $f_0 \ge 0, \dots, f_n\ge 0$
\enspace; 

\medskip

\noindent (d) \ 
Soit $ \mu_i \le \min(f_i,0) $
  pour $i=0, \dots, n$, alors 
\vspace{-2mm} 
$$
\MI(f_0, \dots, f_n) = \MV_{n+1} (Q_{f_0, \mu_0}, \dots, Q_{f_n, \mu_n}) 
+ \sum_{i=0}^n \mu_i \, \MV_n(Q_0, \dots, Q_{i-1}, Q_{i+1}, \dots,
Q_n) \enspace.
$$
\end{prop} 

\vspace{-4mm} 

La notion d'int{\'e}grale mixte g{\'e}n{\'e}ralise donc celle d'int{\'e}grale d'une
fonction concave.  
La partie (d) montre qu'on peut en ramener le calcul {\`a} celui de
volumes 
mixtes (et donc de volumes standard). 

\begin{demo}

Pour la partie~{(a)}, la sym{\'e}trie de l'int{\'e}grale mixte 
est claire {\`a} partir de sa d{\'e}finition. 
En outre, la lin{\'e}arit{\'e} et la positivit{\'e} (alin{\'e}a~{(c)}) sont
imm{\'e}diatement 
v{\'e}rifi{\'e}es {\`a}
partir de~{(d)}, de l'identit{\'e}~(\ref{somme}) ci-dessous et
des propri{\'e}t{\'e}s analogues pour le volume mixte.  
Pour~{(b)}, on v{\'e}rifie  $\boxplus_{i=0}^j f \, ((j+1)
\, x) = 
(j+1) \, f(x) $ pour tout $x \in Q$ gr{\^a}ce {\`a} la concavit{\'e} de $f$; 
cet alin{\'e}a est donc une cons{\'e}quence de la formule (elle-m{\^e}me 
cons{\'e}quence de l'identit{\'e}~(\ref{polarisation})
appliqu{\'e}e {\`a} $(S, +):= (\R, +)$ et $L(z_0, \dots, z_n):= z_0 \cdots
z_n$) 
$$
(n+1)! = \sum_{j=0}^n (-1)^{n+1 - j} \,  {n+1
  \choose j+1} \, (j+1)^{n+1}
\enspace. 
$$

Pour la partie~{(d)}, soient $f : Q \to \R$ et $g : R \to \R$ des
fonctions concaves et $ \mu \le 
\min(f) $,  $\nu \le
\min(g) $.
On v{\'e}rifie alors $\mu + \nu \le \min( f\boxplus g)$ 
et 
\begin{equation} \label{somme} 
Q_{f\boxplus g, \mu + \nu} = Q_{f, \mu} + Q_{g, \nu} \enspace, 
\end{equation} 
\cad $f\boxplus g$ est la param{\'e}trisation de la toiture 
de la somme de Minkowski $Q_{f, \mu} + Q_{g, \nu}$. 
Soit maintenant  $\mu_i \le \min(f_i,0)$, on en d{\'e}duit 
\begin{eqnarray} \label{mui} 
\int_{Q_{i_0} + \cdots + Q_{i_j}} f_{i_0} \boxplus 
\cdots \boxplus 
f_{i_j} \, dx_1 \cdots dx_n
& = & 
\Vol_{n+1} \Big(Q_{f_{i_0} \boxplus 
\cdots \boxplus 
f_{i_j}, \, \mu_{i_0} + 
\cdots +
\mu_{i_j}}) \nonumber  \\[0mm] 
& & + (\mu_{i_0} + 
\cdots +
\mu_{i_j}) \, \Vol_n(Q_{i_0} + \cdots + Q_{i_j} )\\[2mm] 
& = & 
\Vol_{n+1} \Big( Q_{f_{i_0}, \mu_{i_0}} + \cdots 
+ Q_{f_{i_j}, \mu_{i_j}} \Big) \ 
\nonumber \\[0mm] 
&& - \Vol_{n+1}\Big(( Q_{i_0} + [0, -\mu_{i_0}]) +
\cdots +
(Q_{i_j} + [0,-\mu_{i_j}] ) \Big) \nonumber
\end{eqnarray}
o{\`u} l'on d{\'e}signe par 
$Q_{i_\ell}$ et $[0, -\mu_{i_\ell}]$  l'image de ces polytopes par les
inclusions
de $ \R^n$ et de $\R$ dans $\R^{n+1}$ {\it via} $(x_1, \dots,
x_n) \mapsto (x_1, \dots, x_n, 0)$ et $ x \mapsto (0, \dots, 0, x)$, 
respectivement. 
On d{\'e}duit de l{\`a} et des d{\'e}finitions de l'int{\'e}grale mixte et 
du volume mixte 
$$
\MI(f_0, \dots, f_n)  =  
\MV_{n+1} (Q_{f_0, \mu_0}, \dots, Q_{f_n, \mu_n}) 
- \MV_{n+1} \Big( Q_{0} + [0, -\mu_{0}], 
\dots, Q_{n} + [0,-\mu_n] \Big) 
\enspace. 
$$
Comme cons{\'e}quence de la formule~(\ref{mui}), on voit que ce dernier volume mixte
est une forme lin{\'e}aire en  $\mu_0, \dots, \mu_n$, \cad   
$$ 
\MV_{n+1} \Big( Q_{0} + [0, -\mu_{0}], 
\dots, Q_{n} + [0,-\mu_n] \Big) =  A_1 \, \mu_1 + \cdots + A_n
\, \mu_n 
\enspace, 
$$
et on a $ \displaystyle
A_i=  \MV_{n+1}( Q_{0} , 
\dots,  Q_{i-1}, [0,1], Q_{i+1}, \dots,  Q_{n} )$.
En applicant la d{\'e}finition du volume mixte on trouve
\begin{eqnarray*} 
A_i &=& \sum_{j=0}^n (-1)^{n+1 - j} 
\sum_{0 \le i_1 < \cdots < i_{j} \le n \atop i_\ell \ne i} 
\Vol_{n+1}([0,1] + Q_{i_1} + \cdots + Q_{i_{j}}) \\[-4mm] 
&& +  \ 
\sum_{j=0}^{n-1} (-1)^{n+1 - j} 
\sum_{0 \le i_0 < \cdots < i_j \le n \atop i_\ell \ne i } 
\Vol_{n+1}( Q_{i_0} + \cdots + Q_{i_j})  \\[-2mm] 
&= & 
\sum_{j=0}^{n-1} (-1)^{n+1 - j} 
\sum_{0 \le i_0 < \cdots < i_j \le n \atop i_\ell \ne i } 
\Vol_n( Q_{i_0} + \cdots + Q_{i_j}) 
\\[0mm] 
&= &  - \MV_{n} (Q_{0}, \dots, Q_{i-1}, Q_{i+1}, \dots, 
Q_{n})
\enspace, 
\end{eqnarray*} 
par 
le fait que $Q_0 + \dots + Q_n$ est de dimension $n$ et donc de
volume $n+1$-dimensionel nul, et 
$\Vol_{n+1}([0,1] + Q)= \Vol_n( Q)$ pour $Q \subset \R^n$.  
\end{demo} 

\medskip 

Revenons au 
cadre monomial: Soit 
$$
\cA_0 \in (\Z^n)^{N_0+1} \enspace, 
\dots \enspace, \enspace
\cA_n \in (\Z^n)^{N_n+1} 
$$
des vecteurs 
tels que $L_{\cA_i}= \Z^n$
pour  tout $ i$ et 
$\tau_0 \in \R^{N_0+1}, \dots, \tau_n \in \R^{N_n+1}$. 
Soit 
$\varphi_{\cA_i} : \T^n \to \P^{N_i}$ le plongement associ{\'e}
et
$\vartheta_{\cA_i, \tau_i} : Q_{\cA_i} \to \R$
la fonction 
affine par morceaux correspondante, pour $ i =0, \dots, n$. 
La proposition suivante explicite le multipoids de Chow 
du tore   $\T^n$ 
relatif {\`a} $(\varphi_{\cA_0}, \tau_0), $ $\dots, (\varphi_{\cA_n},
\tau_n)$:

\begin{prop} \label{e=MI}
$
\displaystyle
e_{\tau_0,  \dots , \tau_n}(\varphi_{\cA_0}, \dots,
\varphi_{\cA_n}; \T^n) 
= \MI(\vartheta_{\cA_0, \tau_0}, \dots, 
\vartheta_{\cA_n, \tau_n})$\enspace.
\end{prop} 

\begin{demo}

Soient $\cB= (b_0, \dots, b_N) \in (\Z^n)^{N+1}$, 
$\cC = (c_0, \dots, c_P) \in (\Z^n)^{P+1}$
des vecteurs et posons 
$$
\cB \oplus \cC := \Big( b_i + c_j \, : \ 0 \le i \le N, \ 0 \le j \le
P \Big) \ \in (\Z^n)^{(N+1)\, (P+1)} \enspace. 
$$
On v{\'e}rifie  $\varphi_\cB \oplus \varphi_\cC = 
\varphi_{\cB \oplus \cC} $. 
Similairement, pour des poids 
$\tau \in \R^{N+1}$ et $\omega \in \R^{P+1}$
on v{\'e}rifie
$\displaystyle \vartheta_{\cB, \tau} \boxplus \vartheta_{\cC, \omega} 
= \vartheta_{ \cB \oplus \cC, \tau \oplus \omega} $. 
Le r{\'e}sultat est donc une cons{\'e}quence directe 
du lemme~\ref{proprie-multipoids}(c), 
du cas non mixte 
(Proposition~\ref{poid-torique}) et de la d{\'e}finition de l'int{\'e}grale
  mixte.
\end{demo}

Soit maintenant  
$
\alpha_0 \in (K^\times)^{N_0+1}, 
\dots , \alpha_n \in (K^\times)^{N_n+1}
$
et 
notons 
$\varphi_{i}:  \T^n \to \P^{N_i}$ 
le plongement monomial associ{\'e} au  
couple  $(\cA_i , \alpha_i)$, 
pour $i=0, \dots, n$.   
On d{\'e}montre enfin le th{\'e}or{\`e}me~\ref{thm2} explicitant la 
multihauteur normalis{\'e}e
$$
\wh{h}((\cA_0, \alpha_0), \dots, (\cA_n, \alpha_n); \T^n)
:= \hnorm(\varphi_0, \dots, \varphi_n ; \T^n) 
$$
du tore $\T^n$ relative {\`a} ces plongements, en termes d'int{\'e}grales
mixtes. 
Pour chaque $v \in M_K$ on note $\vartheta_{i, v}:
 Q_{\cA_i} \to \R$  la fonction param{\'e}trant la toiture 
du polytope $Q_{i,v} \subset \R^{n+1}$ 
associ{\'e} au vecteur $\cA_i$ et au poids 
$\tau_{\alpha_i \, v}$.

\begin{demo}[D{\'e}monstration du th{\'e}or{\`e}me~\ref{thm2}]
Soit 
$\beta \in (K^\times)^{N+1}$ et $\gamma \in (K^\times)^{P+1}$
et posons 
$$
\beta \otimes \gamma  := \Big( \beta_i \, \gamma_j \, 
: \ 0 \le i \le N, \ 0 \le j \le
P \Big) \ \in (K^\times)^{(N+1) \, (P+1) } \enspace. 
$$
On v{\'e}rifie sans peine 
  $\varphi_{\cB, \beta}  
\oplus \varphi_{\cC, \gamma} = 
\varphi_{\cB \oplus \cC, \beta \otimes \gamma} $
pour des vecteurs  
$\cB \in (\Z^n)^{N+1}$ 
et $\cC \in (\Z^n)^{P+1}$.  
Posons $\cA:= (\cA_0, \dots, \cA_n)$ et
$\alpha:= (\alpha_0, \dots, \alpha_n)$, 
on consid{\`e}re alors  l'application 
$$
\Phi_{\cA, \alpha}: \T^n \to \P^{N_0} \times \cdots \times \P^{N_n} 
\enspace , \quad \quad 
s \mapsto (\varphi_0(s), \dots, \varphi_{n}(s)) 
$$
et on note $Z:= \ov{\Phi_{\cA, \alpha}(\T^n)}$ 
l'adh{\'e}rence de Zariski de son image. 
Pour  $D=(D_0, \dots, D_n) \in (\N^*)^{n+1} $ on rappelle le
plongement mixte 
$\Psi_D :  \P^{N_0} \times \cdots \times \P^{N_m} \to  
\P^{ {D_0+N_0 \choose N_0} \cdots {D_m+N_m \choose N_m} -1}$
introduit au paragraphe~\ref{hauteur normalisee}; 
on v{\'e}rifie 
$
\Psi_D \circ \Phi_{\cA, \alpha}  = 
D_0 \, \varphi_0 \oplus \cdots \oplus D_n \, 
\varphi_{n}=
\varphi_{D \cdot\cA, \alpha^{\otimes D}} $
avec 
$
D \cdot \cA:= D_0 \, \cA_0 \oplus \cdots \oplus D_n \, \cA_n $ et 
$\alpha^{\otimes D} := \alpha_0^{\otimes D_0} \otimes \cdots 
\otimes \alpha_n^{\otimes
  D_n} $.
Donc  
$\Psi_D(Z)  = X_{D \cdot \cA, \alpha^{\otimes D} } $ 
et en appliquant l'identit{\'e}~(\ref{decomp hnorm}) 
on trouve 
\begin{equation} \label{hnorm DA} 
\hnorm(X_{D \cdot\cA, \alpha^{\otimes D} }) = 
\hnorm(\Psi_D(Z)) =
\sum_{c \in \N_{n+1}^{n+1}}
{n+1 \choose c} \, \hnorm_c(Z) \, D^{c} \enspace. 
\end{equation}

Soit $v \in M_K$, 
la vari{\'e}t{\'e} torique 
 $X_{D \cdot \cA} \subset 
\P^{ {D_0+N_0 \choose N_0} \cdots {D_m+N_m \choose N_m} -1}$
associ{\'e}e au vecteur $D \cdot \cA$
satisfait 
$$
e_{\tau_{\alpha^{\otimes D} \,v}} (X_{D\cdot \cA}) = 
e_{D \cdot \tau_{\alpha v}, \dots, 
D \cdot \tau_{\alpha v}} (D \cdot \varphi_\cA, 
\dots, D \cdot \varphi_{\cA}; \T^n)
$$
avec $
D \cdot \varphi_\cA:= D_0 \, \varphi_{\cA_0} \oplus \cdots \oplus 
D_n \, \varphi_{\cA_n}$ et $D\cdot \tau_{\alpha v} := 
D_0 \, \tau_{\alpha_0\, v} \oplus \cdots \oplus D_n \, 
\tau_{\alpha_n\, v}= \tau_{\alpha^{\otimes D} \,v}$. 
La sym{\'e}trie et la multilin{\'e}arit{\'e} 
des multipoids de Chow impliquent 
\begin{equation} \label{e DA} 
e_{D \cdot\tau_{\alpha v}} (X_{D\cdot \cA}) = 
\sum_{c \in \N_{n+1}^{n+1}}
{n+1 \choose c} \, 
e_c (\cA, \alpha, v) 
\, D^{c} 
\end{equation} 
o{\`u} 
$e_c (\cA, \alpha, v) $
d{\'e}signe le multipoids de Chow 
$e_{\boldsymbol{\tau}}(\boldsymbol{\varphi}; \T^n)$
du tore $\T^n$ relatif {\`a}
$$
(\boldsymbol{\varphi}, \boldsymbol{\tau}):= 
\overbrace{\Big(\varphi_{\cA_0},  \tau_{\alpha_0\, v}\Big), \dots,  
\Big(\varphi_{\cA_0},  \tau_{\alpha_0\, v}\Big)}^{c_0 \mbox{ \scriptsize fois}}
\enspace ,
\dots
, \enspace
\overbrace{\Big(\varphi_{\cA_n},  \tau_{\alpha_n\, v}\Big), \dots,  
\Big(\varphi_{\cA_n},  \tau_{\alpha_n\, v}\Big)}^{c_n 
\mbox{ \scriptsize fois}} \enspace. 
$$
Comme cons{\'e}quence des d{\'e}compositions~(\ref{hnorm DA}),~(\ref{e DA})
et 
du cas non
mixte (th{\'e}or{\`e}me~\ref{thm1}) on d{\'e}duit 
$$
\sum_{c \in \N_{n+1}^{n+1}}
{n+1 \choose c} \, \hnorm_c(Z) \, D^{c} 
=
\sum_{c \in \N_{n+1}^{n+1}}
{n+1 \choose c} \,
\Bigg(
\sum_{v \in M_K} 
\frac{[K_v:\Q_v]}{[K:\Q]} 
e_c(\cA, \alpha, v) 
\Bigg) 
\, D^{c} \enspace. 
$$
Cette identit{\'e} polynomiale {\'e}tant valable pour {\it tout}
$D \in (\N^*)^{N+1}$, on en
d{\'e}duit que les coe\-ffi\-cients respectifs doivent co{\"\i}ncider. 
La multihauteur  $ \hnorm((\cA_0, \alpha_0), \dots, (\cA_n,
  \alpha_n) ; \T^n)$
est par d{\'e}finition le coefficient de $(n+1)! \, D_0 \cdots D_n$ 
(correspondant {\`a} $c=(1, \dots, 1) \in \N^{n+1}$) 
dans
cette expression, donc en particulier 
$$
\hnorm((\cA_0, \alpha_0), \dots, (\cA_n,
  \alpha_n) ; \T^n) 
= \sum_{v \in M_K} 
\frac{[K_v:\Q_v]}{[K:\Q]} \,
e_{\tau_{\alpha_0 \, v}, \dots, \tau_{\alpha_n \, v}} (\varphi_{\cA_0},
\dots, \varphi_{\cA_n}; \T^n) 
\enspace,
$$
d'o{\`u} on conclut par application de la proposition pr{\'e}c{\'e}dente. 
\end{demo} 

\begin{rem} 
Plus g{\'e}n{\'e}ralement, on a d{\'e}montr{\'e}
pour tout $c \in \N^{n+1}_{n+1}$:
$$
\hnorm_c( \Phi_{\cA, \alpha}(\T^n)) 
= \sum_{v \in M_K} 
\frac{[K_v:\Q_v]}{[K:\Q]} \,
e_c(\cA,\alpha, v)
= \sum_{v \in M_K} 
\frac{[K_v:\Q_v]}{[K:\Q]} \,
\MI_c( \vartheta_{\cA, \tau_{\alpha\, v}}) 
$$
avec $\MI_c( \vartheta_{\cA, \tau_{\alpha v}})
:= 
\MI\Big( \overbrace{\vartheta_{\cA_0, \tau_{\alpha_0\, v}}, \dots,  
\vartheta_{\cA_0, \tau_{\alpha_0\, v}}}^{c_0 \mbox{ \scriptsize fois}}
\enspace, 
\dots
\enspace , \enspace  
\overbrace{\vartheta_{\cA_n, \tau_{\alpha_n\, v}}, \dots,  
\vartheta_{\cA_n, \tau_{\alpha_n\, v}}}^{c_n \mbox{ \scriptsize fois}} \Big)$.
\end{rem}

\medskip 

On illustre ce r{\'e}sultat sur un exemple. 
Consid{\'e}rons des plongements monomiaux $\varphi: \T^1 \to \P^1$
, $ \displaystyle s \mapsto \Big( \frac{1}{2} : 4\, s\Big)$ et 
$\psi: \T^1 \to \P^1$
, $ \displaystyle s \mapsto \Big( \frac{1}{3} :  \frac{s}{2}\Big)$. 
La figure suivante montre
pour chaque $v \in M_\Q$, les graphes des fonctions $\vartheta_v$
associ{\'e}es, 
en trait discontinu pour celles correspondant {\`a} $\varphi$ 
et $\psi$ et en gras pour 
celle correspondant {\`a} la somme $\varphi \oplus \psi$: 

\vfill 
\pagebreak


\vspace*{56mm} 

\begin{figure}[htbp]

\begin{picture}(0,0) 

\put(40,50){
\epsfig{file=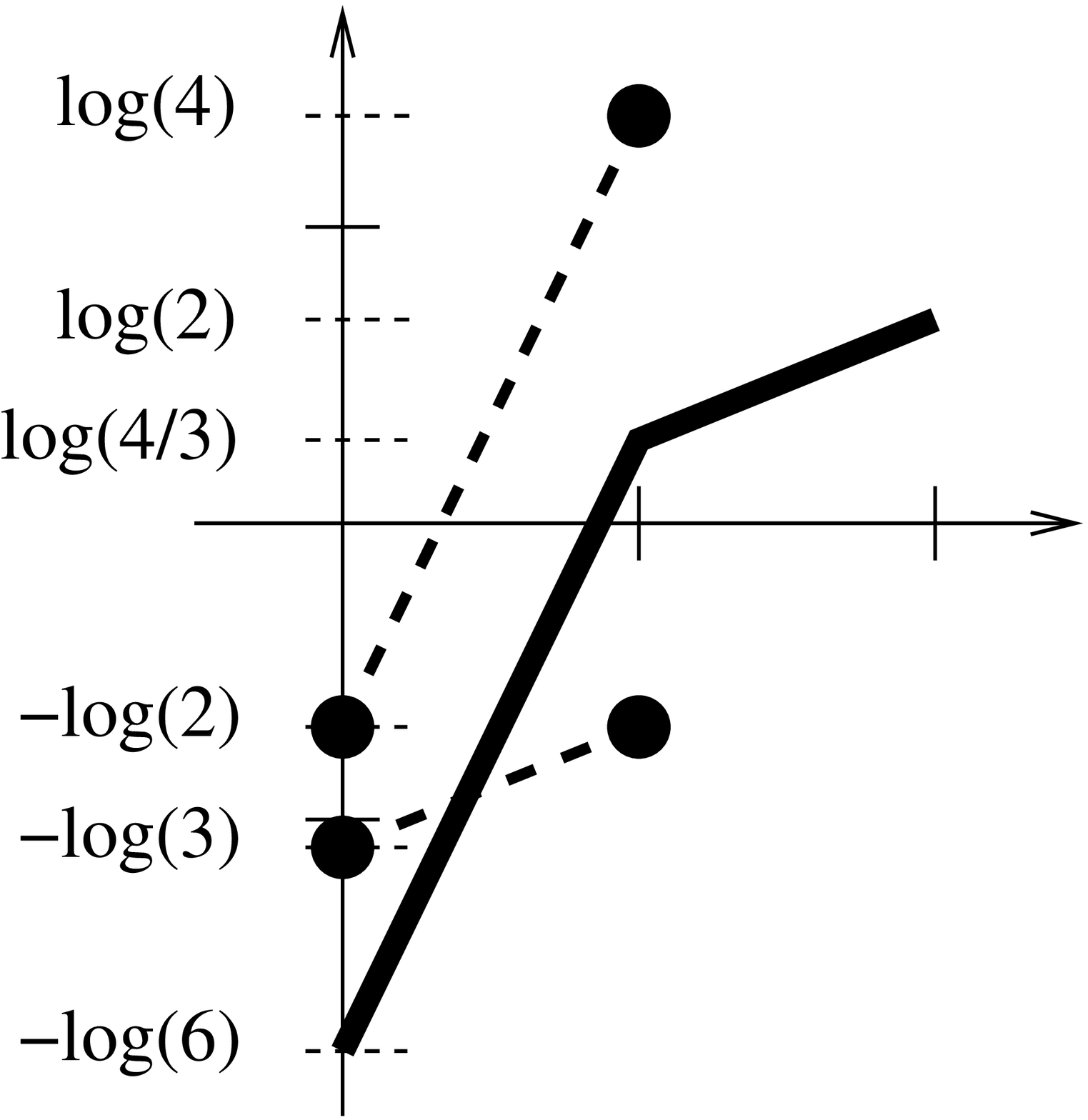, height= 46 mm} 
} 
\put(160,98){$v=\infty$} 

\put(275,66){
\epsfig{file=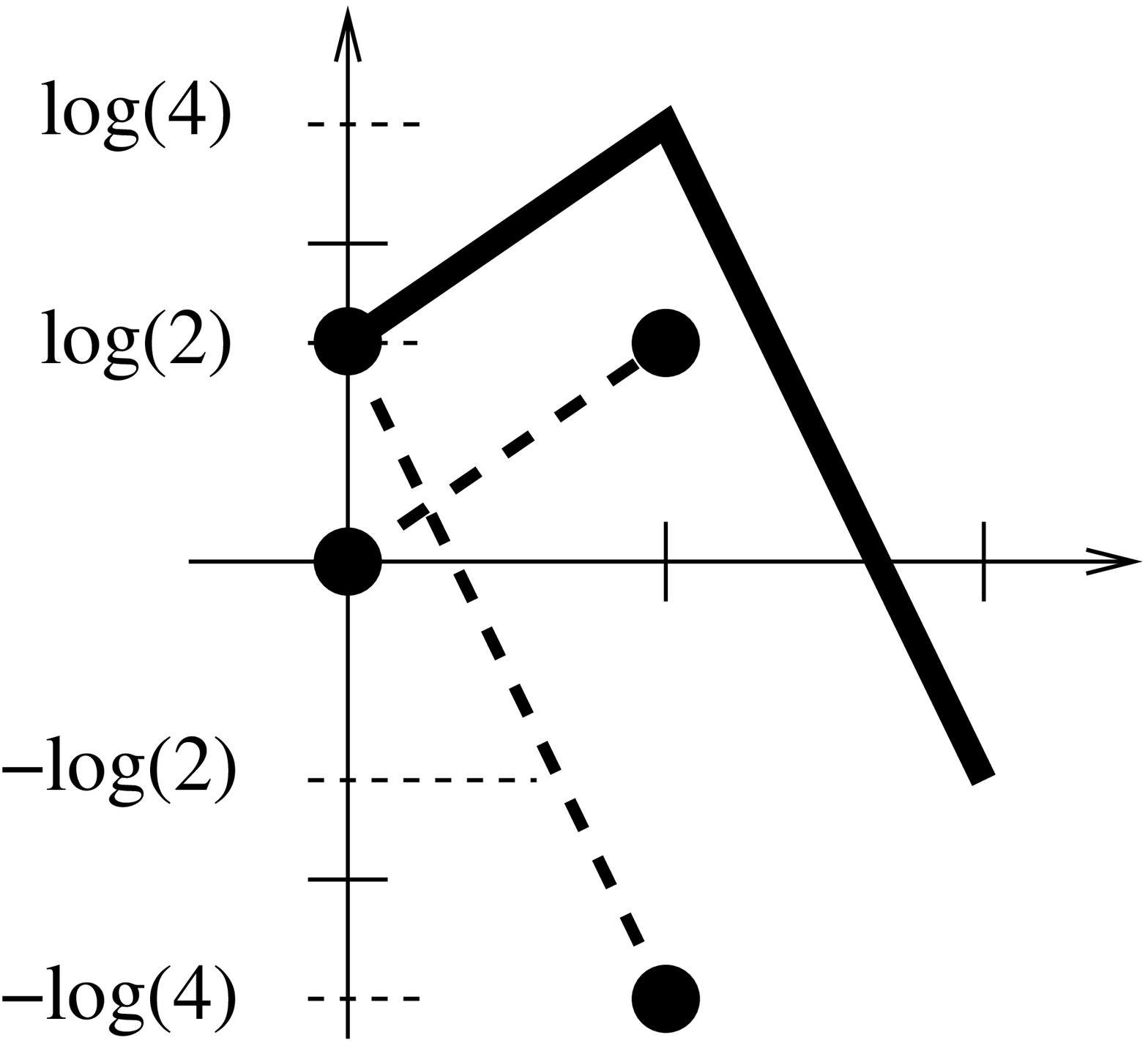, height= 40 mm} 
} 
\put(405,98){$v=2$} 

\put(143,0){
\epsfig{file=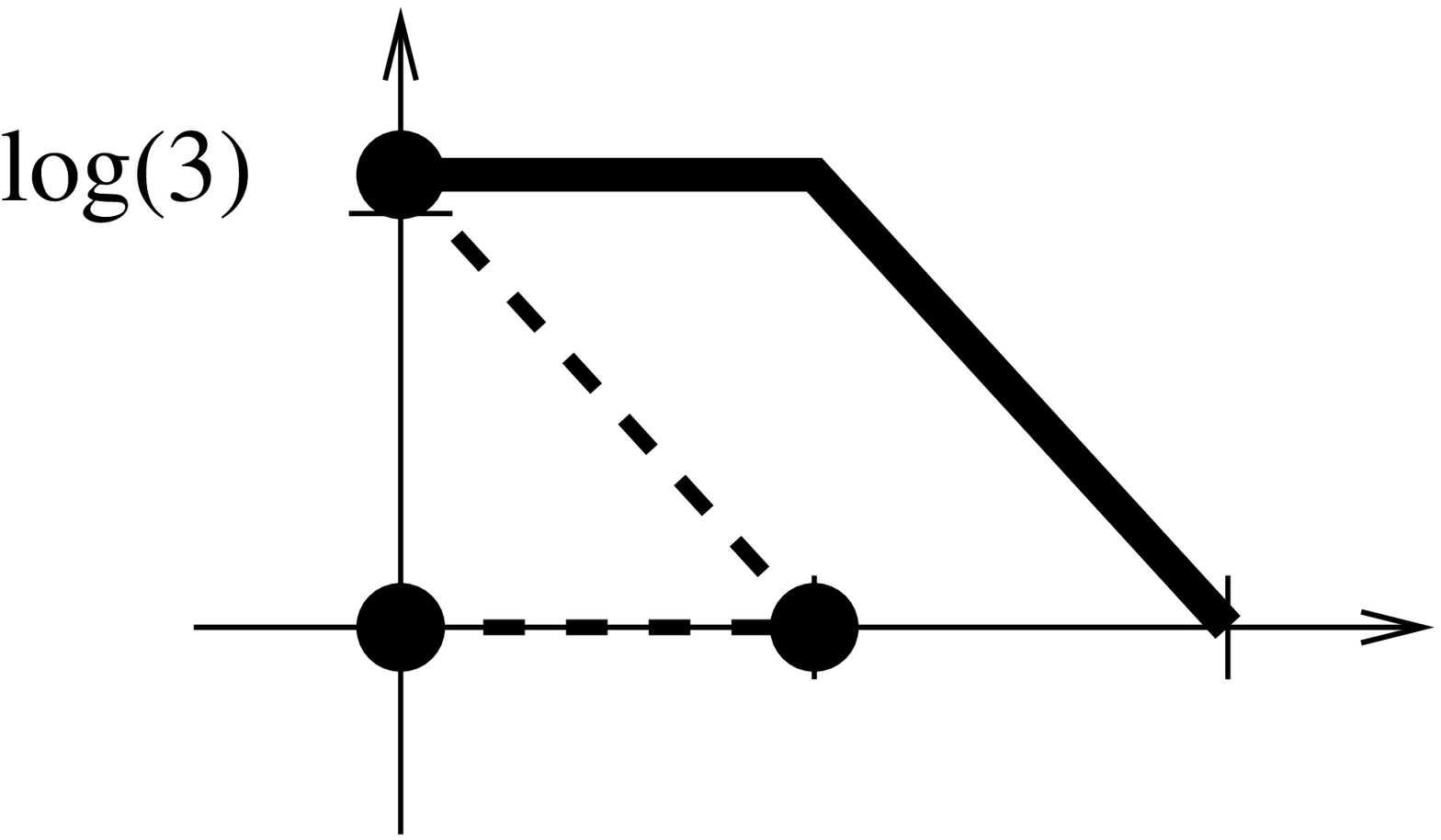, height= 25 mm} 
} 
\put(255,37){$v= 3$}

\end{picture} 

\end{figure}

\vspace{-5mm} 


\noindent et $\vartheta_v \equiv 0$ pour $v \ne \infty, 2,3$. 
D'o{\`u} 
\begin{eqnarray*} 
\hnorm(\varphi, \psi ;\T^1) &=& 
\MI(\vartheta_{\varphi, \infty}, \vartheta_{\psi, \infty}) 
+ \MI(\vartheta_{\varphi, 2}, \vartheta_{\psi, 2}) 
+ \MI(\vartheta_{\varphi, 3}, \vartheta_{\psi, 3}) 
\\[2mm] 
& = &  
\sum_{v=\infty, 2,3} 
\bigg( \int_0^2 \, \vartheta_{\varphi \oplus \psi, v}  \, dx  \ - 
\int_0^2 \, \vartheta_{\varphi, v}  \, dx  \ -  
\int_0^2 \, \vartheta_{\psi, v}  \, dx \bigg) \\[2mm] 
& =&  \Big( 2\, \log(2) -\log(3) ) + 2\, \log(2) + \log(3) 
= 4 \, \log(2)  \enspace.
\end{eqnarray*} 
Remarquons $ \ov{\varphi(\T^1)} = \ov{\psi(\T^1)} = \P^1$,  
donc $ \hnorm(\varphi, \varphi ;\T^1) = \hnorm(\psi, \psi ; \T^1)
  =0$\enspace.


\typeout{Joints, plongements de Segre et de Veronese} 

\subsection{Joints, 
plongements de Segre et de Veronese}

\label{joints}

La classe des vari{\'e}t{\'e}s toriques est ferm{\'e}e par
formation de 
joints, produits de Segre et plongements de Veronese. 
De plus, 
ces constructions se transportent de mani{\`e}re naturelle en des
constructions 
standard sur les polytopes associ{\'e}s. 
Dans la suite, on montre 
que la famille de fonctions $\Theta_\arith$ 
associ{\'e}e {\`a} une vari{\'e}t{\'e}
torique 
se comporte de mani{\`e}re aussi naturelle par rapport {\`a} 
ces op{\'e}rations.

\medskip 

Soient $X \subset \P^N$ et $Y\subset \P^P$ des vari{\'e}t{\'e}s
de dimension $n$ et $p$ respectivement.  
Le {\em joint} 
est par d{\'e}finition
$$
X \# Y = \Big\{ (v\, x : w \,y  ) \, : \ 
x \in X, \ y \in Y, \ (v:w) \in \P^1 \Big\} 
\ \subset \P^{N+P+1}
\enspace.  
$$
Soient $i(X)$ et $ j(Y)$ respectivement, 
les plongements de $X$ et de $Y$ 
{\it via} les inclusions standard 
$i : \P^N \hookrightarrow \P^{N+P+1}$ et 
$j : \P^P \hookrightarrow \P^{N+P+1}$.
Le joint est donc la r{\'e}union des droites joignant les points de $i(X)$ {\`a}
ceux de $j(Y)$; 
c'est une vari{\'e}t{\'e} de 
dimension $n+p+1$ et degr{\'e} $\deg(X) \, \deg(Y)$. 

\smallskip 

Soient
$
\cA= 
 \Big( a_0, \dots, a_N\Big) 
\in ( \Z^n)^{N+1}$,  
$\alpha=(\alpha_0, \dots, \alpha_N)  \in (K^\times)^{N+1}$, 
$\cB=\Big( b_0, \dots, b_P\Big) \in (\Z^p)^{P+1}$ et 
$\beta=(\beta_0, \dots, \beta_P) \in (K^\times)^{P+1}$ 
tels que $L_\cA = \Z^n$ et $L_\cB = \Z^p$,  
on pose alors 
\begin{eqnarray*} 
\cA \# \cB &:=& 
\Big(  (1, a_0, \mathbf{0}_p), \dots,  (1, a_N, \mathbf{0}_p), 
(0,\mathbf{0}_n,b_0), \dots, (0,\mathbf{0}_n,b_P) 
\Big) \in(\Z^{n+p+1} )^{N+P+2}
\enspace, \\[1mm] 
\alpha \# \beta &:=& (\alpha_0, \dots , \alpha_N, 
\beta_0, \dots, \beta_P) \in (K^\times)^{N+P+2}
\enspace. 
\end{eqnarray*}
On v{\'e}rifie 
$
X_\arith^\circ \# X_{\cB, \beta}^\circ =
\Big\{ (u\, \varphi_{\cA, \alpha}(s) : \varphi_{\cB, \beta} (t)) \, : 
\ s \in \T^n, \ t \in \T^p, 
\ 
u\in \T^1
\Big\} 
= X_{\cA \# \cB, \alpha \#\beta}^\circ 
$
donc $$
X_\arith \# X_{\cB, \beta}= X_{\cA \# \cB, \alpha \#\beta} 
\enspace .
$$ 
Posons
$Q\subset
\R^n$ et  $R \subset
\R^p$ les polytopes associ{\'e}s {\`a} $\cA$ et $\cB$ respectivement. 
Le polytope $ Q_{\cA \# \cB} \subset \R^{n+p+1} $ 
s'{\'e}crit 
comme  l'image de l'application 
$$
j_{Q, R} : [0,1] \times Q \times R  \to \R^{n+p+1} 
\enspace , \quad \quad  
(u, s,t) \mapsto (u, u \, s, (1-u) \, t)
\enspace;  
$$
autrement-dit $Q_{\cA \# \cB} $ est la r{\'e}union des 
segments {joignant} les points de 
$\{ 1\} \times Q \times \{ \mathbf{0}_p\} \cong Q$ {\`a} ceux 
de
$\{ 0\} \times \{ \mathbf{0}_n \} \times R \cong R$. 

\smallskip 

Soit maintenant
$$
s_{N,P} : \P^N\times \P^P \to \P^{(N+1)(P+1)- 1}
\enspace ,\quad \quad 
(x,y) \mapsto \Big( x_i\, y_j \, : \ i=0, \dots, N, \, j = 0, \dots, P\Big)
$$
le  {plongement de Segre}. 
Le {\it produit de Segre} $X\times Y \subset \P^{(N+1)(P+1)- 1}$ est
d{\'e}fini comme
l'image par $s_{N,P}$ 
du produit cart{\'e}sien de $X$ et $Y$; 
c'est une vari{\'e}t{\'e} de dimension $n+p$ et degr{\'e} ${{n+p}\choose
  n}\, \deg(X)\, \deg(Y)$. 
On pose alors
\begin{eqnarray*}
\cA \times \cB &:=&  \Big( (a_i, b_j)\, : \ 
i=0, \dots, N, \, j =
0, \dots, P \Big) \in   (\Z^{n+p})^{(N+1)(P+1)}
\enspace, \\[1mm] 
\alpha \otimes \beta &:=& \Big( \alpha_i\, \beta_j \, :  \ 
i=0, \dots, N, \, j =
0, \dots, P \Big) \in   (K^\times)^{(N+1)(P+1)}
\enspace,
\end{eqnarray*}
notons que cette derni{\`e}re notation a d{\'e}j{\`a} {\'e}t{\'e} introduite au cours de la preuve 
du th{\'e}or{\`e}me~\ref{thm2} dans la section~\ref{normalisee}.
Une v{\'e}rification directe montre  
$
X_\arith \times  X_{\cB, \beta} = s_{N,P} (X_\arith, X_{\cB, \beta} ) 
= X_{\cA \times \cB, \alpha \otimes \beta} $ et que
le polytope associ{\'e} $Q_{\cA \times \cB}\subset \R^{n+p} $ 
est 
le produit cart{\'e}sien
$Q \times R$. 

\smallskip 

Pour $D \in \N^*$ on pose 
$$
v_{N,D} : \P^N \to \P^{{D+N \choose N}-1} 
\enspace, \quad  \quad 
x \mapsto \Big( x^b \, : \ b \in \N^{N+1}_D \Big)
$$
le {\em plongement de Veronese}  de degr{\'e} $D$ en $N+1$ variables
homog{\`e}nes. 
L'image 
$v_{N,D} (X_\arith) \subset \P^{{D+N \choose N}-1}$
est une vari{\'e}t{\'e} de dimension $n$ et degr{\'e} $D^n \deg(X)$. 
Posons 
\begin{eqnarray*} 
V(\cA)  & := &  \Big( b_0 \, a_0 + \cdots +b_N \, a_N 
\, : \ b \in \N^{N+1}_D \Big)
\in \Big( \Z^n)^{N+D \choose N}
\enspace,  \\[1mm] 
V(\alpha)  & := & 
\Big( 
\alpha_0^{b_0} \cdots \alpha_N^{b_N} \, : \ b \in \N^{N+1}_D \Big)
\in \Big( K^\times)^{N+D \choose N}
\enspace.
\end{eqnarray*} 
On v{\'e}rifie ais{\'e}ment 
$ v_{N,D} (X_\arith) = X_{V(\cA), V(\alpha)}$  
et que le polytope associ{\'e} $Q_{V(\cA)} \subset \R^n$ 
est l'homoth{\'e}tique
$
D\, Q$. 

\smallskip 

Nous collectons dans l'{\'e}nonc{\'e} suivant le comportement des fonctions
$\vartheta_v$ par rapport {\`a} ces op{\'e}rations. 
La d{\'e}monstration suit des observations pr{\'e}c{\'e}dentes, 
nous laissons les d{\'e}tails au lecteur int{\'e}ress{\'e}.

\begin{prop} \label{theta(joint-segre-veronese)}
Avec les notations introduites, pour $v \in M_K$ on a 

\smallskip

\noindent (a) \ 
$
\vartheta_{\cA \# \cB, \tau_{\alpha \# \beta \, v}} 
\circ j_{Q, R} (u, s,t) = 
u \,
\vartheta_{\cA, \tau_{\alpha  v}} (s) 
+ (1-u) \, \vartheta_{\cB, \tau_{\beta  v}} (t) 
$\enspace
pour $u\in [0,1]$, $s \in Q$, $t\in R$\enspace; 

\medskip 
 
\noindent (b) \  
$
\vartheta_{\cA \times \cB, \tau_{\alpha \times \beta \, v}}  (s,t) = 
\vartheta_{\cA , \tau_{\alpha v}} (s) 
+ \vartheta_{\cB, \tau_{\beta  v}} (t) 
$\enspace
pour $s \in Q$, $t \in R$\enspace; 

\medskip 

\noindent (c) \  
$ 
\vartheta_{ V(\cA), \tau_{V(\alpha) \, v}} (D\,s) = 
D \, \vartheta_{ \cA, \tau_{\alpha  v}} (s)  $ \enspace pour $s \in Q$
\enspace.

\end{prop}

\smallskip 

En calculant l'int{\'e}grale de ces fonctions 
on obtient, {\`a} partir de
cette proposition et du th{\'e}or{\`e}me~\ref{thm1}, des formules pour la 
hauteur normalis{\'e}e des joints, produits de Segre
et plongements de Veronese des vari{\'e}t{\'e}s toriques. 
En fait, ces formules  
sont valables pour des vari{\'e}t{\'e}s quelconques
et sont 
cons{\'e}quence
de r{\'e}sultats analogues pour la hauteur projective:

\begin{prop} \label{joint-segre-veronese} 
Soient $X \subset \P^N$ et $Y\subset \P^P$ des sous-vari{\'e}t{\'e}s
de dimension $n$ et $p$ respectivement, 
alors 

\smallskip 

\noindent (a) \ 
$\hnorm(X \# Y)  =  \hnorm(X) \, \deg(Y) + 
\deg(X)  \, \hnorm(Y)$
\enspace;  

\medskip 
 
\noindent (b) \ 
$ \hnorm(X \times Y)  =   {{n+p+1}\choose p} \, 
\hnorm(X) \, \deg(Y) 
+ {{n+p+1}\choose n} \,\deg(X) \, \hnorm(Y)$
\enspace; 

\smallskip 

\noindent (c) \ 
$\hnorm(v_{N,D}(X)) = D^{n+1} \, \hnorm(X)$ pour $D
  \in \N^*$
\enspace.   

\end{prop} 

\begin{demo} 
Pour la premi{\`e}re partie on applique
l'{\'e}galit{\'e} 
$[k] \, (X \# Y) = ([k] \, X ) \# ([k] \, Y) $ qu'on v{\'e}rifie 
de fa{\c c}on directe {\`a} partir de la d{\'e}finition du joint. 
D'apr{\`e}s~\cite[Prop. 2]{Phi95}  
$$
\frac{h\Big([k] \, (X \# Y) \Big)}{\deg \Big([k] \, (X \# Y) \Big)} 
= \frac{h([k] \, X)}{\deg ([k] \, X )}
+  \frac{h([k] \, Y)}{\deg ([k] \, Y  )}
+ c(n,p) 
$$
avec $
\displaystyle c(n,p) := \sum_{i=0}^n \sum_{j=0}^p \frac{1}{2 \,
  (i+j+1)}$.  
En prenant la limite on en d{\'e}duit   
$$
\deg (X \#  Y) \, \hnorm(X \# Y)
=
\deg (X \#  Y) \cdot \lim_{k \to \infty} 
\frac{h\Big([k] \, (X \# Y) \Big)}{k\, \deg \Big([k] \, (X \# Y) \Big)} 
= \hnorm(X) \, \deg(Y) + 
\deg(X)  \, \hnorm(Y)  
$$
car $\deg(X \# Y) = \deg(X) \, \deg(Y)$. 
La partie~(b) se d{\'e}montre de la m{\^e}me fa{\c c}on, 
en utilisant 
$[k] \, (X \times  Y) 
= ([k] \, X) \times  ([k] \, Y)$ et~\cite[Prop. 1]{Phi95}. 
La partie~(c) correspond {\`a} l'identit{\'e}~(\ref{decomp hnorm})
appliqu{\'e}e {\`a} $m:=0$ et $D_0:=D$. 
\end{demo}

En g{\'e}n{\'e}ral, il est int{\'e}ressant de savoir calculer la hauteur des
vari{\'e}t{\'e}s projectives produites par diverses constructions. 
Comme g{\'e}n{\'e}ralisation des plongements de Veronese, 
on peut consid{\'e}rer 
les applications monomiales $\varphi_{\cB, \beta} : \P^N \to \P^M$, 
$ x \mapsto (\beta_0 \, x^{b_0} : \cdots : \beta_M \, x^{b_M} ) $
avec $b_j \in \N^{N+1}_D$  et
$\beta_j \in K^\times$
pour un certain $D \in \N^*$. 

L'image
d'une vari{\'e}t{\'e} torique 
par une telle application est aussi torique: 
soit $\cA \in (\Z^n)^{N+1}$ et $\alpha \in \P^N$, alors 
$$
\ov{\varphi_{\cB, \beta}(X_\arith)} = X_{\cC, \gamma}
$$ 
o{\`u} $\cC=( c_0, \dots, c_M) \in (\Z^n)^{M+1}$ avec 
$c_j:= M_\cA(b_j) = b_{j \, 0} \, a_0 + \cdots + b_{j \, N} \, a_N$, et
$\gamma = \Phi_{\cB, \beta}(\alpha) = ( \beta_0 \, \alpha^{b_0}:
\cdots: \beta_M \, \alpha^{b_M}) \in \P^M$. 
Remarquons que lorsque la classe des vari{\'e}t{\'e}s toriques est ferm{\'e}e par ces
constructions, 
elles sont de bons candidats pour 
tester des formules conjecturales.





%
%


\typeout{References}

\vfill

\pagebreak


\end{document}